\documentclass[11pt]{article}

\usepackage{psfrag,cite,amsmath,graphicx,epsfig,latexsym,subfigure,color}
\usepackage[scriptsize,bf]{caption}
\usepackage{color,xcolor,multirow}
\usepackage{amssymb}
\usepackage{bbm}
\usepackage{setspace}
\usepackage{array}
\usepackage{booktabs}
\usepackage{hyperref}

\DeclareMathOperator{\argmin}{\mbox{argmin}}
\def\bn{\hfill \\ \smallskip\noindent}
\def\argmin{\mathop{\rm argmin}}

\newcommand{\beq}{\begin{equation}}
\newcommand{\eeq}{\end{equation}}
\newcommand{\st}{{\rm s.t.}}

\newcommand{\cL}{\mathcal{L}}
\newcommand{\hL}{\widehat{\mathcal{L}}}
\newcommand{\ctL}{\widetilde{\mathcal{L}}}
\newcommand{\cP}{\mathcal{P}}
\newcommand{\barh}{\bar{h}}
\newcommand{\barf}{\bar{f}}

\newcommand{\cO}{{\mbox{$\mathcal{O}$}}}

\newcommand{\cG}{{\mbox{$\mathcal{G}$}}}
\newcommand{\cA}{{\mbox{$\mathcal{A}$}}}
\newcommand{\cV}{{\mbox{$\mathcal{V}$}}}

\newcommand{\tL}{{\mbox{$\widetilde{L}$}}}

\newcommand{\tcL}{{\mbox{$\widetilde{\mathcal{L}}$}}}

\newcommand{\norm}[1]{\left\|{#1}\right\|} 

\newcommand{\floor}[1]{\left\lfloor #1 \right\rfloor}

\begin{document}
\def\pn {\par\smallskip\noindent}
\def \bn {\hfill \\ \smallskip\noindent}
\newcommand{\fs}{f_1,\ldots,f_s}
\newcommand{\f}{\vec{f}}
\newcommand{\hf}{\hat{f}}
\newcommand{\hx}{\hat{x}}
\newcommand{\hy}{\hat{y}}
\newcommand{\hz}{\hat{z}}
\newcommand{\hw}{\hat{w}}
\newcommand{\tw}{\tilde{w}}
\newcommand{\hlambda}{\hat{\lambda}}
\newcommand{\hbeta}{\hat{\beta}}
\newcommand{\tG}{\widetilde{G}}
\newcommand{\tg}{\widetilde{g}}
\newcommand{\barhx}{\bar{\hat{x}}}
\newcommand{\vecx}{x_1,\ldots,x_m}
\newcommand{\xoy}{x\rightarrow y}
\newcommand{\barx}{{\bar x}}
\newcommand{\bary}{{\bar y}}
\newcommand{\tR}{{\widetilde R}}
\newcommand{\td}{{\widetilde t}}
\newcommand{\hrho}{\widehat{\rho}}
\newtheorem{theorem}{Theorem}[section]
\newtheorem{lemma}{Lemma}[section]
\newtheorem{corollary}{Corollary}[section]
\newtheorem{proposition}{Proposition}[section]
\newtheorem{definition}{Definition}[section]
\newtheorem{claim}{Claim}[section]
\newtheorem{remark}{Remark}[section]

\def\br{\break}
\def\smskip{\par\vskip 5 pt}
\def\proof{\bn {\bf Proof.} }
\def\QED{\hfill{\bf Q.E.D.}\smskip}
\def\qed{\quad{\bf q.e.d.}\smskip}
\newcommand{\cE}{\mathcal{E}}
\newcommand{\cM}{\mathcal{M}}
\newcommand{\cN}{\mathcal{N}}
\newcommand{\cJ}{\mathcal{J}}
\newcommand{\cT}{\mathcal{T}}
\newcommand{\bx}{\mathbf{x}}
\newcommand{\bp}{\mathbf{p}}
\newcommand{\bX}{\mathbf{X}}
\newcommand{\bY}{\mathbf{Y}}
\newcommand{\bP}{\mathbf{P}}
\newcommand{\bA}{\mathbf{A}}
\newcommand{\bB}{\mathbf{B}}
\newcommand{\bfM}{\mathbf{M}}
\newcommand{\bL}{\mathbf{L}}
\newcommand{\bz}{\mathbf{z}}
\newcommand{\cF}{\mathcal{F}}
\newcommand{\cR}{\mathcal{R}}
\newcommand{\bzero}{\mathbf{0}}

\newcommand{\blue}{\color{blue}}
\newcommand{\red}{\color{red}}

%% PUT YOUR TITLE PAGE INFORMATION HERE %%%
\title{Distributed Non-Convex First-Order Optimization and Information Processing: Lower Complexity Bounds and  Rate Optimal Algorithms}
\author{{Haoran Sun and Mingyi Hong \thanks{H. Sun and M. Hong are with the Department of Electrical and Computer Engineering (ECE), University of Minnesota, Minneapolis, MN 55414, USA. Email: \texttt{\{sun00111,mhong\}@umn.edu}}}}
\maketitle
\begin{abstract}
 We consider a class of popular distributed non-convex optimization problems,  in which agents connected by a network $\mathcal{G}$ collectively optimize a sum of smooth (possibly non-convex) local objective functions. We address the following question: if the agents can only access the gradients of local functions, what are the {\it fastest} rates that any distributed algorithms  can achieve, and how to achieve those rates. 

First, we show that there exist difficult problem instances, such that it takes a class of distributed first-order methods at least  $\mathcal{O}(1/\sqrt{\xi(\cG)} \times \bar{L} /{\epsilon})$  communication rounds to achieve certain $\epsilon$-solution [where $\xi(\cG)$ denotes the {spectral} gap of the graph Laplacian matrix,  and $\bar{L}$ is some Lipschitz constant]. Second, we propose {(near)} optimal methods whose rates match the developed lower rate bound (up to a ploylog  factor). The key in the algorithm design is to properly embed the classical polynomial filtering techniques into modern first-order algorithms. To the best of our knowledge, this is the first time that lower rate bounds and optimal methods have been developed for distributed non-convex optimization problems. 
\end{abstract}

\section{Introduction}\label{sec:introduction}
\subsection{Problem and  motivation}
In this work, we consider the following distributed optimization problem over a network
\begin{align}\label{eq:global:consensus}
\min_{y\in\mathbb{R}^S} \; \bar{f}(y):= \frac{1}{M}\sum_{i=1}^{M} f_i(y),
\end{align}
where $f_i(y): \mathbb{R}^{S}\to\mathbb{R}$ is a smooth and possibly non-convex function accessible to agent $i$. There is no central controller, and the $M$ agents are connected by a network defined by an {\it undirected} and {\it unweighted} graph ${\cG}=\{{\mathcal{V}}, {\mathcal{E}}\}$, with $|{\mathcal{V}}|=M$ vertices and $|{\mathcal{E}}|=E$ edges. Each agent $i$ can only communicate with its immediate neighbors, and it can access one component function $f_i$ (by ``access" we meant that it will be able to query the function and obtain its values and gradients; this notion will be defined precisely shortly).

A common way to reformulate problem \eqref{eq:global:consensus} in the distributed setting is given below.  Introduce $M$ local variables $x_1, \cdots, x_M \in \mathbb{R}^{S}$ and a concatenation of $M$ variables $x := [x_1; \cdots; x_M] \in \mathbb{R}^{SM\times1}$, and suppose the graph $\{\cV, \cE\}$ is connected,  then the following formulation is equivalent to the global consensus problem
\begin{align}\label{eq:global:consensus:equiv:2}
\min_{x\in\mathbb{R}^{SM}} \; f(x):=\frac{1}{M}\sum_{i=1}^{M} f_i(x_i),\quad \st\; x_i=x_j, \forall~(i,j)\in \cE.
\end{align}
The main benefit of the above formulation is that the objective function is now separable, and the linear constraint encodes the network connectivity pattern.  

%\vspace{-0.5cm}
\subsection{Distributed non-convex optimization}
  Distributed non-convex  optimization has gained considerable attention recently. For example, it finds applications in  training neural networks \cite{Lian17decentralized}, clustering   \cite{Forero11}, and dictionary learning\cite{Wai15_icassp}, just to name a few.  

The problem \eqref{eq:global:consensus} and \eqref{eq:global:consensus:equiv:2} have been studied extensively in the literature when $f_i$'s are all convex; see for example \cite{Nedic09subgradient, nedic2015distributed, shi2014extra}. Primal based methods such as distributed subgradient (DSG) method \cite{Nedic09subgradient}, the EXTRA method \cite{shi2014extra}, as well as primal-dual based methods such as distributed augmented Lagrangian method  \cite{jakovetic2015linear}, Alternating Direction Method of Multipliers (ADMM) \cite{BoydADMMsurvey2011,Schizas09} have been proposed. 

On the contrary, only recently there have been works addressing the more challenging problems without assuming convexity of  $f_i$; see \cite{bianchi2013convergence, Zhu-Martinez2,hong14nonconvex_admm,Lorenzo16,Hajinezhad17inexact,hong17icml,Wai15_icassp,Hajinezhad17zeroth, Daneshmand_et_al_Asilomar16,Daneshmand_et_al_ICASSP15, Zeng19distributedGD, Jiang2017,Lian17decentralized,vlaski2019distributed,vlaski2019distributed2, swenson2019annealing}. 
The convergence behavior of the distributed consensus problem \eqref{eq:global:consensus} has been studied in \cite{bianchi2013convergence,Zhu-Martinez2, Wai15_icassp}.
Reference \cite{hong14nonconvex_admm} develops a non-convex ADMM based methods for solving the distributed consensus problem \eqref{eq:global:consensus}. However the network considered therein is a star network in which the local nodes are all connected to a central controller. 
References \cite{hong17icml,Hajinezhad17inexact} propose a primal-dual based method for unconstrained problem over a connected network, and derives a global convergence rate for this setting. 
In \cite{Lorenzo16, Daneshmand_et_al_Asilomar16,Daneshmand_et_al_ICASSP15}, the authors utilize certain gradient tracking idea to solve a constrained nonsmooth distributed problem over possibly time-varying networks. The work \cite{Zeng19distributedGD} summarizes a number of recent progress in extending the DSG-based methods for non-convex problems. References \cite{Jiang2017,Lian17decentralized,Hajinezhad17zeroth} develop methods for distributed stochastic zeroth and/or first-order non-convex  optimization.  It is worth noting that the distributed algorithms proposed in all these works converge to first-order stationary solutions, which contain local maximum, local minimum and saddle points.

Recently, the authors of \cite{hong172ndorder, daneshmand2018second, swenson2019distributed, vlaski2019distributed2} have developed first-order distributed algorithms that are capable of computing second-order stationary solutions (which under suitable conditions become local optimal solutions).  Other second-order distributed algorithms such as \cite{duenner2018distributed, fang2018distributed} are design for convex problems, and they utilize high-order Hessian information about local problems.

%\vspace{-0.5cm}
\subsection{Lower and upper rate bounds analysis}\label{sub:intro:lower}
Despite all the recent interests and contributions in this field, one major  question remains open:
%\vspace{-0.5cm}
\begin{center}
	\noindent\fcolorbox{black}[rgb]{0.9,0.9,0.9}{\begin{minipage}{0.8\columnwidth}
			\begin{center}
				{\bf (Q)} ~~What is the {\it best} convergence rate achievable by {\it any} distributed algorithms for the non-convex problem \eqref{eq:global:consensus}?
			\end{center}
		\end{minipage}}
%		\vspace{-0.2cm}
	\end{center}
	Question ${\rm \bf (Q)}$ seeks to find a ``best convergence rate", which is a characterization of the smallest number of iterations  required to achieve certain {\it high-quality solutions}, among all distributed algorithms. Clearly, understanding ${\rm \bf (Q)}$  provides fundamental insights to distributed optimization and information processing. For example, the answer to ${\rm \bf (Q)}$ offers meaningful optimal estimates on the total amount of communication and computation effort required to achieve a given level of accuracy. Further, the identified optimal strategies capable of attaining the best convergence rates will also help guide the practical design of distributed information processing algorithms.

Question ${\rm \bf (Q)}$ is easy to state, but formulating it rigorously is quite involved and a number of delicate issues have to be clarified. Below we provide a high level discussion on some of these issues.	
 	
 	\noindent{\bf(1)} {\bf Fix Problem and Network Classes.} A class of problems $\mathcal{P}$ and networks $\mathcal{N}$ of interest should be fixed. Roughly speaking, in this work, we will fix $\mathcal{P}$ to be the family of  smooth unconstrained problem \eqref{eq:global:consensus}, and $\mathcal{N}$ to be the set of {\it connected} and {\it unweighted} graphs with finite number of nodes.
 	%\vspace{-0.3cm}
 	
 	 	\noindent{\bf(2)}  {\bf Characterize High-Quality Solutions.} For a properly defined error constant $\epsilon>0$,  one needs to define a  {\it high-quality solution} in distributed and non-convex setting. Differently from the centralized case, the following questions have to be addressed: Should the solution quality be evaluated based on the {\it averaged} iterates among all the agents, or on the individual iterates? Shall we include some {\it consensus measure} in the solution characterization? Different solution notion could potentially lead to different lower and upper rate bounds.%A few different notions of solutions will be discussed in detail shortly;
 	%	\vspace{-0.3cm}
 	
 	 	\noindent{\bf(3)}  {\bf Fix Algorithm Classes.} A class of algorithms $\mathcal{A}$ has to be fixed. In the classical complexity analysis in (centralized) optimization, it is common to define the class of algorithms by the information structures that they utilize \cite{nesterov05}. In the distributed and non-convex setting, it is necessary to specify {\it both} the function information that can be used by individual nodes, as well as the communication protocols that are allowed.
 %	\vspace{-0.3cm}
 
 	 	\noindent{\bf(4)}  {\bf Develop Sharp Upper Bounds.} It is necessary to develop algorithms within class $\mathcal{A}$, which possess provable and sharp global convergence rate for problem/network class $(\mathcal{P},\mathcal{N})$. These algorithms provide achievable {\it upper bounds} on the global convergence rates.
 	%\vspace{-0.3cm}
 	
 	 	\noindent{\bf(5)}  {\bf Identify  Lower Bounds.} It is important to characterize the {\it worst} rates achievable by {\it any} algorithm in class $\mathcal{A}$ for problem/network class $(\mathcal{P},\mathcal{N})$. This task involves identifying  instances in $(\mathcal{P},\mathcal{N})$ that are difficult for algorithm class $\mathcal{A}$.
 	%	\vspace{-0.3cm}
 	
 	 	\noindent{\bf(6)}  {\bf Match Lower and Upper Bounds.} The key task is to investigate whether the developed algorithms are {\it rate optimal}, in the sense that rate upper bounds derived in {\bf(4)}
match the {\it worst-case} lower bounds developed in {\bf (5)}. Roughly speaking, matching two bounds requires that for the class of problem and networks $(\mathcal{P},\mathcal{N})$, 
the following quantities should be matched between the lower and upper bounds: {\it i)} the order of the  error constants $\epsilon$; {\it ii)} the order of problem parameters such as $M$, or that of network parameters such as the spectral gap, diameter, etc.

 Convergence rate analysis (aka iteration complexity analysis) for convex problems dates back to Nesterov, Nemirovsky and Yudin  \cite{nesterov83, nemirovsky83}, in which lower bounds and optimal {\it first-order} algorithms have been developed; also see \cite{nesterov04}. In recent years, many accelerated {\it first-order} algorithms achieving those lower bounds for different kinds of convex problems have been derived; see e.g., \cite{Beck:2009:FIS:1658360.1658364,Ouyang15,tseng08acc}, including those developed for distributed convex optimization \cite{Jakovetic14}. In those works,   the problem is to optimize $\min_x f(x)$ with convex $f$, the optimality measure used is $f(x)-f(x^*)$, and the lower bound can be expressed as \cite[Theorem 2.2.2]{nesterov04}
 \begin{align}
 f(x^t)-f(x^*)\le \frac{\|x^0-x^*\|L}{(t+2)^2},
 \end{align}
 where $L$ is the Lipschitz constant for $\nabla f$; $x^*$ (resp. $x^0$) is the global optimal solution (resp. the initial solution); $t$ is the iteration index. 
 Therefore to achieve $\epsilon$-optimal solution in which $f(x^t)-f(x^*)\le \epsilon$, one needs  $\sqrt{\frac{\|x^*-x^0\|L}{\epsilon}}$ iterations.
  Recently the above approach has been extended to distributed strongly convex optimization in \cite{scaman2017optimal}. In particular,  the authors consider problem \eqref{eq:global:consensus} in which each $f_i$ is strongly convex, and they provide lower and upper rate bounds for a class of algorithms in which the local agents can utilize both $\nabla f_i(x)$ and its Fenchel conjugate $\nabla^* f_i(x)$. 
  We note that this result is not directly related to the class of ``first-order" method, since beyond the first-order gradient information, the Fenchel conjugate $\nabla^* f_i(x)$ is also needed, but computing this quantity requires performing certain exact minimization, which itself involves solving a strongly convex optimization problem.  
   Other related works in this direction also include \cite{uribe2017optimal} and \cite{scaman2018optimal}. In particular, the work \cite{scaman2018optimal} is a non-smooth extension of \cite{scaman2017optimal}, where the lower complexity bound under the  Lipschitz continuity of the global and local objective function are discussed and the optimal algorithm is proposed.

  	\begin{table*}[t]
  		\centering
  		\begin{tabular}{|c|c|c|c|}
  			\hline
  			\multirow{2}{*}{{\sf{Network Instances}}}&\multicolumn{2}{|c|}{{\sf{Problem Classes}}}&\\
  			\cline{2-3}
  			& {\sf {Uniform Lipschitz}} $U$& {\sf {Non-uniform Lipschitz}} $\{L_i\}$ &  {\sf{Rate Achieving Algorithm}}\\
  			\cline{1-4}
  			{\sf{Complete/Star}} & $\mathcal{O}(U/{\epsilon})$ & $\cO(1/\epsilon\times \sum_{i} L_i/M)$   & D-GPDA (proposed)\\
  			\hline
  			{\sf{Random Geometric}} & $\mathcal{\widetilde{O}}(U\sqrt{M}/{(\sqrt{\log{M}}\epsilon)})$ &   ${\widetilde \cO}(\sqrt{M}/(\sqrt{\log(M)}\epsilon)\times \sum_{i} L_i/M)$&  xFILTER (proposed)\\
  			\hline
  			{\sf {Path/Circle}} & $\mathcal{\widetilde{O}}(U M /{\epsilon})$ &  ${\widetilde \cO}(M/\epsilon\times \sum_{i} L_i/M)$ &  xFILTER (proposed)\\
  			\hline
  				{\sf {Grid}} & $\mathcal{\widetilde{O}}(U \sqrt{M} /{\epsilon})$ &  ${\widetilde \cO}(\sqrt{M}/\epsilon\times \sum_{i} L_i/M)$ &  xFILTER (proposed)\\
  			\hline
  			{\sf Centralized} & $\mathcal{O}(U/{\epsilon})$ & $\cO(1/\epsilon\times \sum_{i} L_i/M)$ & Gradient Descent\\
  			\hline
  		\end{tabular}
  		\caption{The main results of the paper when specializing to a few popular graphs. The entries show the best rate bounds achieved by the proposed algorithms (either D-GPDA or {xFILTER}) for a number of specific graphs and problem class; $L_i$ is the Lipschitz constant for $\nabla f_i$ [see \eqref{eq:Lip}]; for the uniform case $U =L_1,\cdots, L_M$. %; for the non-uniform case we assume that $\ell\le L_i\le u \ell$ for some constant $u$ and $\ell$; 
  			For the uniform Lipschitz  the  lower rate bounds derived for the particular graph matches the upper rate bounds (we only show the latter in the table). The last row shows the rate achieved by the centralized gradient descent algorithm. The notation $\tilde{\cO}$ denotes big $\cO$ with some polynomial in logarithms, {i.e, use $\widetilde{\mathcal{O}}$ to denote $\mathcal{O}(\log(M))$ where $M$ is the problem dimension}.} \label{table:result}
%  		\vspace{-0.4cm}
  	\end{table*}
  	
  When the problem becomes non-convex, the size of the gradient function can be used as a measure of solution quality. In particular, let $h_T^*:=\min_{0\le t\le T}\|\nabla f(x^t)\|^2$, then it has been shown that the classical (centralized) gradient descent (GD) method achieves the following rate \cite[page 28]{nesterov04}
  \begin{align*}
 h_T^*\le \frac{c_0 L(f(x^0)-f(x^*))}{T+1}, \; \mbox{where $c_0>0$ is some constant.}
  \end{align*}
 {It has been shown in {\cite{cartis2010complexity}} that the above rate is (almost) tight for GD. Recently, \cite{carmon2017lower} has further shown that the above rate is optimal for {\it any} first-order methods that only utilize the gradient information, when applied to problems with Lipschitz gradient.} However, no lower bound analysis has been developed for distributed non-convex problem \eqref{eq:Lip:original}; there are even not many algorithms that provide achievable upper rate bounds (except for the recent works  \cite{hong14nonconvex_admm,hong17icml,tian2018asy,Daneshmand18}), not to mention any analysis on the tightness/sharpness of these upper bounds.

%\vspace{-0.3cm}
\subsection{Contribution of this work}
In this work, we address various issues that arise in answering ${\rm \bf (Q)}$. Our main contributions are given below:

\noindent{\bf 1)} We identify a class of non-convex problems and networks $(\cP,\cN)$, a class of distributed first-order algorithms $\cA$, and rigorously define the 
$\epsilon$-optimality gap that measures the progress of the algorithms;

\noindent{\bf 2)} We develop the first lower complexity bound for class $\cA$ to solve class $(\cP,\cN)$: To achieve $\epsilon$-optimality, it is necessary for any $a\in\cA$ to perform  $\mathcal{O}(1/\sqrt{\xi(\cG)} \times \bar{L} /{\epsilon})$ rounds of communication among all the nodes, where $\xi(\cG)$ represents certain spectral gap of the graph Laplacian matrix,  and $\bar{L}$ is the averaged Lipschitz constants of the gradients of local functions.  On the other hand, it is necessary for any such algorithm to perform  $\mathcal{O}(\bar{L} /{\epsilon})$ rounds of computation among all the nodes.

\noindent{\bf 3)} We design two algorithms belonging to $\cA$, one based on primal-dual optimization scheme, the other based on a novel  {\it appro\underline{x}imate \underline{filt}ering -then- pr{\underline e}dict and t{\underline r}acking} (xFILTER) strategy, both of which achieve $\epsilon$-optimality condition with provable global rates [in the order of $\cO(1/\epsilon)$];

\noindent{\bf 4)} %We analyze the rate bounds of the proposed methods, as well as a number of its refinements. 
We show that the xFILTER is an optimal method in $\cA$ for problem class $(\cP, \cN)$ as well as a number of its refinements, in that they precisely achieve the lower complexity bounds that we derived (up to a ploylog factor).  

In Table \ref{table:result}, we specialize some key results developed in the paper to a few popular graphs. 

\noindent{\bf Notations.} For a given symmetric matrix $B$, we use $\lambda_{\max}(B)$, $\lambda_{\min}(B)$ and $\underline{\lambda}_{\min}(B)$ to denote the maximum, the minimum and the minimum {\it nonzero} eigenvalues; We use $I_P$ to denote an identity matrix with size $P$, and use $\otimes$ to denote the Kronecker product. We use $[M]$ to denote the set $\{1,\cdots, M\}$. {For a vector $x$ we use $x[i]$ to denote its $i$th element.} We use $\widetilde{\mathcal{O}}$ to denote $\mathcal{O}(\log(M))$ where $M$ is the problem dimension. We use $i\sim j$ to denote two connected nodes $i$ and $j$, i.e., for a graph $\cG:=\{\cV, \cE\}$, $i\sim j$ if $i\ne j$, and $(i,j)\in \cE$. 

%\vspace{-0.2cm}
\section{Preliminaries}\label{sub:prelim}
\subsection{The class $\cP$, $\cN$, $\cA$}\label{sub:class}
We present the classes of problems, networks and algorithms to be studied, as well as some useful results. We parameterize these classes using a few key parameters
so that we can specify their subclasses when needed.

\noindent{\bf Problem Class.} A problem is in class $\cP^{M}_L$  if it satisfies the following conditions. 
\begin{itemize}
%	\vspace{-0.1cm}
	\item[A1.] The objective is an average of $M$  functions; see \eqref{eq:global:consensus}. 
	\item [A2.] Each component function $f_i(x)$'s has Lipschitz gradient:  
%	\begin{subequations}\label{eq:Lip}
			\begin{align}
			&\hspace{-0.5cm}\|\nabla f_i(x_i) - \nabla f_i(z_i)\|\le L_i\|x_i-z_i\|, \; \forall~x_i,z_i\in\mathbb{R}^S,\;\forall~i, \label{eq:Lip}
			\end{align}
%	\end{subequations}}	
where $L_i\ge 0$ is the {\it smallest} positive number such that the above inequality holds true. Define $\bar{L}:=\frac{1}{M}\sum_{i=1}^{M}L_i$, $L_{\max}:=\max_{i}L_i$, and $L_{\min}$ similarly. 

 Define the matrix of Lipschitz constants as:
 \begin{align}\label{eq:L}
 L:= \mbox{diag}([L_1,\cdots, L_M])\otimes I_S\in\mathbb{R}^{M S\times MS}. 
 \end{align}
%Throughout the paper we will assume that $L_i\ge 0, \; \forall~i$, so $L\succ 0$.  
	\item [A3.] The function $f(x)$ is lower bounded over $x\in\mathbb{R}^{MS}$, i.e., 
	\begin{align}\label{eq:inf}
	\underline{f} := \inf_x f(x)>-\infty.
	\end{align}
\end{itemize}
%\vspace{-0.3cm}
These assumptions are rather mild. For example an $f_i$ satisfies [A2-A3] is not required to be second-order differentiable.  Below we provide a few non-convex functions that satisfy Assumption [A2-A3], and each of those can be the component function $f_i$'s. Note that the first four functions are of particular interest in learning neural networks, as they are commonly used as  activation functions. 

	\noindent{\bf (1)} The $\mbox{sigmoid}$ function is given by
	$\mbox{sigmoid}(x) = \frac{1}{1+e^{-x}}.$
	%We have {$\mbox{sigmoid}(x)\ge 1$, $\mbox{sigmoid}''(x)\in[-1,1]$}, therefore [A2-A3] are true  with $L\le 1$.
	We have {$\mbox{sigmoid}(x)\ge 0$, $\mbox{sigmoid}''(x)\in(-1,1)$}, therefore [A2-A3] are true  with $L\le 1$.  
	
	\noindent{\bf (2)} The  $\arctan$ function satisfies {$\arctan(x)\in (-\frac{\pi}{2}, \frac{\pi}{2})$}, $\arctan''(x) = \frac{-2x}{(x^2 +1)^2} \in [-1,1]$.  So [A2-A3] hold with $L\le 1$. 
	
	\noindent{\bf (3)} The  $\tanh$ function satisfies
	{$\tanh(x)\ge -1, \; \tanh''(x) \in [-1,1],$} so [A2-A3] hold with {$L\le 1$}.  
	
	\noindent{\bf (4)} The $\mbox{logit}$  function is related to the $\tanh$ function as follows
	$$2\mbox{logit}(x) =\frac{2e^{x}}{e^x+1}=1+ \tanh(x/2),$$ then  Assumptions [A2-A3] are again satisfied.
	
	\noindent{\bf (5)} The $\log(1+ x^2)$ function has applications in structured matrix factorization \cite{fu16}. Clearly it is lower bounded. Its second-order derivative is also bounded.

	\noindent{\bf (6)} Other functions like $\sin(x)$, $\mbox{sinc}(x)$, $\cos(x)$  are easy to verify.
		Consider $f(x):= -x_1 x_2 + (x_1-1)^2_{+}+(-x_1-1)^2_{+}$ where $(z)^2_{+}:=\max\{0,z\}^2$. This function is interesting because it is not second-order differentiable; nonetheless we can verify that [A2-A3] are satisfied with $L=\sqrt{2}+1$.

\noindent{\bf Network Class.} Let $\cN$ denote a class of networks represented by an {\it undirected} and {\it unweighted} graph ${\cG}=\{{\mathcal{V}}, {\mathcal{E}}\}$, with $|{\mathcal{V}}|=M$ vertices and $|{\mathcal{E}}|=E$ edges, and edge weights all being $1$. In this paper the term `network' and `graph' will be used interchangeably. Also, we use $\cN^{{M}}_{D}$ to denote a class of network similarly as above, but with $M$ nodes and a diameter of $D$, defined below [where $\mbox{dist}(\cdot)$ indicates the distance between two nodes]: 
\begin{align}\label{eq:dist}
D:= \max_{u,v\in \cV}\mbox{dist}(u,v). 
\end{align}
Following the convention in \cite{Chung97}, we define a number of graph related quantities below. 
First, define the {\it degree} of node $i$ as $d_i$, and define the averaged degree as:
\begin{align}\label{eq:average:degree}
\bar{d} := \frac{1}{M}\sum_{i=1}^{M}d_i.
\end{align}
Define the incidence matrix (IM)  $A\in\mathbb{R}^{E\times M}$ as follows: if $e\in\mathcal{E}$ and it connects vertex $i$ and $j$ with $i>j$, then $A_{ev}=1/\sqrt{d_v}$ if $v=i$, $A_{ev}=-1/\sqrt{d_v}$ if $v=j$ and $A_{ev}=0$ otherwise; see the definition in \cite[Theorem 8.3]{Chung97}. Using these definitions, the  {\it graph Laplacian matrix}  and the {\it degree matrix} are defined as follows (see \cite[Section 1.2]{Chung97}):
\begin{align}\label{eq:NL}
\cL:=A^T A\in\mathbb{R}^{M\times M}, \quad \mbox{and}\quad P:=\mbox{diag}[d_1,\cdots, d_M]\in{\mathbb{R}^{M\times M}}.
\end{align}
In particular, the elements of the Laplacian are given as:
\begin{align*}
[{\cL}]_{ij} = \left\{
\begin{array}{ll}
1 & \mbox{if}~i=j\\
-\frac{1}{\sqrt{d_i d_j}} & \mbox{if}~i\sim j, i\ne j\\
0 & \mbox{otherwise}.
\end{array}
\right.
\end{align*}
We note that the graph Laplacian defined here is sometimes known as the {\it normalized} graph Laplacian in the literature, but throughout this paper we follow the convention used in the classical work \cite{Chung97} and simply refer it as the {\it graph Laplacian}. For convenience, we also define a scaled version of the IM:
\begin{align}\label{eq:nIM}
F: = A P^{1/2} \in\mathbb{R}^{E\times M}. 
\end{align}
It is known that IM and scaled IM satisfy the following (where $\mathbbm{1}\in\mathbbm{R}^M$ is an all one vector):
\begin{align}\label{eq:null}
F \mathbbm{1} = A P^{1/2} \mathbbm{1} = 0.
\end{align}
Define the second smallest eigenvalue of $\cL$, as $\underline{\lambda}_{\min}(\cL)$: 
\begin{align}
	\underline{\lambda}_{\min}(\cL) = \inf_{x: \sum_{i=1}^{M}x_i d_i  =0, x\ne 0} \frac{x^T \cL x}{\sum_{i=1}^{M}x^2_i d_i}.
	\end{align}
Then the {\it spectral gap} of the graph $\cG$ can be defined below:  
\begin{align}\label{eq:gap}
\xi(\cG)=\frac{\underline{\lambda}_{\min}(\cL)}{\lambda_{\max}(\cL)}\le 1.
\end{align}

\noindent{\bf Algorithm Class.}  Define the {\it neighbor set} for node $i\in\cE$ as 
\begin{align}\label{eq:neighbor}
\cN_i := \{i\mid i\sim j, j\ne i\}.
\end{align}We say that a {\it distributed, first-order} algorithm is in class $\cA$ if it satisfies the following conditions. 
\begin{itemize}
\item [B1.] At iteration $0$, each node can obtain some network related constants, such as $M$, $D$, eigenvalues of the graph Laplacian $\cL$, etc. 

\item [B2.] At iteration $t+1$, {\it each node} $i\in[M]$ first conducts a communication step by broadcasting the local $x^t_i$ to all its neighbors, through a function $Q^t_i(\cdot):\mathbb{R}^{S}\to \mathbb{R}^{S}$. 
Then each node will generate the new iterate, by combining the received message with its past gradients using a function $W^t_i(\cdot)$: 
\begin{align}\label{eq:span}
	\hspace{-0.3cm} v^{t}_i =\hspace{-0.7cm}\underbrace{Q^t_i(x^t_i)}_{\mbox{communication step}}\hspace{-0.7cm}, \;	x^{t+1}_i \in \hspace{-0.2cm}\ \underbrace{W^t_i\left(\{\{v^k_j\}_{j\in\cN_i}, \nabla f_i(x^k_i), x^k_i\}_{k=1}^t\right)}_{\mbox{computation step}}.
	\end{align} 
	In this work, we will focus on the case where the $Q^t_i(\cdot)$'s and $W^t_i(\cdot)$'s are linear operators. 
\end{itemize}

Clearly $\cA$ belongs to the class of {\it first-order} methods because only local gradient information is used. It is also a class of {\it distributed}  algorithms because at each iteration the nodes only communicate with their immediate neighbors. %Since the linear combination coefficients can be arbitrarily chosen in computing $v^r_i$, {\it at each iteration} the nodes has the flexibility in choosing the subset of its neighbors to communicate, as well as how to combine their outputs. 

Additionally, in practical distributed algorithms such as DSG, ADMM or EXTRA,  nodes are dictated to use a {\it fixed strategy} to linearly combine all its neighbors' information.  To model such a requirement, below we consider a slightly restricted algorithm class $\cA'$, where we require each node to use the same coefficients to combine its neighbors (note that allowing the nodes to use a fixed but {\it arbitrary} linear combination is also possible, but the resulting analysis will be more involved). 

In particular, we say that a {\it distributed, first-order} algorithm is in $\cA'$ if it satisfies B1 and the following:
\begin{itemize}
\item [B2'.] At iteration $t+1$, {\it each node} $i\in[M]$ performs:	
\begin{align}\label{eq:span:2}
	\hspace{-0.4cm} v^{t}_i =Q^t_i(x^t_i), \;	x^{t+1}_i \in \hspace{-0.2cm}\ W^t_i\left(\{\sum_{j\in\cN_i}v^t_j, \nabla f_i(x^k_i), x^k_i\}_{k=1}^t\right).	
	\end{align} 
 \end{itemize}

We remark that, in both algorithm classes, {\it one round} of communication occurs at each iteration, where {\it each} node broadcasts its local variable $x^t_i$ once. Therefore, the total iteration number is the same as the total communication rounds. However, the total times that the entire gradient $\{\nabla f_i(x_i)\}_{i=1}^{M}$ is evaluated could be {\it smaller} than the total iteration number/communication rounds. This is because when we compute $x^{t+1}_i$, the operation $W^t_i(\cdot)$ can set the coefficient in front of $\nabla f_i(x^r_i)$ to be zero, effectively {\it skipping} the local gradient computation.

\subsection{Solution Quality Measure}\label{sub:error}
Next we provide  definitions for the quality of the solution. Note that since we consider using first-order methods to solve non-convex problems, it is expected that in the end some first-order stationary solution with small $\|\nabla f\|$ will be computed. 

Our first definition is related to a {\it global variable} $y^t\in\mathbb{R}^{S}$. We say that $y^t$ is a {\it global} $\epsilon$-solution if the following holds:
%\begin{subequations}
\begin{align}\label{eq:criteria:global}
&y^t \in \mbox{span}\left\{x^t_i\right\}_{i=1}^{M}, \; \min_{t\in[T]}\|\nabla g(y^t)\|^2 \le \epsilon.
\end{align}
%\end{subequations}
This definition is conceptually simple and it is identical to the centralized criteria in Section \ref{sub:intro:lower}. However it has the following issues. First, no global 
variable $y^t$ will be formed in the entire network, so criteria \eqref{eq:criteria:global} is difficult to evaluate. Second, there is no characterization of how close the local variables $x^t_i$'s are. To see the second point, consider the following toy example.

\noindent{\it Example 1:} Consider a network with $M=2$ and $f_1 (y)= - y^2$ and $f_2(y) = y^2$. Suppose that the local variables take the following values: $x^T_1 = -10$ and $x^T_2 = 10$. Then if we pick $y^T = (x^T_1+x^T_2)/2 =0$, we have 
$$\nabla g(y^T) =\frac{1}{2}(\nabla f_1(y^T)+ \nabla f_2(y^T))=0.$$
This suggests that at iteration $T$, there exists one linear combination that makes measure  \eqref{eq:criteria:global} precisely zero. However one can hardly say that the current solution $(x^T_1,x^T_2)=(-10,10)$ is a good solution for problem \eqref{eq:global:consensus:equiv:2}.  \hfill $\square$
 
To address the above issue, we provide a second definition which is directly related to {\it local variables} $\{x_i\in\mathbb{R}^{S}\}_{i=1}^{M}$. At a given iteration $T$, we say that $\{x^T_i\}$ is a {\it local} $\epsilon$-solution if the following holds:
\begin{align}\label{eq:criteria:local}
\hspace{-0.6cm}h^*_T:=\min_{t\in[T]}\bigg\|\sum_{i=1}^{M}\frac{\nabla f_i(x^t_i)}{M}\bigg\|^2 +  \frac{1}{M \underline{\lambda}_{\min}(P^{1/2}\cL P^{1/2})}\sum_{(i,j): i \sim j}{\sqrt{L_i L_j}}\|x^t_i-x^t_j\|^2\le \epsilon.
\end{align}

Clearly this definition takes into consideration the consensus error as well as the size of the local gradients. 
When applied to Example 1, this measure will be large.   Note that the constant $\frac{1}{M \underline{\lambda}_{\min}(P^{1/2}\cL P^{1/2})}$ is needed to balance the two different measures. Also note that the $``\min_{t\in[T]}"$ operation is needed to track the best solution obtained before iteration $T$, because the quantity inside this operation may not be monotonically decreasing. 

In our work we will focus on providing answers to the following specific version of question $\bf (Q)$: 
\begin{center}
	\noindent\fcolorbox{black}[rgb]{0.9,0.9,0.9}{\begin{minipage}{0.98\columnwidth}
			\begin{center}
	For {\it any} given $\epsilon>0$, what is the minimum iteration $T$ (as a function of $\epsilon$) needed for any algorithm in class $\cA$  (or class $\cA$') to solve instances in 	classes $(\cP,\cN)$, so to achieve $h_{T}^*\le \epsilon$?
		\end{center}
	\end{minipage}}
	%		\vspace{-0.2cm}
\end{center}
 
\subsection{Some Useful Facts and Definitions}\label{sub:facts}
Below we provide a few facts about the above classes.

\noindent{\bf On Lipschitz constants.} Assume that each $f_i$ has Lipschitz continuous gradient with constant $L_i$ in \eqref{eq:Lip}. 
Then we have : 
	\begin{align}\label{eq:Lip:original}
	\|\nabla \barf(y_1)-\nabla \barf(y_2)\|\le \sum_{i=1}^{M}\frac{1}{M}L_i\|y_1-y_2\|:=\bar{L} \|y_1-y_2\|, \; \forall~y_1,\;y_2\in\mathbb{R}^{S},
	\end{align}
where $\bar{L}$ is the average of the local Lipschitz gradients. We also have the following 
\begin{align*}
\|\nabla f(x) -\nabla f(z)\|^2& =\frac{1}{M^2}\sum_{i=1}^{M}\|\nabla f_i(x_i) - \nabla f_i(z_i)\|^2,\; \forall~x_i,\;z_i\in\mathbb{R}^S\nonumber%\le \frac{1}{M^2}\sum_{i=1}^{M}L^2_i\|x_i-z_i\|^2 \le \frac{1}{M^2} (x-z)^T L^2 (x-z)
\end{align*}
which implies 
\begin{align}\label{eq:Lip:extended}
\hspace{-0.5cm}\|\nabla f(x)-\nabla f(z)\|\le \frac{1}{M}\|L (x-z)\|, \quad \forall~x,z\in\mathbb{R}^{M S},
\end{align}
where the matrix $L$ is defined in \eqref{eq:L}.

\noindent{\bf On Quantities for Graph $\mathcal{G}$.} This section presents a number of properties for a given graph $\cG$. 
Define the following matrices:{
	\begin{align}\label{eq:Upsilon}
	\hspace{-0.3cm}\Sigma:=\mbox{diag}[\sigma_1,\cdots, \sigma_E]\succ 0, \; \Upsilon:=\mbox{diag}([\beta_1,\cdots, \beta_M])\succ 0. 
	\end{align}}
Define $B\in\mathbb{R}^{E\times M} =|F|$ where the absolute value is taken component-wise. Then we have the following:{
\begin{align}\label{eq:sum:A:B}
\hspace{-0.2cm}&\frac{1}{2}\left(F^T F + B^T  B\right) = P = \mbox{diag}[d_1,\cdots, d_M]\in\mathbb{R}^{M\times M} \\ 
\hspace{-0.2cm}	&\frac{1}{2}\left(F^T \Sigma^2 F + B^T  \Sigma^2 B\right) =\mbox{diag}\left(\bigg\{\sum_{j: i\sim j} \sigma^2_{ij} \bigg\}_{j\in\cN}\right):= \Delta \nonumber, %= P
	\end{align}}
where ${P}$ is the {\it degree matrix} defined in \eqref{eq:NL}. 

For two diagonal matrices $\Upsilon^2$ and $\Sigma^2$ of appropriate sizes, the {\it generalized Laplacian} (GL) matrix is defined as:
\begin{align}\label{eq:normalized::gen:L}
\cL_G = \Upsilon^{-1} F^T \Sigma^2 F \Upsilon^{-1},
\end{align}
and its elements are given by:
\begin{align*}
[{\cL_G}]_{ij} = \left\{
\begin{array}{ll}
\frac{\sum_{q:i\sim q}\sigma^2_{iq}}{\beta^2_i} & \mbox{if}~i=j\\
-\frac{\sigma^2_{ij}}{{\beta_i\times \beta_j}} & \mbox{if}~(ij)\in\cE, i\ne j\\
0 & \mbox{otherwise}
\end{array}
\right..
\end{align*}
Define a diagonal matrix $K\in\mathbb{R}^{E\times E}$ as below:
\begin{align}\label{eq:K}
[K]_{e, q} = \left\{
\begin{array}{ll}
\sqrt{L_i L_j}& \mbox{if}~e=q, \;\mbox{and}\; e=(i,j)\\
0 & \mbox{otherwise}
\end{array}
\right..
\end{align}
Then when specializing $\Upsilon = P^{1/2} L^{1/2}$ and $\Sigma^2 = K$, the GL matrix becomes:
\begin{align}\label{eq:tilde:L}
{\ctL} := L^{-1/2} P^{-1/2} F^T K F P^{-1/2} L^{-1/2}. 
\end{align}
Note that if any diagonal element in the matrix $L$ is zero, then $L^{-1}$ denotes the Moore - Penrose  matrix pseudoinverse.  
Similarly, when specializing $\Upsilon = L^{1/2}$ and $\Sigma^2 = K$, then the GL matrix  becomes:
\begin{align}\label{eq:hat:L}
\widehat{\cL} := L^{-1/2} F^T K F L^{-1/2}. 
\end{align}
These matrices will be used later in our derivations.

Below we list some useful results about the Laplacian matrix \cite{Butler2016,Chung97,Duchi12}. First, all eigenvalues  of ${\cL}$ lie in the interval $[0, \; 2]$. %and that $\underline{\lambda}_{\min}(\cL_N)\le M/(M-1)$.  
Also because $\underline{\lambda}_{\min}(\cL) = \underline{\lambda}_{\min}(P^{-1/2} F^T F P ^{-1/2})$, we have
\begin{align}
\underline{\lambda}_{\min}(\cL)  \le \underline{\lambda}_{\min}(F^T F)\label{eq:compare:lambda:min}. 
\end{align}
Also we have that \cite[Lemma 1.9]{Chung97} 
\begin{align}\label{eq:path:star:lambda}
\underline{\lambda}_{\min}(\cL)\ge \frac{1}{D \sum_{i}d_i}. 
\end{align}

The eigenvalues of $\cL$ for a number of special graphs are given below:

\noindent {\bf 1)} {\bf Complete Graph:}  The eigenvalues are $0$ and $M/(M-1)$ (with multiplicity $M-1$), so  $\xi(\cG)=1$;

\noindent{\bf 2)}  {\bf Star Graph:} The eigenvalues are $0$ and $1$ (with multiplicity $M-2$), and $2$, so $\xi(\cG)={1}/{2}$;

\noindent{\bf 3)}  {\bf Path Graph:} The eigenvalues are $1-\cos(\pi m/(M-1))$ for $m={0,}\; 1,\cdots, M-1$, and {$\xi(\cG)\ge 1/M^2$.} 

\noindent{\bf 4)}  {\bf Cycle Graph:} The eigenvalues are $1-\cos(2 \pi m/M)$ for $m={0,}\; 1,\cdots, M-1$, and {$\xi(\cG)\ge 1/M^2$.} 

\noindent{\bf 5)}  {\bf Grid Graph:} The grid graph is obtained by placing the nodes on a $\sqrt{M}\times \sqrt{M}$ grid, and connecting nodes to their nearest neighbors. We have
$\xi(\cG)\ge 1/M$.

\noindent{\bf 6)}  {\bf Random Geometric Graph:} Place  the nodes uniformly in $[0,1]^2$ and connect any two nodes separated by a distance less than a radius $R\in(0,1)$. Then if the connectivity radius $R$ satisfies \cite{Duchi12}
\begin{align}
R = \Omega\left(\sqrt{\log^{1+\epsilon}(M)/M}\right), \quad \mbox{for any}~\epsilon>0,
\end{align} 
then with high probability
\begin{align}
\xi(\cG)=\mathcal{O}\left({\frac{\log(M)}{M}}\right).
\end{align}

\section{Lower Complexity Bounds}\label{sec:lower}
In this section we develop the lower complexity bounds for algorithms in class $\mathcal{A}$ to solve problems $\mathcal{P}^{M}_L$ over network $\mathcal{N}$. We will mainly focus on the case where $f_i$'s have uniform Lipschitz constants, that is, we assume that 
$$L_i =U, \quad \forall~i\in[M],$$
and we denote the resulting problem class as $\cP^{M}_U$.  
At the end of this section, generalization to the {non-uniform} case will be briefly discussed.

Our proof combines ideas from the classical proof in Nesterov \cite{nesterov05}, as well as two recent constructions \cite{carmon2017lower} (for centralized non-convex problems) and \cite{scaman2017optimal} (for strongly convex distributed problems). {Our construction differs from the previous works in a number of ways, in particular, the constructed functions are only first-order differentiable, but not second-order differentiable. Further, we use the local-$\epsilon$ solution \eqref{eq:criteria:local} to measure the quality of the solution, which makes the analysis more involved compared with the existing {\it global} error measures in \cite{nesterov05,carmon2017lower,scaman2017optimal}. }

To begin with, we construct the following two non-convex functions
\begin{align}\label{eq:h}
{h}(x):=\frac{1}{M}\sum_{i=1}^{M} h_i(x_i), \quad f(x):=\frac{1}{M}\sum_{i=1}^{M} f_i(x_i),
\end{align}
as well as the corresponding versions that evaluate on a ``centralized" variable $y$
\begin{align}\label{eq:bar_h}
{\bar{h}}(y):=\frac{1}{M}\sum_{i=1}^{M} h_i(y), \quad \bar{f}(y):=\frac{1}{M}\sum_{i=1}^{M} f_i(y).
\end{align}
Here  we have $x_i\in\mathbb{R}^T$, for all $i$, $y\in\mathbb{R}^{T}$, and $x:=(x_1,\cdots x_M)\in\mathbb{R}^{TM\times 1}$.  Later we make our construction so that functions $h$ and $\bar{h}$ are easy to analyze, while $f$ and $\bar{f}$ will be in the desired function class in $\mathcal{{P}}^{M}_U$. {Without loss of generality, in the construction we will assume $\nabla f_i$ will be Lipschitz with constant $U \in(0,1)$, for all $i\in[M]$.}

\subsection{Path Graph (${D} = M-1$)} 
\begin{figure}	
	\centering
	\includegraphics[width=0.6\linewidth]{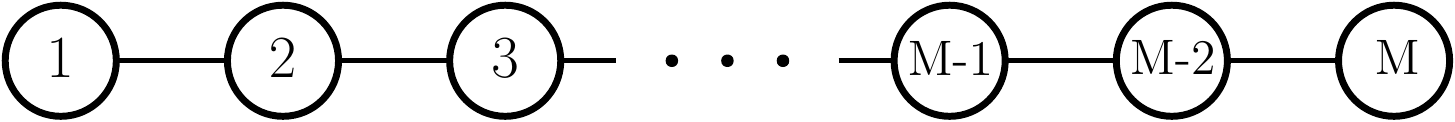}
	\caption{The path graph used in our construction.}
	\label{fig:pathgraph}
\end{figure}

First we consider the extreme case in which the nodes form a path graph with $M$ nodes and each node $i$ has its own local function $h_i$, shown in Figure \ref{fig:pathgraph}. 
\begin{figure}
	\begin{minipage}[c]{0.5\linewidth}
		\includegraphics[width=\linewidth]{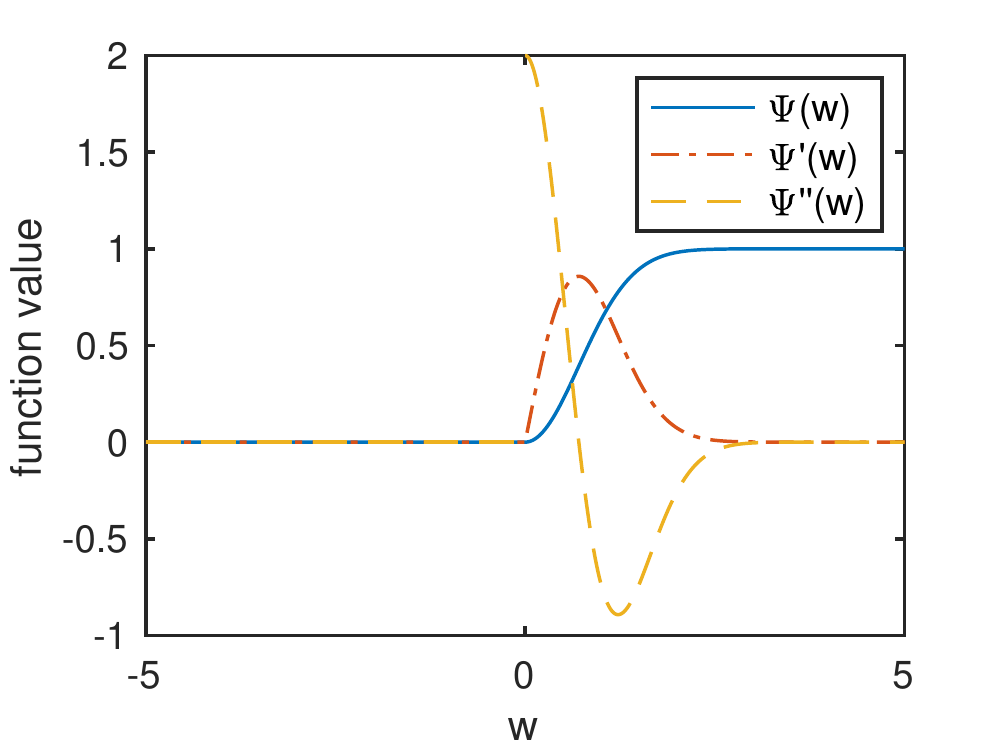}
		\caption{The functional value, and derivatives of $\Psi$.}
		\label{fig:Psi}
	\end{minipage}%
	\hfill
	\begin{minipage}[c]{0.5\linewidth}
		\includegraphics[width=\linewidth]{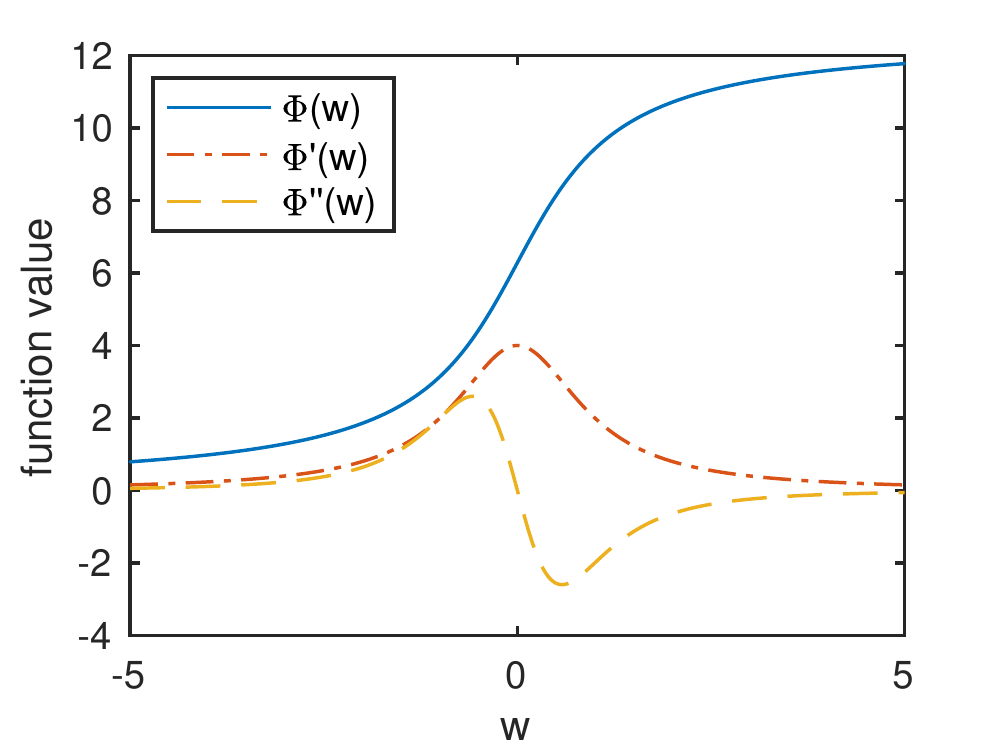}
		\caption{The functional value, and derivatives of $\Phi$. }\label{fig:Phi}
	\end{minipage}
\end{figure} 
For notational simplicity assume that $M$ is a multiple of $3$, that is $M = 3 C$ for some integer $C>0$. 
Also assume that $T$ is an odd number without loss of generality. 

Let us define the component functions $h_i$'s in \eqref{eq:h} as follows. 	
\begin{equation}  \label{eq:construct:h}
h_i(x_i) = \left\{
\begin{aligned}{}
&\Theta(x_i, 1)+3\sum_{j=1}^{\floor{T/2}} \Theta(x_i, 2j), &i\in \left[1, \frac{M}{3}\right]  \\
&\Theta(x_i, 1), &i \in \left[\frac{M}{3}+1, \frac{2M}{3}\right]  \\
&\Theta(x_i, 1)+3\sum_{j=1}^{\floor{T/2}} \Theta(x_i, 2j+1), &i \in \left[\frac{2M}{3}+1, M\right]\\
\end{aligned}
\right.
\end{equation}
{where we have defined the following functions
\begin{subequations}\label{eq:Theta}
		\begin{align}
		\Theta(x_i, j) & := \Psi(-x_i[j-1])\Phi(-x_i[j]) -\Psi(x_i[j-1])\Phi(x_i[j]), \; \forall~j\ge 2\\
		\Theta(x_i, 1)&:= -\Psi(1)\Phi(x_i[1]).
		\end{align}
\end{subequations}
}
The component functions $\Psi, \Phi :\mathbb{R}\to \mathbb{R}$ are given as below
\begin{equation*}
\Psi(w)
:= \begin{cases}
0 & w \le 0\\
1-e^{-w^2}  & w>0,
\end{cases}
\quad \mbox{and}\quad 
\Phi(w): =4 \arctan w+2\pi.
\end{equation*}
%Also as convention we have fixed  $x_i[0] := 1$. 

Suppose $x_1=x_2=\cdots =x_M = y$, then the {average} function becomes:
\begin{align*} \label{eq:hbar}
{\barh}(y)&:= \frac{1}{M}\sum_{j=1}^{M} h_i(y)  = \Theta(y, 1)+\sum_{i=2}^{T} \Theta(y, i) \\
&= -\Psi(1)\Phi\left(y[1]\right)
+ \sum_{i=2}^{T}\left[
\Psi\left(-y[i-1]\right)\Phi\left(-y[i]\right)-
\Psi\left(y[i-1]\right)\Phi\left(y[i]\right)
\right].\nonumber
\end{align*}

\begin{figure}
	\centering
	\begin{minipage}[c]{0.5\linewidth}
		\includegraphics[width=\linewidth]{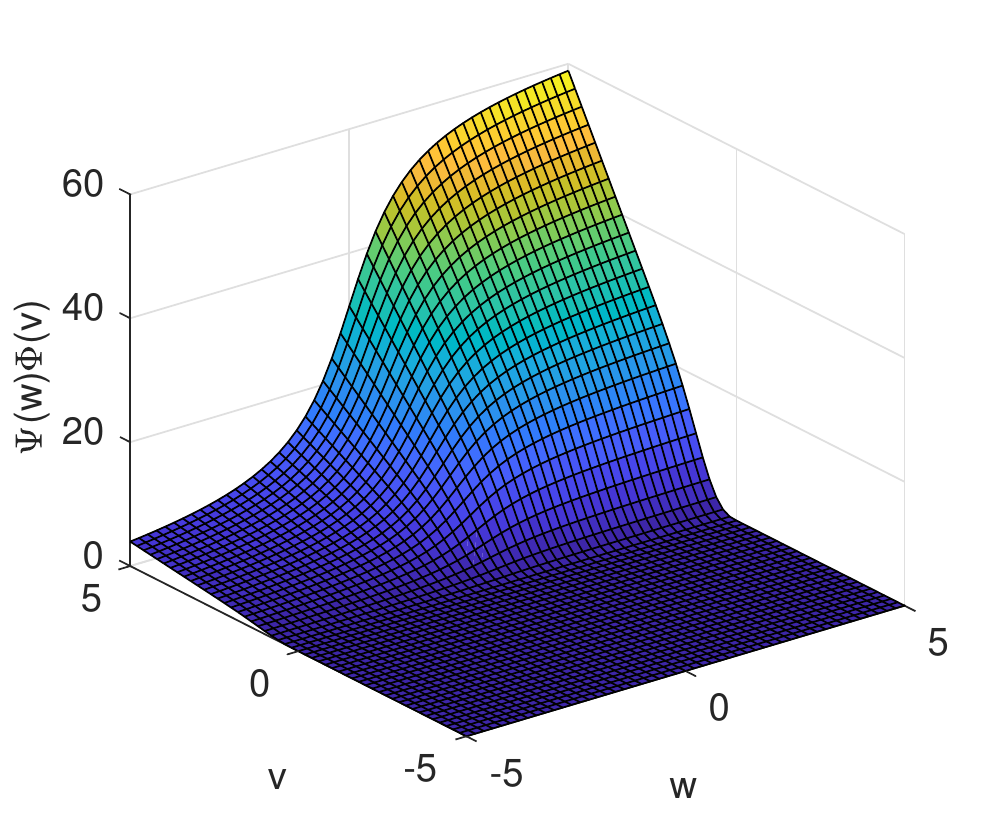}
		\caption{The functional value for $\Theta(w, v) = \Psi(w)  \Phi(v)$.}
	\end{minipage}%
\end{figure}

Further for a given error constant $\epsilon>0$ and a given averaged Lipschitz constant $U\in (0,1)$,  let us define 
\begin{equation}\label{eq:f:construct}
f_i(x_i) :=\frac{150 \pi \epsilon}{U} {h}_i\left(\frac{x_i U}{75 \pi \sqrt{2\epsilon}}\right). 
\end{equation}
Therefore we also have, if $x_1=x_2=\cdots =x_M = y$, then 
\begin{equation}\label{eq:f:construct:2}
\barf(y) :=\frac{1}{M} \sum_{i=1}^{M} f_i(y) = \frac{150 \pi \epsilon}{U} {\barh}\left(\frac{y U}{75 \pi \sqrt{2\epsilon}}\right). 
\end{equation}

First we present some properties of the component functions $h_i$'s. 
\begin{lemma}  \label{lem:func_property}
	The functions $\Psi$ and $\Phi$ satisfy the following.  
	\begin{enumerate}
		\item  For all $w \le 0$, $\Psi(w) = 0$, $\Psi'(w) = 0$. %$\Psi''(w) = 0$.
		\item The following bounds hold for the functions and their first and second-order derivatives:
		%$\Psi, \Psi', \Psi'', \Phi, \Phi'$ and $\Phi''$ are bounded, with
		\begin{equation*}
		0 \le \Psi(w) < 1,
		~~ 0 \le \Psi'(w) \le \sqrt{\frac{2}{e}},
		~~ -\frac{4}{e^{\frac{3}{2}}} \le \Psi''(w) \le 2,	 
		~~ \forall w>0
		\end{equation*}
		\begin{equation*}
		~~ 0 < \Phi(w) < 4\pi,
		~~ 0 < \Phi'(w) \le 4,
		~~ -\frac{3\sqrt{3}}{2} \le \Phi''(w)  \le \frac{3\sqrt{3}}{2}, 
		~~ \forall w\in \mathbb{R}
		\end{equation*}
		\item  The following key property holds: 
		\begin{align}\label{eq:key:3}
		\Psi(w)\Phi'(v) > 1, \quad \forall~w\ge 1, \; |v|<1. 
		\end{align}	
		\item  The function $h$ is lower bounded as follows:
		$${h_i}(0) - \inf_{x_i} {h_i}(x_i) \le {10\pi T }, \; {h}(0) - \inf_{x} {h}(x) \le {10\pi T }.$$
		\item  The first-order derivative of ${\barh}$ (resp. $h_j$) is Lipschitz continuous with constant $\ell = 75\pi$ (resp. $\ell_j = 75\pi$, $\forall~i$).
	\end{enumerate}
\end{lemma}

{\bf Proof.} Property 1) is obviously true.

To prove Property 2), note that following holds for $w>0$:  
\begin{align}
\Psi(w) = 1-e^{-w^2},
~~ \Psi'(w) =  2e^{-w^2}w,
~~ \Psi''(w) =  2e^{-w^2}-4e^{-w^2}w^2, \; \forall~w>0.
\end{align}
Obviously, $\Psi(w)$ is an increasing function over $w>0$, therefore the lower and upper bounds are $\Psi(0) = 0, \Psi(\infty) = 1$; $\Psi'(w)$ is increasing on $[0, \frac{1}{\sqrt{2}}]$ and decreasing on $[\frac{1}{\sqrt{2}}, \infty]$, where $\Psi''(\frac{1}{\sqrt{2}}) = 0$, therefore the lower and upper bounds are $\Psi'(0) = \Psi'(\infty) = 0, \Psi'(\frac{1}{\sqrt{2}}) = \sqrt{\frac{2}{e}}$; $\Psi''(w)$ is decreasing on $(0, \sqrt{\frac{3}{2}}]$ and increasing on $[\sqrt{\frac{3}{2}}, \infty)$ [this can be verified by checking the signs of  $\Psi'''(w) =4e^{-w^2}w(2w^2-3)$ in these intervals]. Therefore the lower and upper bounds are $\Psi''(\sqrt{\frac{3}{2}}) = -\frac{4}{e^{\frac{3}{2}}}, \Psi''(0^{+}) = 2$, i.e.,
\begin{equation*}
0 \le \Psi(w) < 1,
~~ 0 \le \Psi'(w) \le \sqrt{\frac{2}{e}},
~~ -\frac{4}{e^{\frac{3}{2}}} \le \Psi''(w) \le 2,	 
~~ \forall w>0.
\end{equation*}
Further,  for all $w\in\mathbb{R}$, the following holds:
\begin{align}
\Phi(w) = 4\arctan w+2 \pi, 
~~ \Phi'(w) = \frac{4}{w^2+1},
%\end{align}
%\begin{align}
~~ \Phi''(w) = -\frac{8w}{(w^2+1)^2}.
\end{align}
%$$$$

Similarly, as above, we can obtain the following bounds:
	\begin{equation*}
~~ 0 < \Phi(w) < 4\pi,
~~ 0 < \Phi'(w) \le 4,
~~ -\frac{3\sqrt{3}}{2} \le \Phi''(w)  \le \frac{3\sqrt{3}}{2}, 
~~ \forall w\in \mathbb{R}.
\end{equation*}
We refer the readers to Fig. \ref{fig:Psi} -- Fig. \ref{fig:Phi} for illustrations of these functions.

To show Property 3), note that for all $w \ge 1$ and $|v| < 1$, 
$$\Psi(w)\Phi'(v)> \Psi(1)\Phi'(1) = 2 (1-e^{-1}) > 1$$
where the first inequality is true because $\Psi(w)$ is strictly increasing and $\Phi'(v)$ is strictly decreasing for all $w>0$ and $v>0$, and that $\Phi'(v) = \Phi'(|v|)$.

Next we show Property 4). Note that $0 \le \Psi(w) < 1$ and $0 < \Phi(w) < 4\pi$. Therefore we have  ${h}(0) = - \Psi(1)\Phi(0) < 0$ and using the construction in \eqref{eq:construct:h} 
\begin{align}
\inf_{x_i} {h_i}(x_i) &\ge -\Psi(1)\Phi(x_i[1])- 3\sum_{j=1}^{\floor{T/2}}  \sup_{w, v} \Psi(w)\Phi(v) \\
&\ge {-4\pi -6\pi T \ge -10\pi T},
\end{align}
where the first inequality follows $\Psi(w)\Phi(v)>0$ and second follows $\Psi(w)\Phi(v)< 4\pi$, we reach the conclusion.

{Finally we show Property 5), using the fact that a function is Lipschitz if it is piecewise smooth with bounded derivative. 
From construction \eqref{eq:construct:h}, the first-order partial derivative of $h_q(y)$  can be expressed below.

\noindent{\bf Case I)} If $i$ is even, we have 
\begin{equation}  \label{derivative_even}
\frac{\partial h_q}{\partial y[i]} = \left\{
\begin{array}{lll}
&3 \left( - \Psi\left(-y[i-1]\right)\Phi'\left(-y[i]\right) -\Psi\left(y[i-1]\right)\Phi'\left(y[i]\right) \right), &q \in [1, \frac{M}{3}]\\
&0,   &q \in [\frac{M}{3}+1, \frac{2M}{3}] \\
& 3\left( -\Psi'\left(-y[i]\right)\Phi\left(-y[i+1]\right) -\Psi'\left(y[i]\right)\Phi\left(y[i+1]\right) \right), &q \in [\frac{2M}{3}+1, M]
\end{array}
\right..
\end{equation}

\noindent{\bf Case II)} If $i$ is odd but not 1, we have 
\begin{equation}  \label{derivative_odd}
\frac{\partial h_q}{\partial y[i]} = \left\{
\begin{array}{lll}
& 3\left( -\Psi'\left(-y[i]\right)\Phi\left(-y[i+1]\right) -\Psi'\left(y[i]\right)\Phi\left(y[i+1]\right) \right), &q \in [1, \frac{M}{3}]\\
&0,   &q \in [\frac{M}{3}+1, \frac{2M}{3}] \\
&3 \left( - \Psi\left(-y[i-1]\right)\Phi'\left(-y[i]\right) -\Psi\left(y[i-1]\right)\Phi'\left(y[i]\right) \right), &q \in [\frac{2M}{3}+1, M]
\end{array}
\right..
\end{equation}

\noindent{\bf Case III)} If $i = 1$, we have  
\begin{equation}  \label{derivative_zero}
\frac{\partial h_q}{\partial y[1]} = \left\{
\begin{array}{lll}
&- \Psi(1)\Phi'(y[1]) + 3\left( -\Psi'\left(-y[1]\right)\Phi\left(-y[2]\right) -\Psi'\left(y[1]\right)\Phi\left(y[2]\right) \right), &q \in [1, \frac{M}{3}]\\
&- \Psi(1)\Phi'(y[1]),   &q \in [\frac{M}{3}+1, M] \\
\end{array}
\right..
\end{equation}

Obviously, $\frac{\partial h_q}{\partial y[i]}$ is a piecewise smooth function for any $i, q$, and it either equals zero or is separated at the non-differentiable point $y[i] = 0$ because of the function $\Psi$. 

Further, fix a point $y\in\mathbb{R}^{T}$ and a unit vector $v\in\mathbb{R}^{T}$ where $\sum_{i=1}^T v[i]^2 = 1$.  Define 
$$g_{q}(\theta; y,v):= h_q(y+\theta v)$$ 
to be the directional projection of $h_q$ on to the direction $v$ at point $y$.  We will show that there exists $\ell>0$ such that $|g_{q}{''}(0;y,v)| \le \ell$ for all $y\ne 0$ (where the second-order derivative is taken with respect to $\theta$). 

First we can compute $g_{q}{''}(0;y,v)$ as follows:
\begin{align*}{}
g_{q}^{''}\left(0;y,v\right)
&=\sum_{i_{1}, i_{2}=1}^{T}
\frac{\partial^2}{\partial y[i_1] \partial y[i_2]} h_q\left(y\right)v[i_1] v[i_2]
=  \sum_{\delta\in\left\{ 0,1, -1\right\}}
\sum_{i=1}^{T}
\frac{\partial^2}{\partial y[i]\partial y[i+\delta] }
h_q\left(y\right)v[i]v[i+\delta],
\end{align*}
where we take $v[0] := 0$ and $v[T+1]:= 0$.

The second-order partial derivative of $h_q(y)$ ($\forall y\ne 0$) is given as follows when $i$ is even: 
\begin{equation} \label{eq:1a}
\frac{\partial^2 h_q}{ \partial y[i]\partial y[i]}  = \left\{
\begin{array}{lll}
&3 \left(  \Psi\left(-y[i-1]\right)\Phi''\left(-y[i]\right) -\Psi\left(y[i-1]\right)\Phi''\left(y[i]\right) \right), &q \in [1, \frac{M}{3}]\\
&0, &q \in [\frac{M}{3}+1, \frac{2M}{3}]      \\
& 3 \left( \Psi''\left(-y[i]\right)\Phi\left(-y[i+1]\right)  -\Psi''\left(y[i]\right)\Phi\left(y[i+1]\right) \right), &q \in [\frac{2M}{3}+1, M]\\
\end{array}
\right.
\end{equation}

\begin{equation}  \label{eq:1b}
\frac{\partial^2 h_q}{ \partial y[i]\partial y[i+1]} = \left\{
\begin{array}{lll}
&0, &q \in [1, \frac{2M}{3}] \\
& 3 \left( \Psi'\left(-y[i]\right)\Phi'\left(-y[i+1]\right)  -\Psi'\left(y[i]\right)\Phi'\left(y[i+1]\right) \right), &q \in [\frac{2M}{3}+1, M]\\
\end{array}
\right.
\end{equation}

\begin{align}   \label{eq:1c}
\frac{\partial^2 h_q}{ \partial y[i]\partial y[i-1]}
= \left\{
\begin{array}{lll}
& 3 \left( \Psi'\left(-y[i-1]\right)\Phi'\left(-y[i]\right)
-\Psi'\left(y[i-1]\right)\Phi'\left(y[i]\right) \right), &q \in [1, \frac{M}{3}]\\
&0,&q \in [\frac{M}{3}+1, M]     \\
\end{array}
\right..
\end{align}
By applying Lemma \ref{lem:func_property} -- i) [i.e., $\Psi(w)= \Psi'(w)= \Psi''(w)=0$ for $\forall \; w\le0$],  we immediately obtain that at least one of the terms $\Psi\left(-y[i-1]\right)\Phi''\left(-y[i]\right)$ or $ -\Psi\left(y[i-1]\right)\Phi''\left(y[i]\right)$ is zero. It follows that $$\Psi\left(-y[i-1]\right)\Phi''\left(-y[i]\right) -\Psi\left(y[i-1]\right)\Phi''\left(y[i]\right)   \le \sup_w|\Psi(w)|\sup_v|\Phi''(v)|.$$
Similarly, 
$$ \Psi''\left(-y[i]\right)\Phi\left(-y[i+1]\right)  -\Psi''\left(y[i]\right)\Phi\left(y[i+1]\right) \le \sup_w|\Psi''(w)|\sup_v|\Phi(v)| $$
$$ \Psi'\left(-y[i]\right)\Phi'\left(-y[i+1]\right)  -\Psi'\left(y[i]\right)\Phi'\left(y[i+1]\right) \le \sup_w|\Psi'(w)|\sup_v|\Phi'(v)|.$$

{Therefore, take the maximum over equations \eqref{eq:1a} to \eqref{eq:1c} and plug in the above inequalities, we obtain
\begin{align*}
\left|\frac{\partial^2 h_q}{ \partial y[i_1]\partial y[i_2]} \right| \nonumber
&\leq 3 \max \{\sup_w|\Psi''(w)|\sup_v|\Phi(v)|, \sup_w|\Psi(w)|\sup_v|\Phi''(v)|, \sup_w|\Psi'(w)|\sup_v|\Phi'(v)|\} \\
&= 3 \max \left\{8\pi, \frac{3\sqrt{3}}{2}, 4\sqrt{\frac{2}{e}}\right\} < 25\pi, \quad {\forall~i_1\; \mbox{being even}},\; \forall~i_2
\end{align*}
{where the equality comes from Lemma \ref{lem:func_property} -- ii).}

We can also verify that the above bound for $i$ being odd but not $1$ is exactly the same. 

When $i=1$ we have following:
	\begin{equation*}  
	\frac{\partial^2 h_q}{ \partial y[1]\partial y[1]}  = \left\{
	\begin{array}{lll}
	&- \Psi(1)\Phi''(y[1]) + 3\left( -\Psi''\left(-y[1]\right)\Phi\left(-y[2]\right) -\Psi''\left(y[1]\right)\Phi\left(y[2]\right) \right), &q \in [1, \frac{M}{3}]\\
	&- \Psi(1)\Phi''(y[1]),   &q \in [\frac{M}{3}+1, M] \\
	\end{array}
	\right.
	\end{equation*}
	\begin{equation*}  
	\frac{\partial^2 h_q}{ \partial y[1]\partial y[2]}  = \left\{
	\begin{array}{lll}
	&3\left( -\Psi'\left(-y[1]\right)\Phi'\left(-y[2]\right) -\Psi'\left(y[1]\right)\Phi'\left(y[2]\right) \right), &q \in [1, \frac{M}{3}]\\
	&0,   &q \in [\frac{M}{3}+1, M] \\
	\end{array}
	\right.
	\end{equation*}
 Again by applying Lemma \ref{lem:func_property}  -- i) and ii), 	
\begin{align*}
\left|\frac{\partial^2 h_q}{ \partial {y[1]}\partial y[i_2]} \right| \nonumber
&\leq  \max \{\sup_w|\Psi(1)\Phi''(w)|+3\sup_w|\Psi''(w)|\sup_v|\Phi(v)|, 3\sup_w|\Psi'(w)|\sup_v|\Phi'(v)|\} \\
&= \max \left\{ \frac{3\sqrt{3}}{2}(1-e^{-1})+24\pi , 12\sqrt{\frac{2}{e}}\right\}  < 25\pi,\;\forall~i_2.
\end{align*}

Summarizing the above results, we obtain:
\begin{align*} 
|g_{q}''\left(0; y,v\right)| &= 
|\sum_{\delta\in\left\{ 0,1, -1\right\}}
\sum_{i=1}^{T}
\frac{\partial^2}{\partial y[i]\partial y[i+\delta] }
h_q\left(y\right)v[i]v[i+\delta]|\\
&\le 25\pi \sum_{\delta\in\left\{ 0,1, -1\right\}}|
\sum_{i=1}^{T}
v[i]v[i+\delta]| \\
& =25\pi\left(  | \sum_{i=1}^{T}  v[i]^2| + 2| \sum_{i=1}^{T} v[i]v[i+1]|\right) \\
&\le 75\pi \sum_{i=1}^{T} | v[i]^2|=  75\pi.
\end{align*}
Overall, the first-order derivatives of $h_q$ are Lipschitz continuous for any $q$ with constant $\ell=75\pi$.

To show the same result for the function $\bar{h}$, we can apply \eqref{eq:Lip:original}. This completes the proof.  \QED

The following lemma is a simple extension of the previous result. 
\begin{lemma}\label{lem:property:f}
	We have the following properties for the functions $f$ and $\bar{f}$ defined in \eqref{eq:f:construct:2} and \eqref{eq:f:construct}. 
	
	\begin{enumerate}
		\item We have $\forall~x\in\mathbb{R}^{TM\times 1}$
		\begin{equation*}
		f(0) - \inf_x f(x) +\frac{1}{M U}\|d_0\|^2\leq \frac{1650\pi^2 \epsilon}{U} T,
		\end{equation*}
		where we have defined 
		\begin{align}{\label{eq:d0}}
		d_0 := [\nabla f_1(0), \cdots, \nabla f_M(0)].
		\end{align}
	
		\item We have
		\begin{align}
		\norm{\nabla \barf(y)}
		= \sqrt{2\epsilon}	\norm{\nabla \bar{h}\left(\frac{y U}{75 \pi \sqrt{2\epsilon}}\right)}, \;  \forall~y\in\mathbb{R}^{T\times 1}.
		\end{align}
		\item The first-order derivatives of ${\barf}$ and that for each $f_j, j\in[M]$ are Lipschitz continuous, with the same constant $U>0$. 
	\end{enumerate}
\end{lemma}

\noindent{\bf Proof.} To show that property 1) is true, note that from the definition of $f_i(x_i)$ we have 
$$\nabla f_i(x_i) = {\sqrt{2\epsilon}}\times \nabla h_i\left(\frac{x_i U}{75\pi \sqrt{2\epsilon}}\right).$$
Therefore the following holds:
\begin{align}
\frac{1}{M}\|d_0\|^2& =   \frac{2\epsilon}{M}\sum_{i=1}^{M} \|\nabla h_i(0)\|^2 \nonumber\\
&= \frac{2\epsilon}{M} \sum_{i=1}^{M} | \Psi(1)\Phi'(0)|^2 =  {32\epsilon} (1-\exp(-1))^2.
\end{align}
Therefore we have  the following:
{$$f(0) - \inf_x f(x) +\frac{\|d_0\|^2}{MU} = \frac{150 \pi \epsilon}{U} \left({h}(0) - \inf_x {h}(x) + \frac{16(1-\exp(-1))^2}{75\pi }\right).$$}

Then by applying Lemma \ref{lem:func_property} we have that for any $T\ge 1$, the following holds
{$$f(0) - \inf_x f(x) +\frac{\|d_0\|^2}{MU} \le \frac{150 \pi \epsilon}{U} \times (10\pi T + 0.03) \le \frac{150 \pi \epsilon}{U} \times 11 \pi T.$$

Property 2) is true due to the definition of $\bar{f}$. 

Property 3) is true because the following
$$\|\nabla \barf(z) - \nabla \barf(y)\| 
= \sqrt{2\epsilon} \left\| \nabla \barh\left(\frac{z U}{75 \pi \sqrt{2\epsilon}}\right) - \nabla \barh \left(\frac{y U}{75 \pi \sqrt{2\epsilon}}\right)\right\|
\le U\|z-y\|$$
where the last inequality comes from Lemma \ref{lem:func_property}  --  (5). This completes the proof.  \QED

Next let us analyze the size of $\nabla \bar{h}$. We have the following result. 
\begin{lemma} \label{lem:bound_gradient}
	If there exists $k\in[T]$ such that $|y[k]| < 1$,  then
	\begin{equation*}
	\norm{\nabla \barh(y)} = \norm{  \frac{1}{M}\sum_{i=1}^{M} \nabla h_i(y)}
	\ge \left|   \frac{1}{M}\sum_{i=1}^{M}    \frac{\partial }{\partial y[k]} h_i(y) \right| > 1.
	\end{equation*}
\end{lemma}
\noindent{\bf Proof.}
The first inequality holds for all $k\in[T]$, since {$\frac{1}{M}\sum_{i=1}^{M}    \frac{\partial }{\partial y[k]} h_i(y)$ is one element of $ \frac{1}{M}\sum_{i=1}^{M} \nabla h_i(y)$}.

We divide the proof for second inequality into two cases. 

\noindent{\bf Case 1.} Suppose $|y[j-1]|<1$ for all $2 \le j \le k$.  {Therefore, we have $|y[1]|<1$.}
Using \eqref{derivative_zero},  we have the following inequalities:
\begin{equation}  
\frac{\partial }{\partial y[1]} h_i(y) \stackrel{\rm (i)} \leq -\Psi(1)\Phi'(y[1]) \stackrel{\rm (ii)}<-1, \forall i
\end{equation}
{where ${\rm (i)}$ is true because $\Psi'(w), \Phi(w)$ are all non-negative from Lemma \ref{lem:func_property} -(2); ${\rm (ii)}$  is true due to Lemma \ref{lem:func_property}  --  (3). 
Therefore, we have the following
\begin{equation*}
\norm{\nabla \barh(y)} = \norm{  \frac{1}{M}\sum_{i=1}^{M} \nabla h_i(y)}
\ge \left|   \frac{1}{M}\sum_{i=1}^{M}    \frac{\partial }{\partial y[1]} h_i(y) \right| > 1.
\end{equation*} 
\noindent{\bf Case 2)} Suppose there exists $2 \le j \le k$  such that $|y[j-1]|\ge 1$.   

We choose $j$ so that $|y[j-1]|\ge 1$ and $|y[j]| < 1$. Therefore, depending on the choices of $(i,j)$ we have three cases

	\begin{equation*}   
	\frac{\partial h_i(y)}{\partial y[j]} = \left\{
	\begin{array}{lll}
	&-3 \left(  \Psi\left(-y[i-1]\right)\Phi'\left(-y[j]\right) +\Psi\left(y[i-1]\right)\Phi'\left(y[j]\right) \right), &i \in [1, \frac{M}{3}]\\
	&0,   &i \in [\frac{M}{3}+1, \frac{2M}{3}] \\
	& -3\left( \Psi'\left(-y[j]\right)\Phi\left(-y[i+1]\right) +\Psi'\left(y[j]\right)\Phi\left(y[i+1]\right) \right), &i \in [\frac{2M}{3}+1, M]
	\end{array}
	\right..
	\end{equation*}
	
If $i \in [1, \frac{M}{3}]$, because $|y[j-1]|\ge 1$ and $|y[j]| < 1$, using Lemma \ref{lem:func_property}  -- (3), and the fact that the negative part is zero for $\Psi$, and $\Phi'$ is even function,  the expression further equals to 
\begin{align}
- 3 \cdot \Psi(|y[j-1]|)\Phi'\left(|y[j]|\right)  ]\stackrel{\eqref{eq:key:3}}< -3,
\end{align}

If $i \in [\frac{2M}{3}+1, M]$ the expression is obviously non-positive because both $\Psi'$ and $\Phi$ are nonnegative. 
Overall, we have
$$   \left|   \frac{1}{M}\sum_{i=1}^{M}   \frac{\partial h_i(y)}{\partial y[j]} \right| > \left| \frac{1}{M} \sum_{i=1}^{M/3} 3 \right| = 1.$$ 
This completes the proof.  \QED
\begin{lemma}\label{lem:bound_gradient:new}
	Define $\bar{x} := \frac{1}{M} \sum_{i=1}^M x_i$, and  assume  that $U\in(0,1)$. Then we have 
	\begin{align*}
	\bigg\|\frac{1}{M}\sum_{i=1}^{M}\nabla f_i(x_i)\bigg\|^2 + \frac{U}{M \underline{\lambda}_{\min}(P^{1/2}\cL P^{1/2})} \sum_{(i,j):i\sim j}\|x_i-x_j\|^2
	% 	>\epsilon.
	\ge \frac{1}{2} \norm{\nabla \bar{f}(\bar{x})}^2.
	\end{align*} 
\end{lemma}
\noindent{\bf Proof.}
First let us derive a useful property. Define $d:= [d_1;d_2;\cdots; d_M]$ where $d_i$ is the degree for node $i$; further define 
$$\bar{x}:=\frac{1}{M}\sum_{i=1}^{M}x_i, \quad \tilde{x_i}:= x_i-\bar{x}, \quad \tilde{x}:=[\tilde{x}_1; \tilde{x}_2;\cdots; \tilde{x}_M]. $$
It is easy to observe that :
$$\tilde{x}^T \mathbbm{1} =0, \quad \mbox{and}\quad  \tilde{x} \notin \mbox{Null}(F^T F).$$ 
Then the following holds:
\begin{align}
x^T F^T F x= \sum_{(i,j):i\sim j}\|x_i -x_j\|^2 = \sum_{(i,j):i\sim j}\|\tilde{x}_i -\tilde{x}_j\|^2 = \tilde{x}^T F^T F\tilde{x}\ge \underline{\lambda}_{\min}(F^T F)\|\tilde{x}\|^2. 
\end{align}
Therefore the following holds:
\begin{align}\label{eq:new:difference}
\sum_{i=1}^{M}\|\bar{x} -x_i\|^2\le \frac{1}{\underline{\lambda}_{\min}(F^T F)} \sum_{(i,j):i\sim j}\|x_i-x_j\|^2 = \frac{1}{\underline{\lambda}_{\min}( P^{1/2}\cL P^{1/2})} \sum_{(i,j):i\sim j}\|x_i-x_j\|^2. 
\end{align}
Based on the above property, we have the following series of inequalities 
%{\small  
\begin{align*}
&\norm{\nabla \barf(\bar{x})}^2 
\le 2 \bigg\|\frac{1}{M}\sum_{i=1}^{M} \left(   \nabla f_i( \bar{x}) - \nabla f_i(x_i) \right)   \bigg\|^2  + 2 \bigg\|\frac{1}{M}\sum_{i=1}^{M}\nabla f_i(x_i)\bigg\|^2 \\
&\stackrel{\rm (i)} \le  \frac{2}{M}\sum_{i=1}^{M}  \bigg\|     \nabla f_i( \frac{1}{M} \sum_{j=1}^M{x_j}) - \nabla f_i(x_i)    \bigg\|^2  + 2 \bigg\|\frac{1}{M}\sum_{i=1}^{M}\nabla f_i(x_i)\bigg\|^2 \\
&\stackrel{\rm (ii)} \le  \frac{2}{M}\sum_{i=1}^{M} U^2 \bigg\|     \frac{1}{M} \sum_{j=1}^M{x_j} - x_i   \bigg\|^2  + 2 \bigg\|\frac{1}{M}\sum_{i=1}^{M}\nabla f_i(x_i)\bigg\|^2 \\
&\stackrel{\rm(iii)}\le  \frac{2 U}{M\underline{\lambda}_{\min}(P^{1/2}\cL P^{1/2})}\sum_{(i,j):i\sim j} \|     x_j - x_i   \|^2  + 2 \bigg\|\frac{1}{M}\sum_{i=1}^{M}\nabla f_i(x_i)\bigg\|^2,
\end{align*}
%} 
where in $\rm (i)$ and $\rm{(iii)}$ we have used the convexity of the function $\|\cdot\|^2$; {in $\rm (ii)$ we used Lemma \ref{lem:property:f}  --  (3)}; in ${\rm (iii)}$ we have also used the assumption that $U\in(0,1)$ and \eqref{eq:new:difference}. 
Overall we have
\begin{align*}
&\bigg\|\frac{1}{M}\sum_{i=1}^{M}\nabla f_i(x_i)\bigg\|^2 + \frac{U}{M \underline{\lambda}_{\min}(P^{1/2}\cL P^{1/2})} \sum_{(i,j):i\sim j} \|x_i-x_j\|^2
\ge \frac{1}{2} \norm{\nabla f(\bar{x})}^2.
\end{align*} 
This completes the proof. 
\QED

\begin{lemma}  \label{lem:bound_steps}
Consider using an algorithm in class $\cA$ or in class $\cA'$ to solve the following problem:
\begin{align}
\min_{x\in \mathbb{R}^{TM \times 1}} \; h(x)=\frac{1}{M}\sum_{i=1}^{M} h_i(x_i),
\end{align}
over a path graph. Assume  the  initial solution: $x_i= 0, \; \forall~i\in[M]$. Let $\bar{x}=\frac{1}{M}\sum_{i=1}^{M} x_i$ denote the average of the local variables.
Then the algorithm needs at least $(\frac{M}{3}+1) T$ iterations to have $x_i[T]\ne 0, \; \forall~i$ and $\bar{x}[T]\ne 0$.
\end{lemma}

{\bf Proof.} For a given $k\ge 2$, suppose that  $x_i[k], x_i[k+1],...,x_i[T]=0$, $\forall i$, that is, $\mbox{support}\{x_i\} \subseteq \{1,2,3,...,k-1\}$ for all $i$. Then  $\Psi'\left(x_i[k]\right) = \Psi'\left(-x_i[k]\right) = 0$ for all $i$, and $h_i$ has the following partial derivative when $k$ is even:
\begin{equation}  \label{eq:even:construction}
\frac{\partial h_i(x_i)}{\partial x_i[k]} = \left\{
\begin{array}{lll}
&- 3   \left( \Psi\left(-x_i[k-1]\right)\Phi'\left(-x_i[k]\right)\right)
+3   \left(\Psi\left(x_i[k-1]\right)\Phi'\left(x_i[k]\right)\right) , &i \in [1, \frac{M}{3}] \\
&0, &i \in [\frac{M}{3}+1, M]     \\
\end{array}
\right.
\end{equation}
and  the following partial derivative when $k$ is odd and $k\ge 3$:
\begin{equation}  \label{eq:odd:construction}
\frac{\partial h_i(x_i)}{\partial x_i[k]} = \left\{
\begin{array}{lll}
&0, &i \in [1, \frac{2M}{3}]     \\
&- 3   \left( \Psi\left(-x_i[k-1]\right)\Phi'\left(-x_i[k]\right)\right)
+3   \left(\Psi\left(x_i[k-1]\right)\Phi'\left(x_i[k]\right)\right) , &i \in [\frac{2M}{3}+1, M]  \\
\end{array}
\right..
\end{equation}

Recall that for any algorithm in class $\cA$ or $\cA'$, each agent is only able to compute linear combination of historical gradient and neighboring iterates [cf. \eqref{eq:span} and \eqref{eq:span:2}]. Therefore, for a given node $i$, as long as the $k$th element of the gradient as well as that of its neighbors have never been updated once, $x_i[k]$ remains to be zero. Combining this observation with the above two expressions for $\frac{\partial h_i(x_i)}{\partial x_i[k]}$, we can conclude that when $\mbox{support}\{x_i\} \subseteq \{1,2,3,...,k-1\}$ for all $i$, then in the next iteration $x_i[k]$ will be possibly non-zero on the node $i \in [1, \frac{M}{3}]$ for even $k$   and $i \in [\frac{2M}{3}+1, M]$ for odd $k$, and all other nodes still have $x_j[k]=0$, $\forall~j\ne i$.
		
{Now suppose that the initial solution is $x_i[k]=0$ for all $(i,k)$. Then at the first iteration only $\frac{\partial h_i (x_i) }{\partial x_i[1]}$ is non-zero for all $i$, 
	due to the fact that $\frac{\partial h_i (x_i) }{\partial x_i[1]} =  \Psi(1)\Phi'(0) =4(1-e^{-1}) $ for all $i$ from \eqref{derivative_zero}. 
	If follows that even if every node is able to compute 	its local gradient, and can communicate with their neighbors, it is only possible to have $x_i[1]\ne 0, \forall i$. {At the second iteration, we can use \eqref{eq:even:construction} to conclude that it is only possible to have $\frac{\partial h_r(x_r)}{\partial x_r[k]}\ne 0$ for some $r\in[1, \; M/3]$, therefore when using an algorithm in class $\cA$, we can conclude that $x_i[2]=0$ for all $i\notin[1,\; M/3]$.}

	Then following our construction \eqref{eq:construct:h}, we know the nodes in the set $[1, \frac{M}{3}]$ and the set $[\frac{2M}{3}+1, M]$ have minimum distance ${M}/{3}$. It follows that using an algorithm in $\cA$ or $\cA'$, it takes at least $M/3$ iterations for the non-zero {$x_r[2]$} and the corresponding gradient vector to propagate to at least one node in set $[2M/3+1, M]$.  {Once we have $x_j[2]\ne 0$ for some $ j\in[2M/3+1, \; M]$, then according to \eqref{eq:odd:construction},  it is possible to have $\frac{\partial h_j (x_j)}{\partial x_j[3]}\ne 0$, and once  this gradient becomes non-zero, the corresponding variable $x_{j}[3], j \in [2M/3+1, \; M]$ can become nonzero. }
	
	Following the above procedure, it is clear that we need at least $\frac{MT}{3}$ iterates and $T$ computations to make $x_i[T]$ possibly non-zero.
	\QED}

\begin{theorem}\label{thm:1}
	{Let $U\in(0,1)$ and $\epsilon$ be positive}. Let {$x^0[i]=0$} for all $i\in[M]$. Then for any distributed first-order  algorithm in class $\cA$ or $\cA'$,  there exists a problem in class $\cP^{M}_U$ and a network in class $\cN$, such that it requires at least the following number of iterations
	\begin{equation}
	t	\ge 	{\frac{1}{3\sqrt{\xi(\cG)}}}  \floor{\frac{\left( f(0) - \inf_x f(x) + \frac{\|d_0\|^2}{M U}\right)  U }{1650\pi^2 }\label{eq:lower:bound}
		\epsilon^{-1}}
	\end{equation}
	to achieve the following error
	\begin{align}
h^*_t= \bigg\|\frac{1}{M}\sum_{i=1}^{M}\nabla f_i(x^t_i)\bigg\|^2 + \frac{U}{M  \underline{\lambda}_{\min}(P^{1/2} \cL P^{1/2})}\sum_{(i,j):i\sim  j}\|x^t_i-x^t_j\|^2<\epsilon.
	\end{align} 
\end{theorem}

{\bf Proof.} 
By Lemma~\ref{lem:bound_steps} we have $\bar{x}[T] = 0$ for all $t<\frac{M+3}{3} T$. Then by applying Lemma~\ref{lem:property:f}  --  (2) and  Lemma~\ref{lem:bound_gradient}, we can conclude that the following holds
\begin{align}
\norm{\nabla \bar{f}(\bar{x}{[T]})}
= \sqrt{2\epsilon}	\norm{\nabla {\bar{h}}\left(\frac{\bar{x}{[T]} U}{75 \pi \sqrt{2\epsilon}}\right)} > \sqrt{2\epsilon},
\end{align}
{where the second inequality follows that there exists $k\in [T]$ such that $|\frac{\bar{x}[k] U}{75 \pi \sqrt{2\epsilon}}|=0< 1$, then we can directly apply Lemma~\ref{lem:bound_gradient}.} Then by applying Lemma~\ref{lem:bound_gradient:new} gives  $h^*_{(M+3)T/3}>\epsilon$.

The third part of  Lemma~\ref{lem:property:f} ensures  that $f_i$'s are $U$-Lipschitz continuous gradient, and the first part shows
\begin{equation*}
f(0) - \inf_x f(x)  + \frac{\|d_0\|^2}{MU}\leq \frac{1650\pi^2 \epsilon}{U} T,
\end{equation*}
Therefore we obtain
%$f(0) - \inf_x f(x) \le \Delta$, that is, 
\begin{equation}\label{lowerbound_T}
T \ge
\floor{\frac{\left( f(0) - \inf_x f(x) + \frac{\|d_0\|^2}{M U}\right)  U }{1650\pi^2}
	\epsilon^{-1}}.
\end{equation}
Summarizing the above argument, we have
\begin{equation*}	
t \ge \frac{M+3}{3} T \ge
\frac{M+3}{3} 
\floor{\frac{\left( f(0) - \inf_x f(x)+ \frac{\|d_0\|^2}{M U} \right)  U }{1650\pi^2 }
	\epsilon^{-1}}.
\end{equation*}
{By noting that for path graph $\xi(\cG)\ge 1/M^2$, this completes the proof. }
\QED

\subsection{Generalization}

The previous section analyzes the lower complexity bounds for problem $\cP^{M}_{U}$ over a path network. The obtained results can be extended in a number of direction.

\subsubsection{Uniform $L_i$, Fixed $D$ and $M$}
In this subsection, we would like to generalize Theorem \ref{thm:1} to a slightly wider class of networks (beyond the path graph used in our construction). Towards this end, consider a {\it path-star graph} shown in Fig. \ref{fig:path-star}. The graph contains a path graph with $D-1$ nodes, and the remaining nodes are divided into $D-1$ groups, each with either $\lfloor M/(D-1)-1 \rfloor$  or $\lfloor M/(D-1)-1 \rfloor+1$ nodes, and each group is connected to the nodes in the path graph by using a star topology. We have the following corollary to Theorem \ref{thm:1}. 
{%First, for any given $D$ and $M$ satisfy $D<M-1$, we have the following corollary. 
	\begin{corollary}\label{coro1}
		Let $U\in(0,1)$ and $\epsilon$ be positive, and fix any $D$ and $M$ such that $D\le M-1$. For any algorithm in class $\cA$ or $\cA'$, there exists a problem in class $\cP^{M}_U$ and a  network in class $\cN^{M}_D$, so that to achieve accuracy  $h^*_t<\epsilon$, it requires at least the following number iterations
		\begin{equation*}
			t	\ge 	\frac{D}{3}  \floor{\frac{\left( f(0) - \inf_x f(x) + \frac{\|d_0\|^2}{MU}\right)  U }{1650\pi^2}
				\epsilon^{-1}}.
		\end{equation*}
		Alternatively, the above bound can be expressed as the following
				\begin{equation*}
		t	\ge 	\frac{\sqrt{(D-1)/(2M)}}{3 \sqrt{\xi(\cG)}}  \floor{\frac{\left( f(0) - \inf_x f(x) + \frac{\|d_0\|^2}{MU}\right)  U }{1650\pi^2}
			\epsilon^{-1}}.
		\end{equation*}
	\end{corollary}
	
	{\bf Proof.}  Fix any $D$ and $M$ such that $D\le M-1$, we can construct a path-star graph as described in Fig.\ref{fig:path-star}, whose diameter is $D$.  %Starting from the two ends, each group are connected to one of the $D-1$ nodes using a star topology. 
	\begin{figure}
		\centering
		\includegraphics[width=0.6\linewidth]{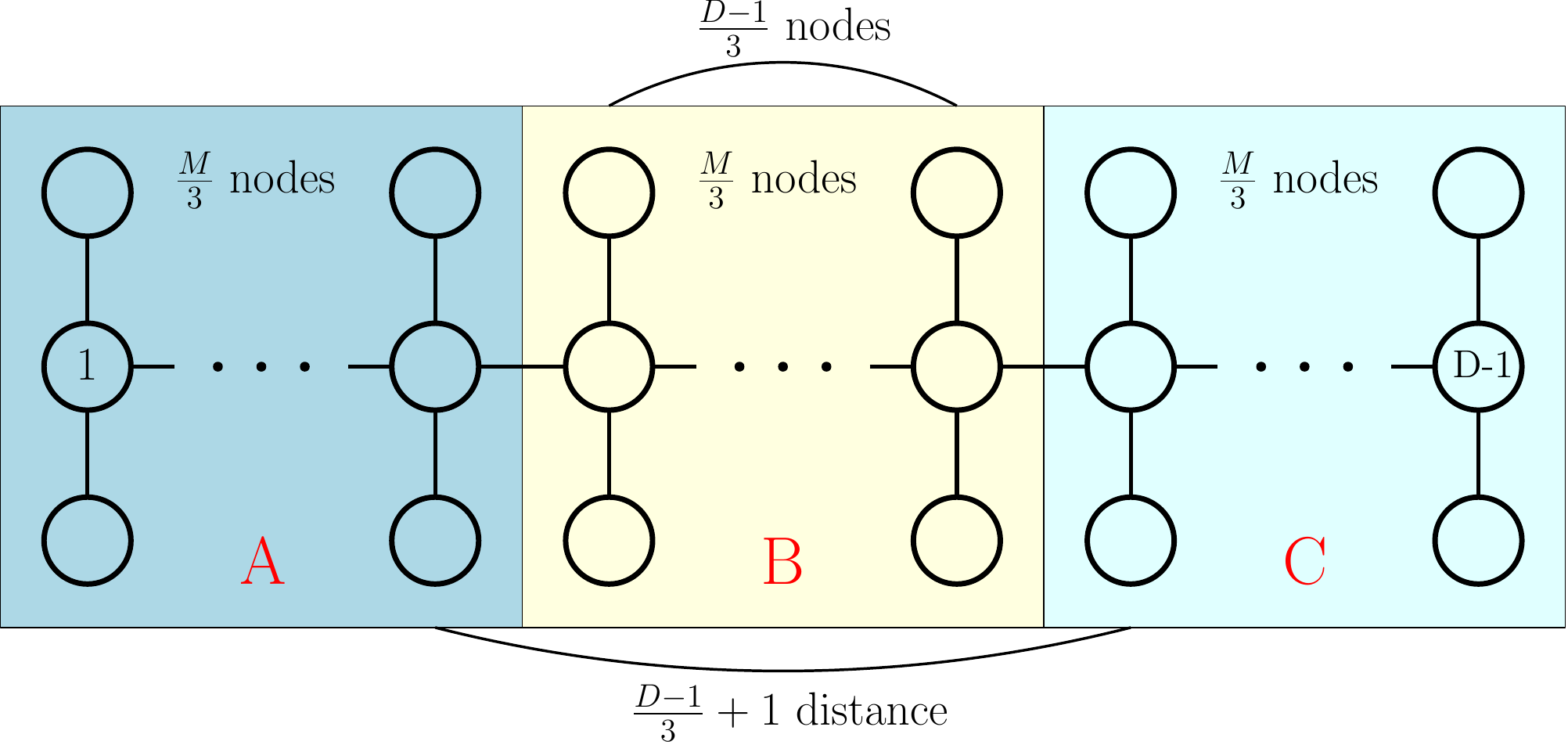}
		\caption{The path-star graph used in our construction.}
		\label{fig:path-star}
	\end{figure}
	
	To show the lower bounds for such a graph, we split all $M$ nodes into three sets $\mathcal{A}, \mathcal{B}, \mathcal{C}$ based on the main path, each with $\frac{M}{3}$ nodes (assume $M$ is a multiple of 3), where $\mathcal{A}$ and $\mathcal{C}$ has minimum $\frac{D+2}{3}$ distance (assume $D-1$ is a multiple of 3). Then we construct the component functions $h_i$'s as follows. 
	
	\begin{equation} 
	h_i(x_i) = \left\{
	\begin{aligned}{}
	&\Theta(x_i, 1)+3\sum_{j=1}^{\floor{T/2}} \Theta(x_i, 2j), &i\in \mathcal{A}  \\
	&\Theta(x_i, 1), &i \in \mathcal{B}  \\
	&\Theta(x_i, 1)+3\sum_{j=1}^{\floor{T/2}} \Theta(x_i, 2j+1), &i \in \mathcal{C}\\
	\end{aligned}
	\right.
	\end{equation}
	
	Since the graph has diameter $D$ in the above construction, and the distance between any two elements in $\mathcal{A}$ and $\mathcal{C}$  is at least $\frac{D+2}{3}$ (assume $D-1$ is a multiple of 3), by a similar step in Lemma~\ref{lem:bound_steps} we can conclude that we need at least $(\frac{D+2}{3}+1) T$ iterations to achieve $x_i[T]\ne 0$. By applying \eqref{lowerbound_T}, we can obtain the desired result.
	
	To show the second result, note that from \eqref{eq:path:star:lambda} we have
	\begin{align}
	\sum_{i}d_i D \ge \frac{1}{\underline{\lambda}_{\min}(\cL)}
	\end{align}
	For the path-star graph under consideration,
		$$\sum_i d_i \le 2 (D-1)-2 + 2 \left( M-(D-1)\right)  \le 2M,$$
		so the following holds: 
		$$D^2 \ge \frac{D/2M}{\underline{\lambda}_{\min}(\cL) } \ge \frac{(D-1)/(2M)}{\underline{\lambda}_{\min}(\cL)}.$$
	The desired result is then immediate. 
	\QED

\subsubsection{Non-uniform $L_i$, Fixed $\cN$}

Finally, for the problem class with non-uniform Lipschitz constants, we can extend the previous result to any network in class $\cN$ (by properly assigning different values of $L_i$'s to different nodes). {In this case the lower bound will be dependent on the spectrum property of $\hL$ as defined in \eqref{eq:hat:L} (expressed below for easy reference)}
\begin{align}
{\hL} := L^{-1/2} F^T K F L^{-1/2}. 
\end{align}

\begin{corollary}\label{thm:2}
	Let $\epsilon$ be positive. For any given network in $\cN^{M}_D$, and for any algorithm in $\cA$, there exists a problem in $\cP^{M}_L$ such that to achieve accuracy  $h^*_t<\epsilon$,  it requires at least the following iterations
	\begin{equation}\label{eq:lower:bound:nonuniform}
	t	\ge 	\frac{1}{3\sqrt{\xi(\hL)}}  \floor{\frac{\left( f(0) - \inf_x f(x) + \frac{\|d_0\|_{L^{-1}}^2}{M}\right)  {\bar{L}}}{1650\pi^2 }
		\epsilon^{-1}}.
	\end{equation}
\end{corollary}

To prove this result, we select the values of the coefficient set $\{L_i\}_{i=1}^{M}$, so that the ``effective" network topology becomes a path. In particular, for any given network in $\cN$, we can construct local functions as follows: First, along the longest path of size $D$, we distributed the functions into three sets $\mathcal{A},  \mathcal{B},  \mathcal{C}$, where $\mathcal{A}$ and $\mathcal{C}$ denotes the first and  last $\frac{D}{3}$ nodes on the path respectively, and $\mathcal{B}$ denotes the  rest nodes on the path. Second, for the rest of the functions not on the path, denoted as set $\mathcal{D}$, set their local functions to zero (or equivalently, set the corresponding $L_i$'s to zero). Then, the local function belongs to each set can be expressed as:  

\begin{equation} 
h_i(x_i) = \left\{
\begin{aligned}{}
&\frac{M}{D} \Theta(x_i, 1)+\frac{3M}{D}\sum_{j=1}^{\floor{T/2}} \Theta(x_i, 2j), &i\in \mathcal{A}  \\
&\frac{M}{D}\Theta(x_i, 1), &i \in \mathcal{B}  \\
&\frac{M}{D}\Theta(x_i, 1)+\frac{3M}{D}\sum_{j=1}^{\floor{T/2}} \Theta(x_i, 2j+1), &i \in \mathcal{C}\\
&0, & i \in \mathcal{D}
\end{aligned}
\right.
\end{equation}
This way the network reduces to a path graph. Note that the Lipschitz constant for the gradient of $h(y) = \frac{1}{M}\sum_{i=1}^{M}h_i(y)$ is  still $1$, and we can use the similar constructions and proof steps leading to Theorem \ref{thm:1} to prove the claim.

\section{The Proposed Algorithms}
In this section, we introduce two different types of algorithms for solving problem \eqref{eq:global:consensus:equiv:2}. The algorithm is {\it near-optimal}, and can achieve the lower bounds derived in Section \ref{sec:lower} except for a multiplicative polylog factor in $M$. To simplify the notation, we utilize the definitions introduced in Section \ref{sub:prelim}, and rewrite problem \eqref{eq:global:consensus:equiv:2} in the following compact form
\begin{align}\label{eq:global:consensus:equiv}
\min_{x\in\mathbb{R}^{SM}} \; f(x):=\frac{1}{M}\sum_{i=1}^{M} f_i(x_i),\quad \st\; (F \otimes I_S) x =0.
\end{align}
It can be verified that, by  using the definition of $F$, the constraint in this problem is equivalent to the ones given in \eqref{eq:global:consensus:equiv:2}. For notational simplicity, in the following we will assume that $S=1$ (scalar variables). All the results presented in subsequent sections extend easily to case with $S>1$. 

\subsection{The D-GPDA Algorithm}
We first present a Distributed Gradient Primal-Dual Algorithm (D-GPDA), which relaxes the linear constraint \eqref{eq:global:consensus:equiv}, and gradually enforces it as the algorithm proceeds.  
To describe the algorithm, let us introduce the augmented Lagrangian (AL) function as 
\begin{align}
\hspace{-0.3cm}{\sf{AL}}(x,\lambda) = f(x)+  \langle \lambda, Fx\rangle +\frac{1}{2}\|\Sigma F x\|^2\label{eq:augmented},
\end{align}
where $\lambda\in\mathbb{R}^{E}$ is the dual variable; $\Sigma = \mbox{diag}([\sigma_1,\cdots, \sigma_{E}])\in\mathbb{R}^{E\times E}$ is a diagonal positive definite matrix. In the following,  we will use the shorthanded notation ${\sf{AL}}^r:= {\sf{AL}}(x^r,\lambda^r)$ where $r$ is the iteration counter.  

Define a {\it penalty matrix} as 
\begin{align}
\Upsilon = \mbox{diag}\{[\beta_1, \cdots, \beta_M]\}\succ 0. 
\end{align}
Then the D-GPDA is described in the following table.
\begin{center}
%	\vspace{-0.2cm}
	\fbox{
		\begin{minipage}{0.8\linewidth}
			\smallskip
			\centerline{\bf {Algorithm 1. The D-GPDA Algorithm}}
			\smallskip
			
			{\bf (S1).} Assign each node $i\in\cN$ with a parameter $\beta_i >0$; Assign each edge $(ij)\in\cE$ with a parameter $\sigma_{ij}>0$;

			{\bf (S2).}  At iteration $r=-1$, initialize $\lambda^{-1} = 0$ and $x^{-1} =0$;
			
			{\bf (S3).}  At iteration $r=0$, set 	$\lambda^0$ and $x^0$ using the following:
			\begin{align}
			%\frac{1}{M}\sum_{i=1}^{M}\nabla f_i (0) -A^T \lambda^0 =0. 
		 	\nabla f(x^{-1}) + (2\Delta + \Upsilon^2) x^0 = 0, \; 	\lambda^0 = \Sigma^2 F x^0; \quad \label{eq:initialization}
			\end{align}
			Equivalently $x^0$ can be written as:
			\begin{align}\label{eq:init:x0}
			x^0_i  = \bigg(2 \sum_{j: j\sim i}\sigma^2_{ij} + \beta^2_i\bigg)^{-1}\nabla f_i(0)/M,\;  \forall~i\in[M];
			\end{align}
			
			{\bf (S4).}  At each iteration $r+1$, $r\ge 0$, update variables  by:
			\begin{subequations}
				\begin{align}
				x^{r+1}& =\arg\min_{x}\; \langle \nabla f(x^r) + F^T\lambda^r , x-x^r\rangle \label{eq:x:update}\\
				&\quad +\frac{1}{2}\|\Sigma F x\|^2 +\frac{1}{2}\|\Sigma B (x-x^r)\|^2+ \frac{1}{2}\|\Upsilon (x-x^r)\|^2\nonumber\\
				\lambda^{r+1}& = \lambda^r +\Sigma^2 F x^{r+1} \label{eq:mu:update}.
				\end{align}
			\end{subequations}
			
		\end{minipage}
	}
\end{center}

We note that each iteration of the D-GPDA performs a gradient descent type step on the AL function, followed by taking a step of dual gradient ascent (with a {\it stepsize matrix} $\Sigma^2\succ 0$). The term $\frac{1}{2}\|\Sigma B(x-x^r)\|^2$ used in \eqref{eq:x:update} is a {\it network proximal term} that regularizes the $x$ update using network structure, and its presence is critical to ensure separability and distributed implementation (see Remark \ref{rmk:distributed} below). 

The D-GPDA is closely related to many classical primal-dual methods, such as the Uzawa method \cite{UZAWA58} (which has been recently utilized to solve linearly constrained {\it convex} problems \cite{Nedic2009saddle}), and the proximal method of multipliers (prox-MM)  \cite{rockafellar1976augmented,wright_proximal}. The latter method has been first developed by Rockafellar in \cite{rockafellar1976augmented}, in which a proximal term has been added to the AL in order to make it strongly convex  in each iteration. However, the theoretical results derived for Prox-MM in \cite{rockafellar1976augmented, wright_proximal} are only valid for convex problems. It is also important to note that when the matrices $\Sigma$ and $\Upsilon$ are specialized as multiples of identity matrices, that is, when $\Sigma =\sigma I_M$ and $\Upsilon=\kappa I_M$ for some $\sigma, \kappa>0$, then the D-GPDA reduces to the Prox-GPDA algorithm briefly discussed in our earlier work \cite[Section 5]{hong17icml}, for solving a general linearly constrained problem.

\subsection{The xFILTER Algorithm}
Despite the fact that D-GPDA is conceptually simple, we will show shortly that it is only optimal for special network classes with small diameter [or large gap function $\xi(\cG)$], such as the complete/star networks (see Table \ref{table:result} and our detailed analysis in Section \ref{sec:tightness:uniform}).  Intuitively, the issue is that having the {\it network proximal term} imposes very heavy regularization, enforcing the new iterates to be close to the old ones. This causes slow information propagation over the network.

In this section, we present a {\it near-optimal} algorithm that can achieve the lower bounds derived in Section \ref{sec:lower} for a number of different graphs (up to some polylog factor in the problem dimension).
{To motivate our algorithm design, observe that the communication lower bound $\mathcal{O}(1/\sqrt{\xi(\mathcal{G})} \times \bar{L} /{\epsilon})$ in Section \ref{sec:lower} can be decomposed into the product two parts,  $\mathcal{O}(1/\sqrt{\xi(\mathcal{G})})$  and $\mathcal{O}( \bar{L} /{\epsilon})$,  corresponding roughly to the communication efficiency and the computational complexity, respectively. Such a product form motivates us to {\it separate} the computation and communication tasks, and  design a {\it double loop} algorithm to achieve the desired lower bound.
	
	Our proposed algorithm is based on a novel  {\it appro\underline{x}imate \underline{filt}ering -then- pr{\underline e}dict and t{\underline r}acking} (xFILTER) strategy, which properly  combines the modern first-order optimization methods and the classical polynomial filtering techniques. It is a ``double-loop" algorithm, where in the outer loop local gradients are computed to extract information from local functions, while in the inner loop some filtering techniques are used to facilitate efficient information propagation. Please see {\bf Algorithm 1}  for the detailed description, from the system perspective.  
	It is important to note that the algorithm contains an outer loop {\bf (S3)}--{\bf (S4)} and an inner loop {\bf (S2)}, indexed by $r$ and $q$, respectively. Further, the local gradient evaluation only appears in the outer loop step  {\bf (S3)}.
	
	To understand the algorithm, we note that one important task of each agent is to update its local variable so that it is close to the average $\frac{1}{M}\sum_{i=1}^{M} x_i$. Let us use $d_i$ to denote a local variable that {\it approximates} the above average. At the beginning of the algorithm, $d_i$ is just a rough estimate of the average, so we have $d_i=\frac{1}{M}\sum_j x_j + e_i$, where $e_i$ is the {\it deviation} from the true average, and it can be viewed as some kind of ``estimation noise". %Then it is natural that the proposed method must have some algorithmic components that can gradually remove such deviation (or ``noises"). %This is where the so-called ``graph based joint bilateral filtering" comes in. 
	To gradually remove such a noise, in step {\bf S1)} we resort to the so-called graph based joint bilateral filtering used for image denoising \cite{Tian14,Gadde13}, which can be formulated as the following regularized least squares problem:
	\begin{align}\label{eq:filter}
	x^{r+1}_{*}:=\arg\min _{x\in\mathbb{R}^{M}}\frac{1}{2}\|x- d^{r}\|_{\Upsilon^2}^2 + \frac{1}{2}x^{\top} F^{\top} \Sigma^2 F x,
	\end{align}	
	where $d^r$ is the noisy signal, $F$ is a penalty high pass filter related to the graph structure (in our case, $F$ is the adjacency matrix), and $\Sigma^2$ is a regularization parameter. 
	Its {solution}, denoted as $x_{*}^{r+1}$ as given below, will be close to the ``unfiltered" signal $d^r$, while having reduced high frequency components, or high fluctuations across the components:
	%(see \cite{Gadde13,Davies01,Smola03} for further discussion):
	\begin{align}\label{eq:opt:coupled}
	Rx^{r+1}_{*} = d^r, \quad \mbox{with}\quad R:=\Upsilon^{-2}F^{\top}\Sigma^2F + I_M. 
	\end{align}
	It is important to note that if $x_*^{r+1}$ indeed achieves consensus, then by \eqref{eq:null} we have $F^{\top}\Sigma^2F x_*^{r+1}=0$, implying $x_*^{r+1} =d^{r}$, which says $d^r$ should ``track"  $x^{r+1}_{*}$. 
	
	Unfortunately,  the system \eqref{eq:opt:coupled} cannot be precisely solved in a distributed manner, because inverting $R$ destroys its pattern about the network structure embedded in the product $F^{\top} \Sigma^2 F$. {More specifically, $F^{\top}\Sigma^2 F$ is the weighted graph Laplacian matrix whose $(i,j)$th entry is nonzero if and only if node $i$, $j$ are connected, but $(\Upsilon^{-2}F^{\top}\Sigma^2F + I_M)^{-1}$ is a dense matrix without such a property.}
	Therefore in {\bf S2)},  we use a degree-$Q$ Chebyshev polynomial to approximate $x^{r+1}_*$. The output, denoted as $x^{r+1}$, stays in a Krylov space ${span}\{d^r, R d^r, \cdots, R^Q d^r\}$. {Specifically, at each iteration, the only step that requires communication is the operation $R u$, which is given by
		\begin{align}\label{eq:Rd}
		(R u_{q-1}) [i] &= (\Upsilon^{-2}F^{\top}\Sigma^2F u_{q-1})[i] + d_{q-1}[i] \\
		&= \frac{1}{\beta^2_i} \sum_{j:j\sim i} \sigma^2_{ij} (u_{q-1}[j] - d^r[i])   + u_{q-1}[i], \; \forall~i,\nonumber
		\end{align}
		so this step can be done distributedly,  via one round of local message exchange.} 
	
	After completing $Q>0$ such Chebyshev iterations \eqref{eq:cheby:iteration} (C-iteration for short), the obtained solution $x^{r+1}$ will be an approximate solution to the system \ref{eq:opt:coupled}, with a residual error vector  $\epsilon^{r+1}$ as given below %denotes the residual error resulted %satisying The error vector can be expressed as below:
	\begin{align}\label{eq:linear:cheby:error}
	&R x^{r+1} = d^r + R \epsilon^{r+1},  \;\mbox{with}\quad \epsilon^{r+1}:=x^{r+1} - x^{r+1}_{*}.
	\end{align}}

 Up to this point, the filtering technique we have discussed aims at removing the ``non-consensus" parts from a vector $d=[d_1,\cdots, d_N]^T$. However, recall that the goal of distributed optimization is not only to achieve consensus, but also to optimize the objective function $\sum_i f_i(x_i)$. Therefore, a {\it prediction} step {\bf (S3)}  is performed to incorporate the most up-to-date local gradient $\nabla f_i(x_i)$. Then a {\it tracking} step {\bf (S4)} is performed to update $d$. Ideally, one would like the new $d^{r+1}_i$ to have the following three properties: 1) It is close to the previous $d^r_i$; 2) it takes into consideration the new local gradient information offered by the ``predicted" $\tilde{x}^{r+1}_i$; 3) it is a ``low frequency" signal, meaning $d^{r+1}_i$ and $d^{r+1}_j$ are relatively close, for all $i\ne j$.  Taking a closer look at the ``tracking" step, we can see that all three components are included: It adds to the previous $d^r$  the differences of the last two predictions, and it removes some non-consensus components among the local variables.
The detailed algorithm is given in the  Algorithm 2. 
%To facilitate understanding the proposed algorithm, we also include a table that explains the xFilter algorithm from the perspective of each agent (i.e., local view).

 To end this subsection, we emphasize that, the use of the polynomial Chebyshev filtering requires $Q$ vector communications steps every time  that {\bf (S2)} is performed. However, such a filtering step is critical to make the proposed algorithm achieve performance lower bounds predicted in Section \ref{sec:lower}.}  Intuitively, it helps to accelerate information propagation across the network. Indeed, as will be shown shortly, the number $Q$ in {\bf (S2)} is directly related to properties of the underlying graph. It is also somewhat surprising that the inner problem \eqref{eq:filter} is not required to be solved with {\it increased} accuracy. On the contrary, only a {\it fixed} number of filtering steps are needed.  
\begin{center} 
	%\vspace{-0.5cm}
	\fbox{ 
		\begin{minipage}{0.9\linewidth}
			\smallskip
			%(Global View)
			\centerline{\bf {Algorithm 2. The xFILTER Algorithm }}
			\smallskip
			
			{\bf (S1) [Initialization].} Assign each node $i\in\cN$ with  $\beta_i >0$; Assign each edge $(ij)\in\cE$ with $\sigma_{ij}>0$;  
			{Initialize $x^{-1}=0$, $d^{-1} = -\Upsilon^{-2}\nabla f(x^{-1})$ and $\tilde{x}^{-1} = x^{-1}-\Upsilon^{-2}\nabla f(x^{-1})$.} {Compute $R$  by \eqref{eq:opt:coupled};}
			
			{\bf (S2) [Filtering].}  At iteration $r+1$, $r\ge -1$: 
			For a fixed constant $Q>0$, run the following C-iterations (with parameters $\{\alpha_{q}, \tau\}$) %to solve 
			%the linear system of equations
			\begin{align}\label{eq:cheby:iteration}
			&u_0= x^r,  \; u_1 = (I - \tau R) u_0 + \tau d^r,\\
			& u_{q} = \alpha_{q}(I-\tau R) u_{q-1} + (1-\alpha_{q})u_{q-2}+\tau \alpha_{q} d^r,  \; q= 2,\cdots, Q\nonumber;
			\end{align}
			Set $x^{r+1} = u_Q$; 
			
			{\bf (S3) [Prediction].} Compute $\tilde{x}^{r+1}$ by:
			\begin{align}\label{eq:x:tilde:update}
			\tilde{x}^{r+1} = x^{r+1}- \Upsilon^{-2}\nabla f(x^{r+1});
			\end{align}
			
			{\bf (S4) [Tracking].}  Compute ${{d}^{r+1}}$ by:
			\begin{align}\label{eq:d:update}
			d^{r+1}  = d^{r}+  (\tilde{x}^{r+1} -\tilde{x}^{r})- \Upsilon^{-2}F^{\top} \Sigma^2 F  x^{r+1}.
			\end{align}
			Set $r=r+1$, go to {\bf (S2)}.
			%			\end{subequations}
		\end{minipage}
		%	\vspace{-0.1cm}
	}
\end{center}

\subsection{Discussion}\label{sub:discussion}
In this subsection, we establish some key connections between the two algorithms discussed so far, and provide some additional remarks. 

First, we provide an important interpretation of the xFILTER strategy, which will help us subsequently provide an unified analysis framework for both D-GPDA and xFILTER. 
First, similarly as in the D-GPDA algorithm, let us introduce an auxiliary variable $\lambda^r\in\mathbb{R}^{E}$, which is updated as follows:
\begin{align}\label{eq:mu:update:+}
\lambda^{r+1} = \lambda^r + \Sigma^2 F x^{r+1}. 
\end{align}
Suppose $\lambda^{-1}=0$, then according to \eqref{eq:d:update} and \eqref{eq:x:tilde:update} we have the following relationship
\begin{align*}
d^0 &:=- \Upsilon^{-2}\nabla f(x^{-1}) + ( x^{0}- \Upsilon^{-2}\nabla f(x^{0}) -( x^{-1}- \Upsilon^{-2}\nabla f(x^{-1})))-\Upsilon^{-2}F^T \lambda^0\nonumber\\
&=  x^{0}- \Upsilon^{-2}\nabla f(x^0)- \Upsilon^{-2}F^T \lambda^0.
\end{align*}
By using the induction argument,  we can show that for all $r\ge 0$, the following holds
\begin{align}\label{eq:R}
d^r &:= x^{r}-\Upsilon^{-2}\nabla f(x^r) - \Upsilon^{-2}F^T \lambda^r.
\end{align}
Combining  \eqref{eq:opt:coupled} and \eqref{eq:R}, we obtain the following useful alternative expressions of \eqref{eq:opt:coupled} and \eqref{eq:linear:cheby:error}: 
  \begin{subequations}
			\begin{align}
			&\Upsilon^{-2} \hspace{-0.1cm}\left(\nabla f(x^r) \hspace{-0.1cm} + \hspace{-0.1cm} F^{\top} (\lambda^r + \Sigma^2 F x^{r+1}_{*})\right) \hspace{-0.1cm} + \hspace{-0.1cm}(x^{r+1}_{*} -x^r)= 0\label{eq:linear:cheby:equiv}\\
			&\Upsilon^{-2} \hspace{-0.1cm}\left(\nabla f(x^r) \hspace{-0.1cm} + \hspace{-0.1cm} F^{\top} \hspace{-0.1cm}(\lambda^r + \Sigma^2 F x^{r+1})\right) \hspace{-0.1cm} + \hspace{-0.1cm}(x^{r+1}\hspace{-0.2cm} -x^r) = R\epsilon^{r+1}\label{eq:linear:cheby:error:equiv}.
			\end{align}
\end{subequations}
Using \eqref{eq:linear:cheby:equiv}, it is clear that $x^{r+1}_{*}$ can be equivalently written as the optimal solution of the following problem:
\begin{align}\label{eq:equivalent:problem}
x^{r+1}_{*}&=\argmin_x\; \langle \nabla f(x^r) + F^T\lambda^r , x-x^r\rangle +\frac{1}{2}\|\Sigma F x\|^2 +  \frac{1}{2}\|\Upsilon (x-x^r)\|^2. 
%&=\argmin_x\; \langle \nabla f(x^r) + F^T\lambda^r , x-x^r\rangle +\frac{1}{2}\|\Sigma F x\|^2 +  \frac{1}{2}\|\Upsilon (x-x^r)\|^2\nonumber\\
\end{align}
The relations \eqref{eq:mu:update:+} and \eqref{eq:equivalent:problem} together show that D-GPDA and xFILTER are closely related. However, we note that when comparing \eqref{eq:equivalent:problem} with \eqref{eq:x:update}, one key difference  is that the {\it network proximal} term $\frac{1}{2}\|\Sigma B(x-x^r)\|^2$ used in D-GPDA is no longer used in xFILTER.

We have the following additional remarks on the proposed algorithms. 
\noindent\begin{remark}{\bf{(Parameters)}}
 It is important to note that in both Alg. 1 and 2, in the update of the primal and dual variables, some {\it ``matrix parameters"} are used instead of scalar ones. In particular, the matrix $\Upsilon^2$ is used as the primal {\it ``proximal parameter"}, while $\Sigma^2$ is used as the {\it ``dual stepsize"}. Using these matrices ensures
that we can appropriately design parameters for each node/link, which 
is one key ingredient in ensuring the optimal rate. %Using these matrices provides flexibility to allow each node and link has its own individual parameter. 
\end{remark}
\begin{remark} {\bf (Initialization)}
	The initialization steps in {\bf (S2)} and {\bf (S3)} of Alg. 1  
	can be done in a distributed manner. Each node $i$ only requires to know the neighbors' $\sigma^2_{ij}$'s in order to update $x^0_i$. Once $x^0$ is updated, $\lambda^0$ can be updated by using:
	\begin{align*}
	\lambda^0_{ij} = \sigma^2_{ij} (x^0_i - x^0_j), \quad\forall~(i,j)\in E.
	\end{align*}
\end{remark}

\begin{remark} {\bf (Distributed Implementation and Algorithm Classes)}\label{rmk:distributed}
To see how the D-GPDA can be executed distributedly, we write down the optimality condition of  \eqref{eq:x:update}. 
For notational simplicity, define: 
\begin{align}\label{eq:H}
H:=B^T \Sigma^2 B + \Upsilon^2.
\end{align}
Then we have{
%\begin{subequations}
	\begin{align}\label{eq:optimality}
	\hspace{-0.3cm}\nabla f(x^r) + F^T \lambda^r +  F^T \Sigma^2 F x^{r+1}+ H (x^{r+1}-x^r) = 0.
	\end{align}}
%\end{subequations}
Rearranging, and using property  \eqref{eq:sum:A:B}, we have
\begin{align*}
\nabla f(x^r) +F^T \lambda^r + (2 \Delta+ \Upsilon^2) x^{r+1} - H x^r=0.
\end{align*}
Subtracting the same equation from the $r$th iteration, 
and use the fact that $F^T (\lambda^{r}-\lambda^{r-1})= F^T \Sigma^2  F x^r$, we have
\begin{align}\label{eq:EXTRA}
x^{r+1} &= x^r -\left(2\Delta+\Upsilon^2\right)^{-1} \bigg(\nabla f(x^r) - \nabla f(x^{r-1}) + (F^T \Sigma^2 F-H)x^r +Hx^{r-1}\bigg).
\end{align}
According to the above update rule, each node $i$ can distributedly implement \eqref{eq:EXTRA} by performing the following{
\begin{align}\label{eq:EXTRA:distributed}
&\hspace{-0.4cm}x_i^{r+1} = x_i^r -\frac{1}{2\sum_{j:j\sim i}\sigma^2_{ij} + \beta^2_i} \bigg(\frac{1}{M}(\nabla f_i(x_i^r) - \nabla f_i(x_i^{r-1}))\\
&\;  - 2\sum_{j: j\sim i}\sigma^2_{ij} x_j^r - \beta^2_i x^r_i + %\sum_{j:j\sim i}(\sigma^2_{ij}/d_i) x_j^{r-1} + 
\beta^2_i x^{r-1}_i + \sum_{j:j\sim i}\sigma^2_{ij}(x^{r-1}_j+x_i^{r-1})\bigg).\nonumber
\end{align}}

It is also easy to see that the Chebyshev iteration in xFILTER can be implemented distributedly, since the $R$ matrix defined in \eqref{eq:opt:coupled} preserves the network structure. To see how we can compute the $d^r$ vector distributedly, we first note that $d^{-1} = -\Upsilon^{-2}\nabla f(0)$. Then suppose we know $d^{r-1}$, by combining \eqref{eq:d:update} and \eqref{eq:x:tilde:update} we have
	\begin{align}
	d^{r} =d^{r-1}+  (x^{r} -x^{r-1})- \Upsilon^{-2}(\nabla f(x^{r}) -\nabla f(x^{r-1})) - \Upsilon^{-2} F^T\Sigma^2 F x^{r}\nonumber.
	\end{align}

Therefore each $d^r_i$ can be updated as
\begin{align}\label{eq:d:distributed}
	d^r_i = d^{r-1}_i +  (x^r_i-x^{r-1}_i) - \frac{1}{M\beta^{2}_i}(\nabla f_i(x^r_i)-\nabla f_i(x^{r-1}_i)) + \sum_{j: j\sim i}\frac{\sigma^2_{ij}}{\beta^2_i} (x^r_i - x^r_j).
	\end{align}
Combining the above expression with the expression in \eqref{eq:Rd} for computing $Rd^r$, it is clear that all the computation only involves in local communication and local gradient computation.

These observations also suggest that for a general choice of parameter matrix $\Sigma^2\succ 0$,  both D-GPDA and xFILTER are in class $\cA$. Further, if $\Sigma^2$ is a multiple of identity matrix (i.e., there exists $\sigma^2>0$ such that $\Sigma^2 = \sigma^2 I_E$), then the computations in \eqref{eq:EXTRA:distributed} and \eqref{eq:d:distributed} only involve the sum of neighboring iterates, therefore both algorithms belong to class $\cA'$ as well.  
\end{remark}

%\vspace{-0.4cm}
\section{The Convergence Rate Analysis}
In this section we provide the analysis steps of the convergence rate of the D-GPDA and xFILTER. All the proofs of the results can be found in the appendix. Note that we use the primal-dual representation discussed in Section \ref{sub:discussion} for xFILTER, so that it can be analyzed together with the D-GPDA. 

\noindent{\bf Step 1.} We first analyze the dynamics of the dual variable. 
\begin{lemma}\label{lm:dual:different}
	Suppose that $f(x)$ is in class $\cP^{M}_L$. Then for all $r\ge 0$, the iterates of D-GPDA satisfy 
	\begin{align}\label{eq:y:diff:bound}
	\hspace{-0.5cm}\|\lambda^{r+1}-\lambda^r\|_{\Sigma^{-2}}^2 
	 \le 2 \kappa \left(\frac{ \|\Upsilon^{-1} L (x^r-x^{r-1})\|^2}{M^2}+  \|w^{r+1}\|^2_{H}\right). 
	\end{align}
Further, for all $r\ge 0$, the iterates of xFILTER satisfy 
	\begin{align}\label{eq:y:diff:bound:+}
	&\|\lambda^{r+1}-\lambda^r\|_{\Sigma^{-2}}^2  \le \widetilde{\kappa} \bigg(\frac{3}{M^2} \|\Upsilon^{-1} L (x^r-x^{r-1})\|^2+3  \|w^{r+1}\|^2_{\Upsilon^2} + 3\|\Upsilon R (\epsilon^{r+1}-\epsilon^r)\|^2\bigg).
	\end{align}
 In the above we have defined the following { 
\begin{subequations}
	\begin{align}
	&\kappa: = \frac{1}{ \underline{\lambda}_{\min}(\Sigma F H^{-1}F^T\Sigma)},	\; \widetilde{\kappa}: = \frac{1}{ \underline{\lambda}_{\min}(\Sigma F \Upsilon^{-2}F^T\Sigma)}=\frac{1}{\underline{\lambda}_{\min}(\cL_G)} \label{eq:kappa}\\
	&w^{r+1}:= (x^{r+1}-x^r)-(x^r-x^{r-1})\label{eq:w}.
	\end{align} 
\end{subequations}}
\end{lemma}

\noindent {\bf Step 2.} In this step we analyze the dynamics of the AL.
\begin{lemma}\label{lm:AL:change}
	Suppose that $f(x)$ is in class $\cP^{M}_L$. Then for all $r\ge 0$, the iterates of  D-GPDA satisfy
	\begin{align}\label{eq:descent}
	&{\sf{AL}}^{r+1}- {\sf{AL}}^r \le -\frac{1}{2}\|x^{r+1}-x^r\|^2_{\Delta+2\Upsilon^2 - L/M} \nonumber\\
	&\quad + \kappa \left(\frac{2}{M^2} \|\Upsilon^{-1} L (x^r-x^{r-1})\|^2+  2  \|w^{r+1}\|^2_{H}\right).
	\end{align}
	Further,  for all $r\ge 0$, the iterates of xFILTER satisfy
	\begin{align}\label{eq:descent+}
	&\hspace{-0.4cm}{\sf{AL}}^{r+1} - {\sf{AL}}^r \le -\frac{1}{2}\|x^{r+1}-x^r\|^2_{\Upsilon^2 R- \frac{L}{M}}+ \langle \Upsilon^2 R \epsilon^{r+1},x^{r+1}-x^r\rangle \\
	& \hspace{-0.4cm}+ \widetilde{\kappa} \left(\frac{3}{M^2} \|\Upsilon^{-1} L (x^r-x^{r-1})\|^2+  3  \|w^{r+1}\|^2_{\Upsilon^2}+ 3\|\Upsilon R(\epsilon^{r+1}-\epsilon^r)\|^2\right).\nonumber
	\end{align}
\end{lemma} 

Before moving forward, we provide bounds for the important parameters $\kappa$ and $\tilde{\kappa}$. 
From \eqref{eq:kappa} we can express $\kappa$ as 
\begin{align}
\hspace{-0.4cm}\kappa 
=& \frac{1}{ \underline{\lambda}_{\min}(\Sigma F \Upsilon^{-1}\left(\Upsilon^{-1}B^T  \Sigma^2 B \Upsilon^{-1} + I \right)^{-1} \Upsilon^{-1}F^T\Sigma)} \nonumber\\
=&  \frac{1}{ \underline{\lambda}_{\min}(\left(\Upsilon^{-1}B^T  \Sigma^2 B \Upsilon^{-1} + I \right)^{-1} \cL_G)}\nonumber\\
\stackrel{\eqref{eq:normalized::gen:L}}=& \frac{1}{ \underline{\lambda}_{\min}\left((-\cL_G+ 2\Upsilon^{-1}\Delta\Upsilon^{-1} + I )^{-1}  \cL_G\right)}. \label{eq:kappa:full}
\end{align}
Similar derivation applies for $\tilde{\kappa}$. In summary we have 
\begin{align}\label{eq:kappa:bounds}
%\begin{align}
\kappa  \le \frac{\lambda_{\max}(2\Upsilon^{-1}\Delta\Upsilon^{-1} + I)}{\underline\lambda_{\min}(\cL_G)}, \quad  \widetilde{\kappa} = \frac{1}{\underline\lambda_{\min}(\cL_G)}. 
%\end{align}
\end{align}

\noindent {\bf Step 3.} In this step, we analyze the error sequences $\{\epsilon^{r+1}\}$ generated by the xFILTER. 
First we have the following well-known result on the behavior of the Chebyshev iteration; see, e.g., \cite[Chapter 6]{Tsynkov07} and \cite[Theorem 1, Chapter 7]{Samarskij89}. 
\begin{lemma}\label{lm:cheby}
	Consider using the Chebyshev iteration \eqref{eq:cheby:iteration} to solve $R x = d^r$. Define $x^{r+1}_*=R^{-1}d^r$,  with 
	\begin{align}
	R:= \Upsilon^{-2}(F^T\Sigma^2F + \Upsilon^2).
	\end{align}
	Define the following constants:
	\begin{align}\label{eq:def:xi}
	\xi(R): = \frac{\lambda_{\min}(R)}{\lambda_{\max}(R)}\le 1, \; \xi(\Upsilon^2): = \frac{\lambda_{\min}(\Upsilon^2)}{\lambda_{\max}(\Upsilon^2)}\le 1, \; \theta(R):= \lambda_{\min}(R)+\lambda_{\max}(R).
	\end{align} 
	Choose the following parameters:
	\begin{align*}
	\tau=\frac{2}{\theta(R)}, \; \alpha_1 = 2, \; \alpha_{t+1}=\frac{4}{4-\rho^2_0\alpha_t}, \; \rho_0=\frac{1-\xi(R)}{1+\xi(R)}.
	\end{align*}
	Then for any $\eta\in(0,1)$, in order to achieve the following accuracy
	\begin{align}\label{eq:u:contract}
	\|u_Q - x^{r+1}_*\|_{\Upsilon^2}^2\le \eta \|u_0-x^{r+1}_*\|_{\Upsilon^2}^2,
	\end{align}
	it requires the following number of iterations
	\begin{align*}
	Q\ge -\frac{1}{4}\ln (\eta/4)\sqrt{1/\xi(R)}. 
	\end{align*}
\end{lemma}
Recall that in Algorithm 2 the initial and final solutions for the Chebyshev iteration are assigned to $x^{r}$ and $x^{r+1}$, respectively. Define $\widetilde\epsilon^r:= u_0-x^{r+1}_*$, we have
\begin{align*}
Rx^r= R u_0 = R(u_0-x^{r+1}_*) + Rx^{r+1}_* := R{\widetilde\epsilon^r} +d^r, \forall~r\ge -1.
\end{align*}
Plugging in the definition of $d^r$ in \eqref{eq:R}, we obtain
\begin{align}\label{eq:epsilon:tilde}
 R{\widetilde\epsilon^r} = R x^r +  \Upsilon^{-2}(\nabla f(x^r) + F^T\lambda^r - \Upsilon^2 x^r).
\end{align}
Using  the definition of $\epsilon^{r+1}$ in \eqref{eq:linear:cheby:error:equiv}, and the fact that $R$ is invertible, we obtain the following key relationship
\begin{align}\label{eq:error:relation}
\epsilon^{r+1} - {\widetilde\epsilon^r}  = x^{r+1}-x^r, \; \forall~r\ge -1. 
\end{align}
Recall that $\epsilon^{r+1}:= x^{r+1}-x^{r+1}_*$, and $x^{r+1} = u_Q$, $x^r = u_0$, then \eqref{eq:u:contract} implies 
\begin{align}\label{eq:error:relation:2}
\|\epsilon^{r+1}\|_{\Upsilon^2}^2\le \eta \|\widetilde{\epsilon}^r\|_{\Upsilon^2}^2. 
\end{align}
By combining Lemma \ref{lm:cheby}, \eqref{eq:error:relation} and \eqref{eq:error:relation:2}, the following result provides some essential relationships between the error sequences $\{\epsilon^{r+1}\}$ incurred by running finite number of C-iterations, with the outer-loop iterations $\{x^{r+1}\}$.
\begin{lemma}\label{lemma:relation:epsilon}
	{\it Choose the inner iteration of xFILTER as 
			\begin{align}\label{eq:Q:specific}
			\hspace{-0.4cm}Q = -\frac{1}{4}\ln \left(\frac{\theta^2}{16+128 M \max\{\lambda_{\max}(\Upsilon^2 R),1\}}\right)\sqrt{1/\xi(R)}.
			\end{align}
		where $\theta: = \xi(\Upsilon^2R)\xi(\Upsilon^2)\times\min\{1, \lambda_{\min}(\Upsilon^2)\}$.
		Then we have the following inequalities   
			\begin{subequations}\label{eq:epsilon:bound}
				\begin{align}
				&\hspace{-0.1cm} \|\Upsilon^2 R\epsilon^{r+1}\|^2 \le \frac{1}{16M}\|x^{r+1}-x^r\|_{\Upsilon^2  R}^2,\label{eq:a}\\
				&\hspace{-0.1cm} \|\epsilon^{r+1}\|_{\Upsilon^2 R}^2 \le \frac{1}{16M}\|x^{r+1}-x^r\|_{\Upsilon^2 R}^2,\label{eq:b}\\
				&\hspace{-0.1cm} \|\Upsilon R\epsilon^{r+1}\|^2 \le \frac{1}{16M}\|x^{r+1}-x^r\|_{\Upsilon^2  R}^2,\label{eq:c}\\
				&\hspace{-0.1cm}\langle \Upsilon^2 R \epsilon^{r+1}\hspace{-0.2cm} , x^{r+1}\hspace{-0.2cm} -x^r\rangle\le \frac{3}{16}\|x^{r+1}-x^r\|_{\Upsilon^2 R}^2, \label{eq:d}\\
				& \hspace{-0.1cm}\langle \Upsilon^2 R \epsilon^{r}\hspace{-0.1cm}, x^{r+1} \hspace{-0.2cm}- x^r\rangle\le \frac{1}{8}\|x^{r}-x^{r-1}\|_{\Upsilon^2 R}^2 +\frac{1}{16}\|x^{r+1}-x^{r}\|_{\Upsilon^2 R}^2.\label{eq:e}
				\end{align}
	\end{subequations}}
\end{lemma}
 
\subsection{Proof of Lemma \ref{lemma:relation:epsilon}}
\noindent{\bf Proof.}
%For  simplicity we define $\xi: =\xi(\Upsilon^2 R)\times\xi(\Upsilon^2)$. 
Let us choose 
\begin{align}\label{eq:sigma}
\eta=\theta^2/(4+32M \max\{\lambda_{\max}(\Upsilon^2 R),1\} ).
\end{align} 
Then from Lemma \ref{lm:cheby}, it is clear that if $Q$ satisfies \eqref{eq:Q:specific}, then 
\begin{align}\label{eqLepsilon:relation}
\|\epsilon^{r+1}\|_{\Upsilon^2 }^2\le \eta \|\tilde{\epsilon}^r\|_{\Upsilon^2 }^2.
\end{align}
Note that $\Upsilon^2 R = F^{\top} \Sigma^2 F + \Upsilon^2\succ 0$, then it follows that {\small
	\begin{align*}
	&\|\Upsilon^2  R \epsilon^{r+1}\|^2 \le \frac{\lambda_{\max}(R \Upsilon^2 \Upsilon^2  R )}{\lambda_{\min}(\Upsilon^2)}\|\epsilon^{r+1}\|_{\Upsilon^2}^2 \nonumber\\
	&\stackrel{\eqref{eq:error:relation:2}}\le \frac{\eta\lambda_{\max}(R \Upsilon^2 \Upsilon^2  R )}{\lambda_{\min}(\Upsilon^2)} \|\tilde{\epsilon}^r\|_{\Upsilon^2 }^2\le \frac{\eta\lambda_{\max}(R \Upsilon^2 \Upsilon^2  R )\lambda_{\max}(\Upsilon^2)}{\lambda_{\min}(\Upsilon^2)} \|\tilde{\epsilon}^r\|^2\nonumber\\
	&\le \frac{\eta\lambda_{\max}(R \Upsilon^2 \Upsilon^2  R )\lambda_{\max}(\Upsilon^2)}{\lambda_{\min}(R \Upsilon^2 \Upsilon^2  R )\lambda_{\min}(\Upsilon^2)} \|\Upsilon^2 R\tilde{\epsilon}^r\|^2\le \eta \theta^{-2}\|\Upsilon^2 R\tilde{\epsilon}^r\|^2. 
	\end{align*}}
Using the above relation, we can then obtain the following{\small
	\begin{align*}
	\|\Upsilon^2  R\epsilon^{r+1}\|^2 &\le 2\eta\theta^{^{-2}}( \|\Upsilon^2 R{\epsilon}^{r+1}\|^2+ \|\Upsilon^2 R (\epsilon^{r+1}-\tilde{\epsilon}^r)\|^2)\nonumber\\
	&\stackrel{\eqref{eq:error:relation}}\le 2\eta\theta^{-2} (\|\Upsilon^2 R{\epsilon}^{r+1}\|^2+ \|\Upsilon^2 R(x^{r+1}-x^r)\|^2).
	\end{align*}}
Therefore, we obtain{\small
	\begin{align*}
	\|\Upsilon^2  R \epsilon^{r+1}\|^2 \le 2\eta\theta^{-2}/(1-2\eta\theta^{-2})\|\Upsilon^2 R(x^{r+1}-x^r)\|^2.
	\end{align*}}
Plugging the definition of $\eta$ in \eqref{eq:sigma}, we have 
\begin{align*}
\|\Upsilon^2  R \epsilon^{r+1}\|^2 &\le \lambda_{\max}(\Upsilon^2 R)2\eta\theta^{-2}/(1-2\eta\theta^{-2})\|x^{r+1}-x^r\|_{\Upsilon^2 R}^2\nonumber\\
&\stackrel{\eqref{eq:sigma}}\le 1/(16M)\|x^{r+1}-x^r\|_{\Upsilon^2 R}^2, \quad\forall~r\ge -1.
\end{align*}
To obtain the second inequality, notice that  
\begin{align}\label{eq:epsilon:R:2}
\|\epsilon^{r+1}\|_{\Upsilon^2 R}^2 \le \theta^{-1} \eta \|\widetilde{\epsilon}^r\|_{\Upsilon^2  R}^2 \le \theta^{-2} \eta \|\widetilde{\epsilon}^r\|_{\Upsilon^2  R}^2
\end{align}
where the last inequality is due to the fact that $\theta\le 1$. 
Then repeating the above derivation we can obtain the desired result. 
The third inequality in \eqref{eq:epsilon:bound} can be derived in a similar way, and the last two in \eqref{eq:epsilon:bound} can be obtained by using Cauchy-Swartz inequality.  
\QED

Clearly, using the Chebyshev iteration is one critical step that ensures fast reduction of the error $\{\epsilon^{r+1}\}$. In particular, to achieve constant reduction of error, the total number of required Chebyshev iteration is proportional to $\sqrt{1/\xi(R)}$, rather than $1/\xi(R)$ in conventional iterative scheme such as the Richardson's iteration \cite{Tsynkov07}. 	Such a choice enables the final bound to be dependent on $\sqrt{1/\xi(\cG)}$, rather than ${1/\xi(\cG)}$.  

\noindent {\bf Step 4.} Let us construct the following potential functions (parameterized by constants $c, \widetilde{c}>0$){
\begin{subequations}
\begin{align}
& \hspace{-0.2cm} P_{c}(x^{r+1},x^r,\lambda^{r+1}):={\sf{AL}}^{r+1}+ \frac{2\kappa}{M^2} \|  \Upsilon^{-1}  L (x^{r+1}-x^{r})\|^2 \nonumber\\
& \hspace{-0.2cm} + \frac{c}{2}\left(\|\Sigma Fx^{r+1}\|^2 +\|x^{r+1}-x^r\|_{H+L/M}^2\right). \\%\quad \forall~r\ge -1.
& \hspace{-0.2cm} \widetilde{P}_{\tilde{c}}(x^{r+1},x^r,\lambda^{r+1}):={\sf{AL}}^{r+1}+ \frac{3\widetilde{\kappa}}{M^2} \|  \Upsilon^{-1}  L (x^{r+1}-x^{r})\|^2 \\
& \hspace{-0.2cm} + \frac{3\widetilde{\kappa}} {8}\|x^{r+1}-x^r\|^2_{\Upsilon^2 R}+ \frac{\widetilde{c}}{2}\left(\|\Sigma Fx^{r+1}\|^2 +\|x^{r+1}-x^r\|_{\Upsilon^2+\frac{\Upsilon^2R}{4}+\frac{L}{M}}^2\right)\nonumber. %\quad \forall~r\ge -1.
\end{align}
\end{subequations}
}
For notational simplicity we will denote them as $P^{r+1}$ and $\widetilde{P}^{r+1}$, respectively. 
In the following we show that when the algorithm parameters are chosen properly, the potential functions will decrease along the iterations. 
\begin{lemma}\label{lm:Potential}
	Suppose that $f(x)$ is in class $\cP^M_{L}$, and that the parameters of D-GPDA are chosen as below 
	\begin{subequations}\label{eq:constants}
		\begin{align}
		\hspace{-0.6cm}	&c=\max\{6\kappa,1\}, \; 	\Upsilon^2 \succeq \frac{L \Upsilon^{-2} L}{M^2},\label{eq:constants:choices:kappa}\\
		& \frac{1}{2}\left(\Delta+\Upsilon^2\right) -  \frac{L}{M} - \frac{4\kappa}{M^2}L \Upsilon^{-2}L - \frac{2c L}{M} \succeq 0.\label{eq:constants:choices}
		\end{align}
	\end{subequations}
	
	Then for all $r\ge 0$, we have 
	\begin{align}\label{eq:final:descent}
	&P^r- P^{r+1}\ge  \frac{1}{4}\|x^{r+1}-x^r\|^2_{\Delta+ \Upsilon^2} + \kappa \|w^{r+1}\|_{H}^2.
	\end{align}
\end{lemma}
\begin{lemma}\label{lm:Potential+}
Suppose that $f(x)$ is in class $\cP^{M}_L$, $Q$ is chosen according to \eqref{eq:Q:specific}, and the rest of the parameters of xFILTER are chosen as below {

	\begin{subequations}\label{eq:constants+}
		\begin{align}
		&	\hspace{-0.6cm}	\tilde{c}=8\tilde{\kappa} = \frac{8}{ \underline{\lambda}_{\min}(\Sigma F \Upsilon^{-2}F^T\Sigma)}, \; 	\Upsilon^2 \succeq \frac{L \Upsilon^{-2} L}{M^2},\label{eq:constants:choices:kappa+}\\
		& 	\hspace{-0.6cm}	 ({1}/{4}- 3\tilde{\kappa}- \tilde{c} )\Upsilon^2 R- (1+2\tilde{c}) L/M - \frac{6\tilde{\kappa}}{M^2}L \Upsilon^{-2}L \succeq 0.\label{eq:constants:choices+}
		\end{align}
	\end{subequations}}
	Then for all $r\ge 0$, we have 
	\begin{align}\label{eq:final:descent+}
	&\widetilde{P}^r- \widetilde{P}^{r+1}\ge  \frac{1}{8}\|x^{r+1}-x^r\|^2_{\Upsilon^2 R} + \tilde{\kappa} \|w^{r+1}\|_{\Upsilon^2}^2.
	\end{align}
\end{lemma}

\noindent{{\bf Step 5.}
Next we show the lower and upper boundedness of the potential function. 

\begin{lemma}\label{lm:bounded}
	Suppose that $f(x)$ is in class $\cP^{M}_L$ and the parameters are chosen according to \eqref{eq:constants}. Then the iterates generated by D-GPDA satisfy
	\begin{subequations}
	\begin{align}
	&P^{r+1}\ge \underline{f}>{ - \infty}, \quad \forall~r>0, \label{eq:lower:bound:case1:lower}\\
	&P^0\le f(x^0) + \frac{2}{M} d_0^T L^{-1} d_0,  \label{eq:lower:bound:case1:P}
	\end{align}
	\end{subequations}
	where $d_0$ is defined in \eqref{eq:d0}.
	
%	where the constant in the last ex is given in \eqref{eq:c1}. 

	Similarly, for xFILTER the function $\widetilde{P}^{r+1}$ has the same expression as in \eqref{eq:lower:bound:case1:lower}, and
	\begin{align}\label{eq:P0:bound+}
&\widetilde{P}^0  \le f(x^0)+   \frac{5}{M} d_0^T L^{-1} d_0. 
	\end{align}
\end{lemma}
\noindent {\bf Step 6.} We are ready to derive the final bounds for the convergence rate of the proposed algorithms. 
\begin{theorem}\label{thm:final:bound}
	Suppose that $f(x)$ is in class $\cP^{M}_L$ and the parameters are chosen according to \eqref{eq:constants}. %Further assume that 
%	\begin{align}\label{eq:upsilon}
%	\end{align}
Let $T$ denote an iteration index in which D-GPDA satisfies
\begin{align}\label{eq:def:eT}
e(T):=\min_{r\in[T]}\bigg\|{1}/{M}\sum_{i=1}^{M}\nabla f_i(x^{r}_i)\bigg\|^2 + \|\Sigma F x^{r}\|^2\le\epsilon.
\end{align}
Then we have the following bound for the error:
\begin{align}\label{eq:overall:bound}
&\epsilon \le  C_1\times\frac{C_2}{T}, \; 
\mbox{with}\; C_1 :=  8\bigg(f(x^0)-\underline{f} +  \frac{2}{M} d_0^T L^{-1} d_0\bigg)\nonumber\\
&\quad\quad\quad  C_2 :=4 \sum_{(i,j):i\sim j} {{\sigma^2_{ij}}} + \sum_{i=1}^{M}\beta^2_i + 4.
\end{align}
Similarly, for xFILTER when the parameters are chosen according to \eqref{eq:constants+} and \eqref{eq:Q:specific}, the same equation
 \begin{align} 
&\epsilon \le  \widetilde{C}_1 \times\frac{\widetilde{C}_2}{T_r}
\end{align}
holds true (with $T_r$ denoting the total number of {\it outer} iterations), with the following constants
\begin{subequations}\label{eq:overall:bound+}
	\begin{align}
	&\widetilde{C}_1 :=  f(x^0)-\underline{f} +  \frac{5}{M}d_0^T L^{-1}d_0\\
	&\widetilde{C}_2 := 128\left(\sum_{i=1}^{M}\beta^2_i + 3 +\frac{1}{32\widetilde{\kappa}}\right).
	\end{align}
\end{subequations}

\end{theorem}

We note that one key difference between the two rates is that, the constant $C_2$ for D-GPDA depends explicitly on $\sigma_e$'s, while its counterpart for xFILTER depends on $1/\widetilde{\kappa}$ instead. Further, for xFILTER, the constant $\widetilde{\kappa}$ in \eqref{eq:kappa:bounds} only depends on $\underline{\lambda}_{\min}(\cL_G)$, while for D-GPDA $\kappa$ is further dependent on $\lambda_{\max}(\Upsilon^{-1}\Delta \Upsilon^{-1})$. These properties will be leveraged later when choosing algorithm parameters to ensure that optimal rates for different problems and networks are obtained.

\section{Rate Bounds and Tightness} \label{sec:tightness:uniform}
In this section we provide explicit choices of various parameters, and discuss the tightness of the resulting bounds for D-GPDA and xFILTER. 

\subsection{Parameter Selection and Rate Bounds for D-GPDA}
Let us pick the following parameters for D-GPDA {
	\begin{align}
	\hspace{-0.3cm}	\sigma^2_{ij} &=\frac{\beta^2 \sqrt{L_i L_j}}{\sqrt{d_i d_j}}, \; \Upsilon^2 = \beta^2 L, \; \beta^2 = \frac{80\max\{\lambda_{\max}(W),1\}}{\min\{\underline{\lambda}_{\min}(\cL_G),1\}M}. \label{eq:beta:gpda}
	\end{align}}
It follows that the following relations hold  %the above choice of parameters, it is easy to check that 
\begin{align}
\hspace{-0.2cm}\Delta=  {\beta^2} W,\quad\ \beta^2_i=\beta^2 {L_i}, \; \forall~i,\quad \kappa \stackrel{\eqref{eq:kappa:bounds}}\le \frac{1+2\lambda_{\max}(W)}{\min\{\underline\lambda_{\min}({\tL}),1\}}.
\end{align}
In the above definitions, we have defined $W\in\mathbb{R}^{M}$ as a diagonal matrix with 
$$[W]_{ii}= \frac{\sqrt{L_i}}{\sqrt{d_i}}\sum_{q:q\sim i}\frac{\sqrt{L_q}}{\sqrt{d_q}},$$ 
and that
\begin{align}
[{\cL_G}]_{ij} = \left\{
\begin{array}{ll}
\sum_{q:q\sim i}\frac{1}{\sqrt{d_q d_i}} & \mbox{if}~i=j\\
-\frac{1}{\sqrt{d_i d_j}} & \mbox{if}~(ij)\in\cE, i\ne j\\
0 & \mbox{otherwise}.
\end{array}
\right.\label{eq:hL}
\end{align}
Note that when $d_i = d_j, \; \forall~i,j$, we have $\cL_G = \cL$. We have the following result. 
\begin{theorem}\label{thm:rate:bounds}
	Consider using D-GPDA to solve problems in class $(\cP^{M}_{L}, \cN^M_D)$, using parameters in \eqref{eq:beta:gpda}. Then the condition \eqref{eq:constants:choices}  will be satisfied. Further, to achieve $e(T)\le \epsilon$, it requires at most the following number of iterations 
	\begin{align}\label{eq:dgpda:rate:bounds}
	T&\le \frac{8}{\epsilon}\bigg(f(x^0)-\underline{f} +  \frac{2}{M} \|d_0\|_{L^{-1}}^2\bigg)\times C_2
	\end{align}
	where $C_2$ is given by [with $W$ and $\tcL$ defined in \eqref{eq:hL}]
	\begin{align}\label{eq:C2}
	C_2 &\le\frac{320  \max\{\lambda_{\max}(W),1\}}{\min\{\underline{\lambda}_{\min}(\cL_G),1\}}\sum_{(i,j):i\sim j}\left(\frac{\sqrt{L_i L_j}}{\sqrt{d_id_j} M} + \frac{\bar{L}}{4}\right) + 4.
	\end{align}
\end{theorem}
\noindent{\bf Proof.} For D-GPDA, use the parameters in \eqref{eq:beta:gpda}, we have 
$$c\le \frac{6+12\lambda_{\max}(W)}{\min\{\underline\lambda_{\min}(\cL_G),1\}}, \quad \Upsilon^2 = \frac{80 \max\{\lambda_{\max}(W),1\}}{\min\{\underline{\lambda}_{\min}(\cL_G),1\} M} L.$$
Therefore to ensure condition \eqref{eq:constants:choices}, it suffices to ensure the following
\begin{align}
&\frac{40  \max\{\lambda_{\max}(W),1\}}{\min\{\underline{\lambda}_{\min}(\cL_G),1\} M} L- \frac{(4+8\lambda_{\max}(W))}{M 80 \max\{\lambda_{\max}(W),1\}}L -  \frac{1}{M} L - \frac{6+12\lambda_{\max}(W)}{\min\{\underline\lambda_{\min}({\cL_G}),1\}}\frac{2}{M} L\succ 0.
\end{align}
It is easy to check that this inequality will be satisfied using the above choice of parameters. Using these choices, we can obtain the desired expression for $C_2$.  \QED

\subsection{Parameter Selection and Rate Bounds for xFILTER}
First, recall that we have defined the matrix  $\tL$ and $\hL$ as follows  [see the definition in \eqref{eq:tilde:L}]
\begin{align*}
\ctL &= L^{-1/2}P^{-1/2} F^T K F P^{-1/2}L^{-1/2}, \\
\hL &= L^{-1/2}F^T K F L^{-1/2}.
\end{align*}	
Below we will provide two different choices of parameters. 

\noindent{\bf Choice I.} We will focus on a class of graphs such that there exists an {\it absolute} constant $k>0$ such that the following holds (i.e., the degrees of the nodes are not quite different from their averages):
\begin{align}\label{eq:def:k}
k P\succeq \bar{d} I_M. 
\end{align}
The above condition says that  the degrees of the nodes are not quite different from their averages. 
For example the following graphs satisfy \eqref{eq:def:k}: Complete graph ($k=1$), star graph ($k=2$), grid graph ($k=2$), cubic graph ($k=1$), path graph ($k=2$), and  any regular graph ($k=1$). 

For the class of graphs satisfy \eqref{eq:def:k}, let us pick the parameters for xFILTER as follows 
\begin{align}\label{eq:choice:uniform:2+}
\Sigma^2= \frac{48\times 96k}{  \sum_{i}d_i \underline{\lambda}_{\min}(\tcL)} K, \quad \Upsilon^2 = \frac{96 k}{\sum_{i}d_i} P^{1/2} L P^{1/2}.
\end{align}
Using the above choice, we have 
\begin{align}\label{eq:beta:specialize}
\beta^2_i = \frac{96 L_i d_i k}{\sum_i d_i}
\end{align}
and that the matrix $\Upsilon$ satisfies the following
\begin{align}\label{eq:upsilon:bound:1}
\Upsilon^2 = \frac{96 k}{\sum_{i}d_i} P^{1/2} L P^{1/2}\succeq \frac{96}{M} L. 
\end{align}
Plugging these choices to the generalized Laplacian $\cL_G$ in \eqref{eq:normalized::gen:L} we obtain
\begin{align}
\cL _G& = \Upsilon^{-1} F^T \Sigma^2 F \Upsilon^{-1} \nonumber\\
& = \frac{48}{\underline{\lambda}_{\min}(\ctL)} L^{-1/2}P^{-1/2} F^T  K F P^{-1/2}L^{-1/2} =  \frac{48}{\underline{\lambda}_{\min}(\ctL)}\ctL.
\end{align}
%Then according to \eqref{eq:choice:uniform:2+},  %and using the fact that $\xi(\cG)\ge 1/M^2$, 
%we have the following relations
%$$\Delta = \frac{48\beta^2}{\underline{\lambda}_{\min}(\cL)}  P, \quad \ctL= \frac{48}{\underline{\lambda}_{\min}(\cL)} {\cL}.$$ 
Therefore by \eqref{eq:kappa:bounds} we have 
\begin{align}\label{eq:kappa:star:line+}
\tilde{\kappa} \ =  \frac{\underline{\lambda}_{\min}(\ctL)}{ 48 \underline{\lambda}_{\min}(\ctL)}= \frac{1}{48}.
\end{align}
Also in this case we have 
\begin{align*}
R& = \Upsilon^{-2}F^T \Sigma^2 F  + I = \frac{48}{\underline{\lambda}_{\min}(\tcL)} P^{-1/2}L^{-1} P^{-1/2} F^T  K F + I.
\end{align*}
By noting that the matrix $P^{-1/2}L^{-1} P^{-1/2} F^T K F $ and $\ctL$ has the same set of eigenvalues, we obtain 
\begin{subequations}\label{eq:spectrum:R}
	\begin{align}
	\lambda_{\max}(R)&\le \left(\frac{48\lambda_{\max}(\ctL)}{\underline{\lambda}_{\min}(\ctL)}  + 1\right) \le \frac{50}{\xi(\ctL)} , \; {\lambda}_{\min}(R)= 1,\\
	\xi(R)& \ge 1/\left(\frac{48\lambda_{\max}(\ctL)}{\underline{\lambda}_{\min}(\ctL)}  + 1\right)\ge \frac{\xi(\ctL)}{50}. 
	\end{align}
\end{subequations}

\noindent{\bf Choice II.} For general graphs not necessarily satisfying \eqref{eq:def:k},  let us pick the parameters for xFILTER as follows 
\begin{align}\label{eq:choice:uniform:2+:2}
\Sigma^2= \frac{48\times 96}{M\underline{\lambda}_{\min}(\hL)} K, \quad \Upsilon^2 = \frac{96}{M} L.
\end{align}
Using the above choice, we have 
\begin{align}\label{eq:beta:specialize:2}
\beta^2_i = \frac{96 L_i}{M}. 
\end{align}
We have that
\begin{align}
\cL_G & = \Upsilon^{-1} F^T \Sigma^2 F \Upsilon^{-1} = \frac{48}{\underline{\lambda}_{\min}(\hL)} L^{-1/2}F^T  K F L^{-1/2} =  \frac{48}{\underline{\lambda}_{\min}(\hL)}\hL.
\end{align}
%Then according to \eqref{eq:choice:uniform:2+},  %and using the fact that $\xi(\cG)\ge 1/M^2$, 
%we have the following relations
%$$\Delta = \frac{48\beta^2}{\underline{\lambda}_{\min}(\cL)}  P, \quad \ctL= \frac{48}{\underline{\lambda}_{\min}(\cL)} {\cL}.$$ 
Therefore by \eqref{eq:kappa:bounds} we have 
\begin{align}\label{eq:kappa:star:line+:2}
\tilde{\kappa} \ =  \frac{\underline{\lambda}_{\min}(\hL)}{ 48 \underline{\lambda}_{\min}(\hL)}= \frac{1}{48}.
\end{align}
Also in this case we have 
\begin{align*}
R& = \Upsilon^{-2}F^T \Sigma^2 F  + I = \frac{48}{\underline{\lambda}_{\min}(\hL)} L^{-1} F^T  K F + I.
\end{align*}
By noting that the matrix $L^{-1} F^T K F $ and $\hL$ has the same set of eigenvalues, we obtain 
\begin{subequations}\label{eq:spectrum:R:2}
	\begin{align}
	\lambda_{\max}(R)&\le \left(\frac{48\lambda_{\max}(\hL)}{\underline{\lambda}_{\min}(\hL)}  + 1\right) \le \frac{50}{\xi(\hL)} , \; {\lambda}_{\min}(R)= 1,\\
	\xi(R)& \ge 1/\left(\frac{48\lambda_{\max}(\hL)}{\underline{\lambda}_{\min}(\hL)}  + 1\right)\ge \frac{\xi(\hL)}{50}. 
	\end{align}
\end{subequations}

\noindent\begin{remark}{\bf{(Choices of Parameters)}}
The main difference between the above two choices of parameters is whether $\Upsilon^2$ is scaled with the degree matrix or not. The resulting bounds are also dependent on the spectral gap for $\widetilde{L}$ and $\hL$, one inversely scaled with the degree matrix, and the other does not. Note that the spectral gap of $\tL$ and $\hL$ may not be the same. For example for a star graph with $L_i=L_j$, $\xi(\hL) =\mathcal{O}(1/M)$ but $\xi(\ctL)=\mathcal{O}(1)$. Therefore one has to be careful in choosing these parameters so that $\xi(R)$ is made as large as possible. 

Additionally, since we are mainly interested in choosing the optimal parameters so that the resulting rate bounds will be optimal in their dependency on problem parameters, the absolute constants in the above parameter choices have not been optimized. 
\end{remark}

The following result is a direct consequence of the second part of Theorem \ref{thm:final:bound}. 
\begin{theorem}\label{thm:rate:bounds+}
	Consider using xFILTER to solve problems in class $(\cP^{M}_{L}, \cN^M_D)$, then the following holds.

\noindent 	{\bf Case I.} Further restricting $\cN^{M}_D$ to a subclass satisfying \eqref{eq:def:k}. If parameters in \eqref{eq:choice:uniform:2+} is used,  then the condition \eqref{eq:constants:choices+} will be satisfied. Further, to achieve $e(T)\le \epsilon$, it requires at most the following number of iterations  (where $T$ denotes the {\it total} iterations of the xFILTER algorithm)

\begin{align}\label{eq:dgpda:rate:bounds+}
	T&\le \frac{1}{\epsilon}\bigg(f(x^0)-\underline{f} +  \frac{5}{M} \|d_0\|_{L^{-1}}^2\bigg)\times \widetilde{C}_2
	\nonumber\\
	&\quad \times \frac{1}{4}\ln \left(\frac{16+128 M \max\{\lambda_{\max}(\Upsilon^2 R),1\}}{\theta^2}\right)\sqrt{1/\xi(R)} \nonumber\\
	&\le \frac{1}{\epsilon}\bigg(f(x^0)-\underline{f} +  \frac{5}{M} \|d_0\|_{L^{-1}}^2\bigg)\times \widetilde{C}_2
	\nonumber\\
	&\quad \times \frac{1}{4}\ln \left(\frac{ 50^2(M L_{\max}/L_{\min})^4\times (16+128 M \max\{ 50\times 96kL_{\max} , 1 \}) }{\xi^3(\ctL)\times\min\{1, 96^2k^2L^2_{\min}/M^2\}}\right)\sqrt{50/\xi(\ctL)}
	\end{align}
	where $\widetilde{C}_2$ is given by 
	\begin{align}\label{eq:c2+}
	\widetilde{C}_2 &\le  128\left({\frac{96 k}{\sum_{i=1}^{M}d_i}\sum_{i=1}^{M}d_i L_i}+ 19\right).
	\end{align}

\noindent 	{\bf Case II.} Suppose parameters in \eqref{eq:choice:uniform:2+:2} are used. Then the condition \eqref{eq:constants:choices+} will be satisfied. Further, to achieve $e(T)\le \epsilon$, it requires at most the following number of iterations 
 \begin{align}\label{eq:dgpda:rate:bounds+:2}
	T&\le \frac{1}{\epsilon}\bigg(f(x^0)-\underline{f} +  \frac{5}{M} \|d_0\|_{L^{-1}}^2\bigg)\times \widetilde{C}_2
	\nonumber\\
	&\quad \times \frac{1}{4}\ln \left(\frac{ 50^2( L_{\max}/L_{\min})^4\times (16+128 M \max\{ 50\times 96L_{\max}/M , 1 \}) }{\xi^3(\ctL)\times\min\{1, 96^2L^2_{\min}/M^2\}}\right)\sqrt{50/\xi(\ctL)}
	\end{align}
where $\widetilde{C}_2$ is given by 
\begin{align}\label{eq:c2+:2}
\widetilde{C}_2 &\le  128\left({\frac{96}{M}\sum_{i=1}^{M}L_i}+ 19\right).
\end{align}

\end{theorem}
We note that compared with the results in Theorem \ref{thm:final:bound}, the additional multiplicative term in \eqref{eq:dgpda:rate:bounds+} accounts for the Chebyshev iterations that are needed for every iteration $t$. It is interesting to observe that comparing with the previous result, the constant $\widetilde{C}_2$ in \eqref{eq:c2+} is independent on any graph parameters.  Such a desirable property turns out to be crucial for obtaining tight rate bounds.

%\vspace{-0.3cm}
\subsection{Tightness of the Upper Rate Bounds}

In this section, we present some tightness results of the upper rate bounds for our proposed D-GPDA and xFILTER. In particular,  we compare the expressions derived in Theorem \ref{thm:rate:bounds} -- \ref{thm:rate:bounds+}, and the lower bounds derived in Section \ref{sec:lower}, over different kinds of graphs and for different problems. We will mainly focus on the case with uniform Lipschitz constants, i.e.,  $L_i=U, \; \forall~i$. WE will briefly discuss the case of non-uniform Lipschitz constants at the end of this section. 

%\subsubsection{Problems with uniform Lipschitz constants}

First,  we consider the problem class $\cP^{M}_U$ with the following properties: 
%To provide detailed analysis, we divide the problem class $\cP^{L}_M$ into two subclasses: The uniform case
\begin{align}\label{eq:uniform}
L_1 = L_2 = \cdots L_M = \frac{1}{M}\sum_{i=1}^{M} L_i:= U, \quad L = U I_M.
\end{align}
{It follows that in this case $\tcL =\cL$, and $\hL = P^{1/2}\cL P^{1/2}$.} Let us first make some useful observations. 
\begin{remark}\label{rmk:gpda}
	Let us specialize the parameter choices for D-GPDA algorithm in \eqref{eq:beta:gpda} and derive the bounds for $C_2$ in \eqref{eq:C2} for two special graphs. 
	
	\noindent{\bf Complete graph.} For complete graphs we have $d_i=d_j= M-1,\forall~i,j$, which implies that $\cL_G = \cL$, so $\underline{\lambda}_{\min}(\cL_G) =M/(M-1)$. Because $L_i=L_i =U, \; \forall~i,j$, we have $W=I_M$. Therefore using the expression \eqref{eq:C2} we obtain the following:
	\begin{align}\label{eq:C:complete}
	C^{\rm comp}_2 &\le 400U+ 4.  
	\end{align}
	
	\noindent{\bf Cycle graph.} For cycle graph we have $d_i=d_j= 2,\forall~i,j$, which implies that $\cL_G = \cL$, and $\underline{\lambda}_{\min}(\cL_G) \ge 1/M^2$. Because  $L_i=L_i =U, \; \forall~i,j$, we have $W=I_M$. Therefore using the expression \eqref{eq:C2} we obtain the following:
	\begin{align}\label{eq:C:circle}
	C^{\rm cycle}_2 &\le{240 U M^2 + 4.} 
	\end{align}		
	It is clear that for cycle graph whose diameter is in $\mathcal{O}(M)$, the rate bounds is very large. 
\end{remark}
\begin{remark}\label{rmk:gpda+}
	Let us specialize the parameter choices for xFILTER algorithm in \eqref{eq:choice:uniform:2+} and derive the bounds for $\widetilde{C}_2\times 1/\sqrt{\xi(\ctL)}$ in \eqref{eq:c2+} for the following special graphs. {Note that because uniform $L_i$'s are assumed, we have $\ctL = \cL$.}
	
	\noindent{\bf Complete graph.} Complete graphs satisfy \eqref{eq:def:k} with $k=1$. It also satisfies $\underline{\lambda}_{\min}(\ctL) =M/(M-1)\ge 1$. Therefore using the expression \eqref{eq:c2+} we obtain the following:
	\begin{align}\label{eq:C:complete+}
	\widetilde{C}^{\rm comp}_2 \times \frac{1}{\sqrt{\xi(\ctL)}} &\le  12500 U + 2560. 
	\end{align}
	\noindent{\bf Grid graph.}   Grid graphs satisfy \eqref{eq:def:k} with $k=2$. It also satisfies $\underline{\lambda}_{\min}(\ctL) \ge 1/M$.  Therefore using the expression  \eqref{eq:c2+} we obtain the following:
		\begin{align}\label{eq:C:grid}
		\widetilde{C}^{\rm grid}_2 \times \frac{1}{\sqrt{\xi(\ctL)}}  &\le(12500 U + 2560)\times \sqrt{M}.
		\end{align}
		
	\noindent{\bf Star graph.} {Star graphs satisfy \eqref{eq:def:k} with $k=2$. It also has $\xi(\ctL) = 1/2$}. Therefore using the expression \eqref{eq:c2+} we obtain the following:
		\begin{align}\label{eq:C:star}
		\widetilde{C}^{\rm star}_2  \times \frac{1}{\sqrt{\xi(\ctL)}}  & \le(12500 U + 2560)\times \sqrt{2}.
		\end{align}
		
{\noindent{\bf Geometric graph.}  {For geometric graphs which place  the nodes uniformly in $[0,1]^2$ and connect any two nodes separated by a distance less than a radius $R\in(0,1)$. Then if the connectivity radius $R$ satisfies \cite{Duchi12}
\begin{align}
R = \Omega\left(\sqrt{\log^{1+\epsilon}(M)/M}\right), \quad \mbox{for any}~\epsilon>0,
\end{align} 
then with high probability
\begin{align}
\xi(\ctL)=\mathcal{O}\left({\frac{\log(M)}{M}}\right).
\end{align}
Further, from the proof of \cite[Lemma 10]{boyd2006randomized}, for any $\epsilon$ and $c>0$, if 
\begin{align}
R = \Omega\left(\sqrt{\log^{1+\epsilon}(M)/(M\pi)}\right)
\end{align} 
then with probability at least $1-2/M^{c-1}$, the following holds
\begin{align}
\log^{1+\epsilon} M -\sqrt{2} c \log M\le d_i \le \log^{1+\epsilon} M + \sqrt{2} c \log M, \; \forall~i. 
\end{align}
This means that \eqref{eq:def:k} is satisfied (with $k=1$) with high probability (also see discussion at the end of \cite[Section V]{Duchi12}). 
Therefore using the expression \eqref{eq:C2} we obtain the following:
\begin{align}\label{eq:C:geometric}
	\widetilde{C}^{\rm geometric}_2  \times \frac{1}{\sqrt{\xi(\ctL)}}  & \le(12500 U + 2560)\times \mathcal{O}\left({\frac{\sqrt{M}}{\sqrt{\log(M)}}}\right) .
\end{align}
}}

	\noindent{\bf Cycle/Path graph.} Cycle/path graphs  satisfy \eqref{eq:def:k} with $k=2$. We also have $\underline{\lambda}_{\min}(\ctL) \ge 1/M^2$ (see the discussion in Sec. \ref{sub:facts}). Therefore using the expression \eqref{eq:c2+} we obtain the following: 
	\begin{align}\label{eq:C:circle+}
	\widetilde{C}^{\rm cycle}_2 \times \frac{1}{\sqrt{\xi(\ctL)}}   &\le(12500 U + 2560)\times M.
	\end{align}		
\end{remark}
From the above comparison, it is clear that the rate bounds for xFILTER is about $\cO({M})$ times better than the D-GPDA for the path/cycle graph. 

We also note that for the xFILTER algorithm, the fact that $L_i=U, \; \forall~i$ implies that the matrix $\Sigma^2$ given in \eqref{eq:choice:uniform:2+}  is a multiple of identity matrix. Therefore by Remark \ref{rmk:distributed}, we can conclude that in this case xFILTER belongs to both $\cA$ and $\cA'$. 

Now we are ready to present our tightness analysis on D-GPDA and xFILTER.  
\begin{theorem}\label{thm:tight:uniform:1}
We have the following tightness results. 

%\vspace{-0.2cm}
\noindent{\bf (1)} Let $D=1$ and consider the class $(P^{M}_{U},\cN^M_D)$. Then D-GPDA is an optimal algorithm, and its convergence rate in \eqref{eq:dgpda:rate:bounds} is tight (up to a universal constant).
%\vspace{-0.2cm}

\noindent{\bf (2)} Let $D=M-1$ and consider the class $(P^{M}_{U},\cN^M_D)$. Then xFILTER is an optimal algorithm, and its convergence rate in \eqref{eq:dgpda:rate:bounds+} is tight (up to a polylog factor).

%\vspace{-0.2cm}
\noindent{\bf (3)} {More generally, consider the problem class $P^{M}_{U}$, and a subclass of $\cN^M_D$ satisfying \eqref{eq:def:k}. Then the convergence rate in \eqref{eq:dgpda:rate:bounds+} is tight (up to a polylog factor). }

\end{theorem}

\noindent {\bf Proof.} We divide the proof into different cases. %The proof combines a few previous lemmas. 

\noindent{\bf Case 1).} The network class is  a complete graph with $M$ nodes. Using the parameters in \eqref{eq:beta:gpda}, $C_2$ is given by \eqref{eq:C:complete}, and we have that 
$\Sigma^2 = \frac{80 U}{(M-1)M} I_E$.
Note that the following holds
$$\|Fx\|^2 = \sum_{(i,j):i\sim j}\|x_i-x_j\|^2.$$
If  \eqref{eq:def:eT} holds, then Theorem \ref{thm:final:bound} and Theorem \ref{thm:rate:bounds} imply
\begin{align}
&T\le 8\big(f(0)-\underline{f} + \frac{2}{MU}\|d_0\|^2 \big)\times\frac{400U+4}{\epsilon}.  \nonumber
\end{align}
For complete graph it is easy to check that $\xi(\cG) \ge 1$. 
Using the definition in \eqref{eq:criteria:local}, we also have
\begin{align*}
h^*_T& =\min_{r\in[T]}\big\|\frac{1}{M}\sum_{i=1}^{M}\nabla f_i(x_i^r)\big\|^2 + \frac{U}{M^2}\|Ax^r\|^2\nonumber\\
& \le \min_{r\in[T]}\big\|\frac{1}{M}\sum_{i=1}^{M}\nabla f_i(x_i^r)\big\|^2 + \frac{1}{80} \|\Sigma Ax^r\|^2\le e(T)\le \epsilon. 
\end{align*}

By comparing the lower bound derived in Lemma \ref{thm:2}, we conclude that the above rate bound is tight (up to some universal constants).

\noindent{\bf Case 2).} The network class is a path graph with $M=D+1$. From Section \ref{sub:facts} we have 
\begin{align}\label{eq:lambda:max:bound}
\xi(\cG)\ge \frac{1}{M^2}.
\end{align}
Further we note that condition \eqref{eq:def:k} satisfies with $k=2$. We have
\begin{align}
\ctL = P^{-1/2} F^T F P^{-1/2} = \cL. 
\end{align}
Therefore we conclude that 
\begin{align}\label{eq:hL:path}
\xi(\ctL) \ge \xi(\cG) \ge \frac{1}{M^2}.
\end{align}

Applying the above estimate to \eqref{eq:choice:uniform:2+}, we can choose 
\begin{align}\label{eq:sigma:express}
\Sigma^2=\frac{4608 U}{4(M-2) \underline{\lambda}_{\min}(\ctL)}  I_E,\quad  \Upsilon^2=\frac{96 U}{4(M-2)} P.
\end{align}

Using these choices, again we will have
\begin{align}\label{eq:R:estimate:path}
\xi(R) \stackrel{\eqref{eq:spectrum:R}}\ge \frac{\xi(\ctL)}{50}\stackrel{\eqref{eq:hL:path}}\ge \frac{1}{50 M^2}.
\end{align}	

Using these constants, and note $D\le M$, we have
\begin{align*}
 h^*_{T_r} %&\min_{r=1,\cdots,T} \big\|\frac{1}{M}\sum_{i=1}^{M}\nabla f_i(x_i^r)\big\|^2 + 50 M \hspace{-0.2cm}\sum_{(i,j): i\sim j}\hspace{-0.2cm}\sqrt{\frac{L_i L_j}{d_i d_j}}\|x^r_i-x^r_j\|^2\nonumber\\
&=\min_{r\in[T_r]}\big\|\frac{1}{M}\sum_{i=1}^{M}\nabla f_i(x_i^r)\big\|^2 +  \frac{U}{M \underline{\lambda}_{\min}(P^{1/2}\cL P^{1/2})} \sum_{(i,j):i\sim j}\|x_i-x_j\|^2\nonumber\\
&\le \min_{r\in[T_r]}\big\|\frac{1}{M}\sum_{i=1}^{M}\nabla f_i(x_i^r)\big\|^2 +   \frac{U}{\underline{\lambda}_{\min}(\cL)M}\|Fx^r\|^2\nonumber\\
&\stackrel{\eqref{eq:sigma:express}}\le \min_{r\in[T_r]}\big\|\frac{1}{M}\sum_{i=1}^{M}\nabla f_i(x_i^r)\big\|^2 +  \frac{1}{2304}\|\Sigma F x^r\|^2\le e(T_r),
\end{align*}
where in the first inequality we have used $P\succeq I_M$. 
Similarly as in the previous case, suppose $e(T_r)\le \epsilon$, then according to Theorem \ref{thm:rate:bounds+}
we have
\begin{align}
\epsilon\le \bigg(f(0)-\inf_{x}f(x) +\|d_0\|^2\frac{5}{M U}I_M\bigg)\times\frac{128(96U+19)}{T_r}.  \nonumber
\end{align}
Recall that for xFILTER, $T_r$ represents the number of times the dual update \eqref{eq:mu:update:+} is performed. Between two dual updates $Q$ primal iterations are performed, where the precise number is given in \eqref{eq:Q:specific}. According to \eqref{eq:R:estimate:path} we have 
\begin{align}
\sqrt{1/\xi(R)}  \le   13M.
\end{align}
Overall, the total number of iterations required is given by
	 \begin{align} 
	T&\le \frac{1}{\epsilon}\bigg(f(x^0)-\underline{f} +  \frac{5}{M U}\|d_0\|^2\bigg)\times 128(96U+19)
	\nonumber\\
	&\quad \times \frac{1}{4}\ln \left(\frac{ 50^2M^{10}\times (16+128 M \max\{ 50\times 192U , 1 \}) }{\min\{1, 96^2 \times 4 U^2/M^2\}}\right)\times 13M.
	\end{align}
This implies that the lower bound obtained in Theorem \ref{thm:1} is tight up to some universal constant and a ploylog factor in $M$, and the bound-achieving algorithm in class $\cA$ is the xFILTER.

\noindent{\bf Case 3).} The proof follows similar steps are in the previous case. When $L_i=L_j, \; \forall~i\ne j$, and when \eqref{eq:def:k} is satisfied, it is easy to verify that the following holds 
\begin{align}
{h^*_{T_r}\le e(T_r)}, \; \mbox{and} \; \tcL = \cL. 
\end{align}
To bound the total number of iteration required to achieve $h^*_{T_r}\le \epsilon$, note that when \eqref{eq:def:k} is satisfied, we can apply the bound \eqref{eq:dgpda:rate:bounds+} in Theorem \ref{thm:rate:bounds+} and obtain
	\begin{align}
	T&\le  \frac{1}{\epsilon}\bigg(f(x^0)-\underline{f} +  \frac{5}{M U} \|d_0\|^2\bigg)\times 128\left({96 k}U+ 19\right)
	\nonumber\\
	&\quad \times \frac{1}{4}\ln \left(\frac{50^2 M^4\times (16+128 M \max \{50\times 96kU, 1\})}{\xi^3(\cG)\times \min\{1, 96^2k^2U^2/M^2\}}\right)\sqrt{50/\xi(\cG)}.
	\end{align}
Comparing with the lower bound obtained in Theorem \ref{thm:1}, it is clear that apart from the multiplicative $\ln(\cdot)$ term, the remaining bound is in the same order as the lower bound given in \eqref{eq:lower:bound}. 
\QED

\begin{remark} {\bf (Optimal Number of Gradient Evaluations)}\label{rmk:gradient:evaluation}
It is important to note that the ``outer" iteration of the  xFILTER required to achieve $\epsilon$-local solution scales with $\mathcal{O}(U/\epsilon)$, which is independent of the network size. Because local gradient evaluation is only performed in the outer iterations, the above fact suggests that the total number of gradient evaluation required is also in this order, which is optimal because it is the same as what is needed for the centralized gradient descent. 
\end{remark}

\begin{remark}\label{rmk:tightness} {\bf (Performance Gap Between D-GPDA and xFILTER)}
	If we apply  D-GPDA to the path or cycle  graph, then according to Remark \ref{rmk:gpda}, the corresponding $C_2$, as well as the final upper bound, will be in $\mathcal{O}(M^2 U)$, which is $\cO(M)$ worse than the lower bound. %This is in sharp contrast to xFILTER which achieves the lower bound except for a polylog factor. 
	Intuitively, this phenomenon happens because of the following: in order to decompose the entire problem into the individual nodes, the $x$-update \eqref{eq:x:update} has to create a proximal term that matches the quadratic penalty $\|\Sigma A x\|^2$. But such an additional term forces the variables to stay close to their previous iteration. In contrast, xFILTER circumvents the above difficulty by leaving the quadratic penalty intact, but instead using a few fast and decomposable iterations to approximately solve the resulting problem. 
%	which is at least $M^2$ time worse than that of the D-GDPA+. 
\end{remark}
 
\begin{remark} {\bf (An Alternative Bound)}
	For problems and graphs in  $(\cP^{M}_U, \cN^{M}_D)$ without additional conditions, it can be verified that the second choice of the parameters \eqref{eq:choice:uniform:2+:2} gives the following convergence rates [cf. \eqref{eq:dgpda:rate:bounds+:2}]
		\begin{align}
	{T}= \mathcal{\widetilde{O}} \left(\big(f(0)-\inf_{x}f(x) +\|d_0\|^2\frac{5}{M U}I_M\big)\times\frac{U}{\epsilon} \times\frac{1}{\sqrt{\xi(P^{1/2} \cL P^{1/2})}}\right),
	\end{align}
	where the notation $\widetilde{\cO}$ denotes $\cO$ with a multiplicative ploylog factor. The above rate is proportional to the square root of the eigengap for the matrix $P^{1/2} \cL P^{1/2}$, which is the {\it unnormalized} Laplacian matrix for graph $\cG$. 
	\end{remark}

{\begin{remark} {\bf (Non-uniform Lipschitz Constants)}
We comment that for the general case $L_i\ne L_j$, $\forall~i,j$, we can use similar steps to  verify that the bound \eqref{eq:dgpda:rate:bounds+:2} derived in Theorem \ref{thm:rate:bounds+} is optimal, in the sense that they achieve the lower bound \eqref{eq:lower:bound:nonuniform}  predicted in Corollary \ref{thm:2}. 
\end{remark}}

\section{Numerical Results}
 This section presents numerical examples to show the effectiveness of the proposed algorithms. Two kinds of problems are considered, distributed binary classification and distributed neural networks training. We use the former one to demonstrate the behavior and scalability of our algorithm and use the latter one to show the practical performance.

  \subsection{Simulation Setup}
 In our simulations, all algorithms are implemented in MATLAB R2017a for binary classification problem and implemented in Python 3.6 for training neural networks, running on a computer node with two 12-core Intel Haswell processors and 128 GB of memory (unless otherwise specified).   Both synthetic and real data are used for performance comparison. For synthetic data, the feature vector   is randomly generated with standard normal distribution with zero mean and unit variance. The label vector  is randomly generated with uniformly distributed pseudorandom integers taking the values $\{-1, 1\}$.
 For real data, we use the breast cancer dataset  {\footnote{\href{https://archive.ics.uci.edu/ml/datasets/Breast+Cancer+Wisconsin+(Diagnostic)}{https://archive.ics.uci.edu/ml/datasets/Breast+Cancer+Wisconsin+(Diagnostic)}}} for binary classification and MNIST{\footnote{\href{http://yann.lecun.com/exdb/mnist/}{http://yann.lecun.com/exdb/mnist/}}} for training neural network. The breast cancer dataset contains a total of 569 samples each with 30 real positive features. The MNIST dataset contains a total of 60,000 handwritten digits, each with a $28\times 28$ gray scale image and a label from ten categories.
 
 \subsection{Distributed Binary Classification}
We consider a non-convex distributed binary classification problem \cite{antoniadis2011penalized}. The global consensus problem \eqref{eq:global:consensus:equiv:2} can be expressed as follows:
\begin{align*}
\min_{x\in\mathbb{R}^{SM}} \; f(x):=\frac{1}{M}\sum_{i=1}^{M} f_i(x_i),\quad \st\; x_i=x_j, \forall~(i,j)\in \cE.
\end{align*}
And each component function $f_i$  is expressed by
\begin{align*} 
f_i(x_i) = \frac{1}{B}   \sum_{j=1}^B \log{\left( 1+\exp(-y_{ij}x_i^Tv_{ij})\right) } +  
\sum_{s=1}^S \frac{\lambda \alpha x^2_{i,s}}{1+\alpha x_{i,s}^2}.
\end{align*}  
Here $v_{ij}\in \mathbb{R}^S$ denotes the feature vector with dimension $S$, $y_{ij}\in \{1, -1\}$ denotes the label for the $j$th date point in $i$th agent, and there are total $B$ data points for each agent. Unless otherwise noted, the graph $\cE$ used in our simulation is generated using the random geometric graph {and the graph parameter {$Ra$} is set to $0.5$.  The regularization parameter is set to $\lambda = 0.001, \alpha = 1$.

   \begin{figure}
% 	\vspace{-0.8cm}
 	\begin{minipage}[c]{0.33\linewidth}
 		\includegraphics[width=\linewidth]{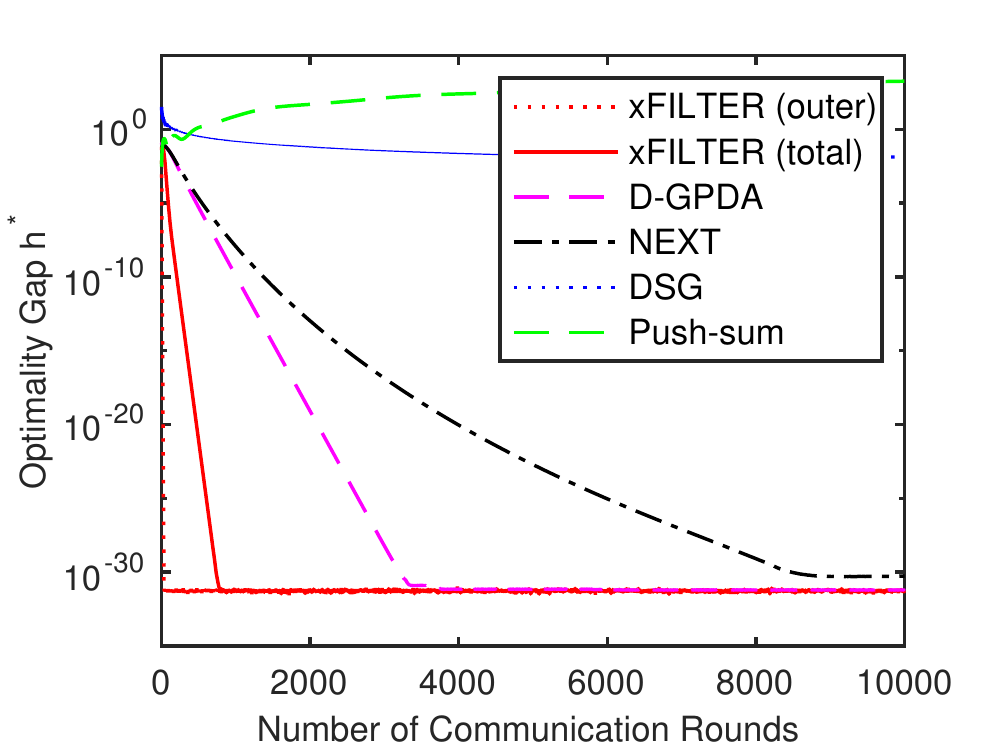}
 		\caption{$M=5, B=200, K=10$} \label{fig:ObjectiveComparison-1}
 	\end{minipage}%
 	\hfill
 	\begin{minipage}[c]{0.33\linewidth}
 		\includegraphics[width=\linewidth]{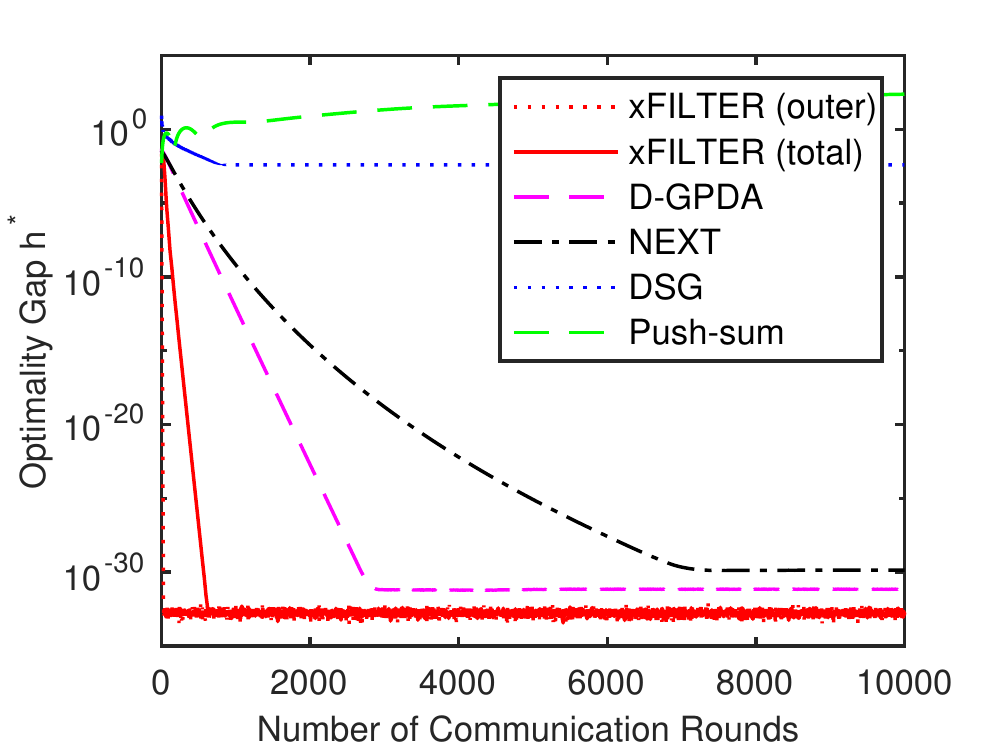}
 		\caption{$M=10, B=200, K=10$} \label{fig:ObjectiveComparison-2}
 	\end{minipage}
 	\hfill
 	\begin{minipage}[c]{0.33\linewidth}
 		\includegraphics[width=\linewidth]{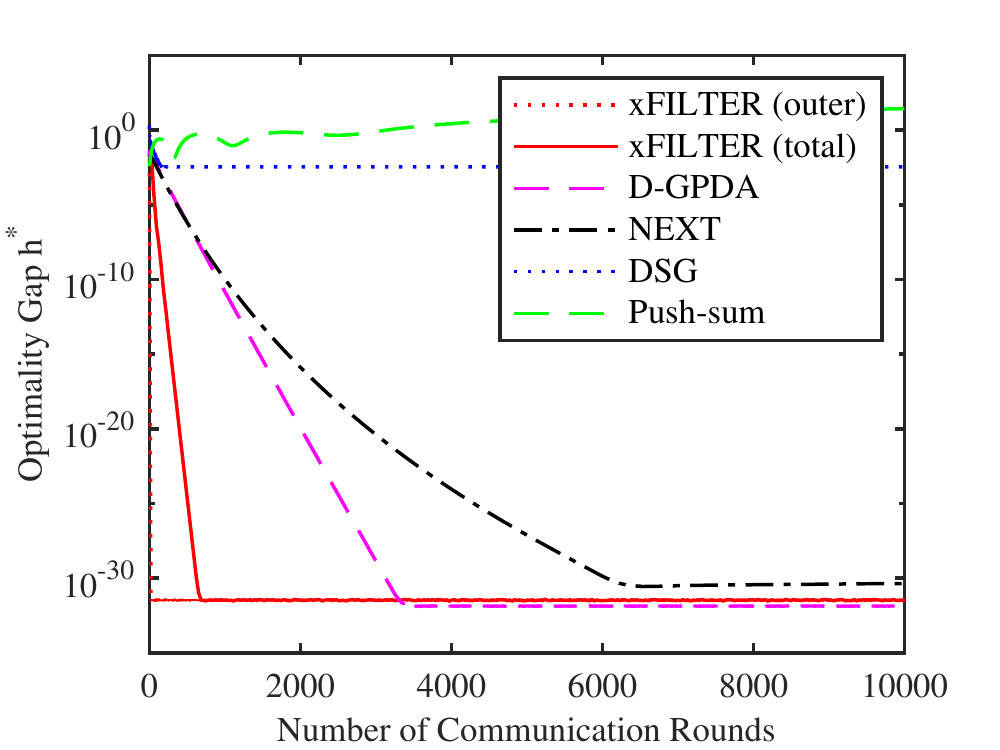}
 		\caption{$M=20, B=200, K=10$} \label{fig:ObjectiveComparison-3}
 	\end{minipage}
 \end{figure}

 \begin{figure}
% 	\vspace{-0.5cm}
 	\begin{minipage}[c]{0.33\linewidth}
 		\includegraphics[width=\linewidth]{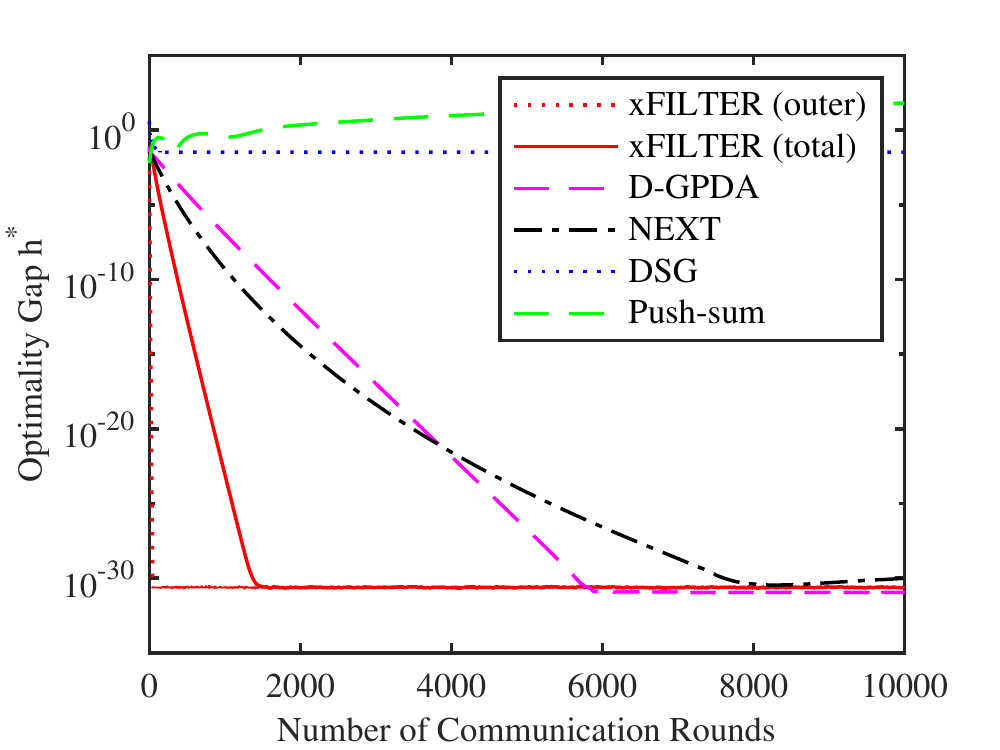}
 		\caption{$M=20, B=50, K=10$} \label{fig:ObjectiveComparison-4}
 	\end{minipage}% 
 	\hfill
 	\begin{minipage}[c]{0.33\linewidth}
 		\includegraphics[width=\linewidth]{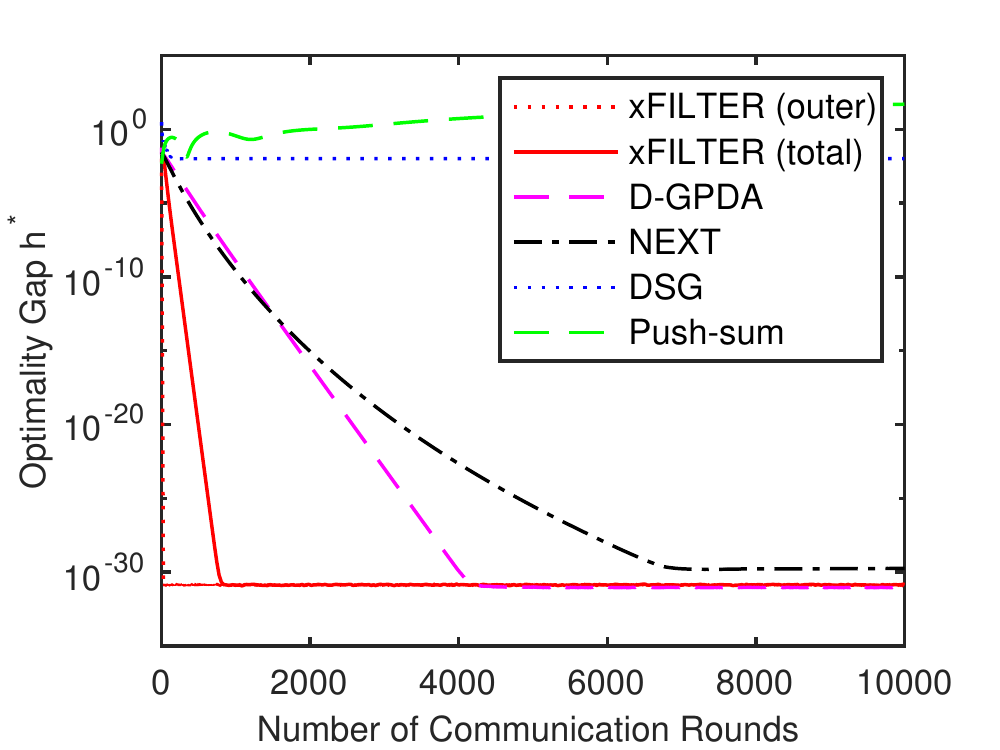}
 		\caption{$M=20, B=100, K=10$} \label{fig:ObjectiveComparison-5}
 	\end{minipage}
 	\hfill
 	\begin{minipage}[c]{0.33\linewidth}
 		\includegraphics[width=\linewidth]{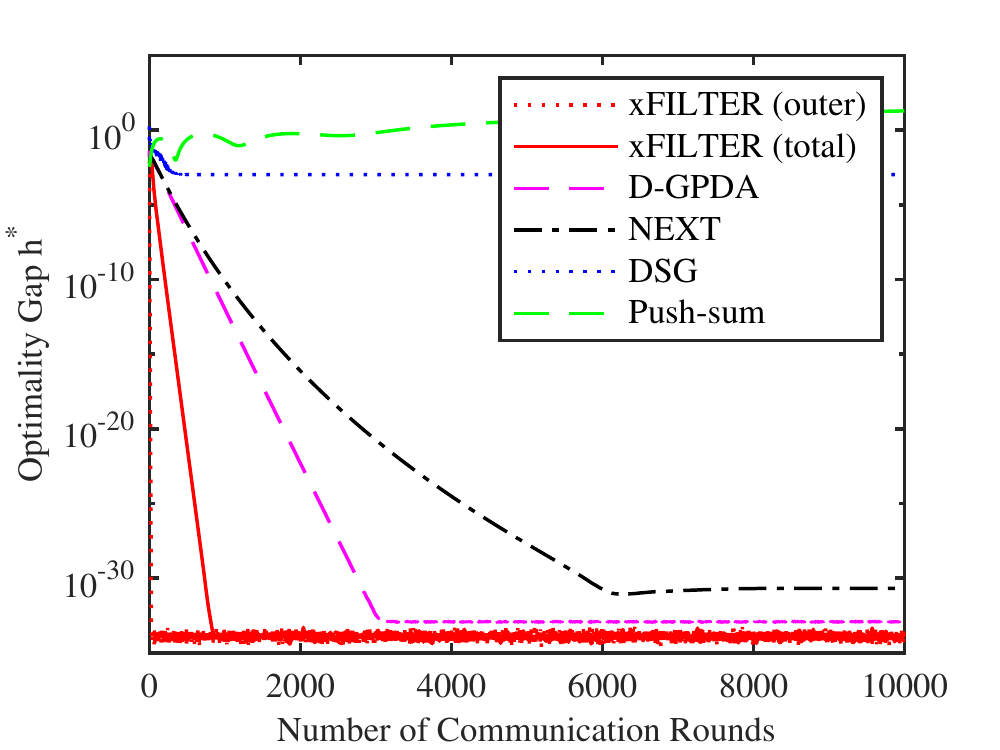}
 		\caption{$M=20, B=400, K=10$} \label{fig:ObjectiveComparison-6}
 	\end{minipage}
 \end{figure}

 \begin{figure}[t]
% 	\vspace{-0.5cm}
 	\begin{minipage}[c]{0.33\linewidth}
 		\includegraphics[width=\linewidth]{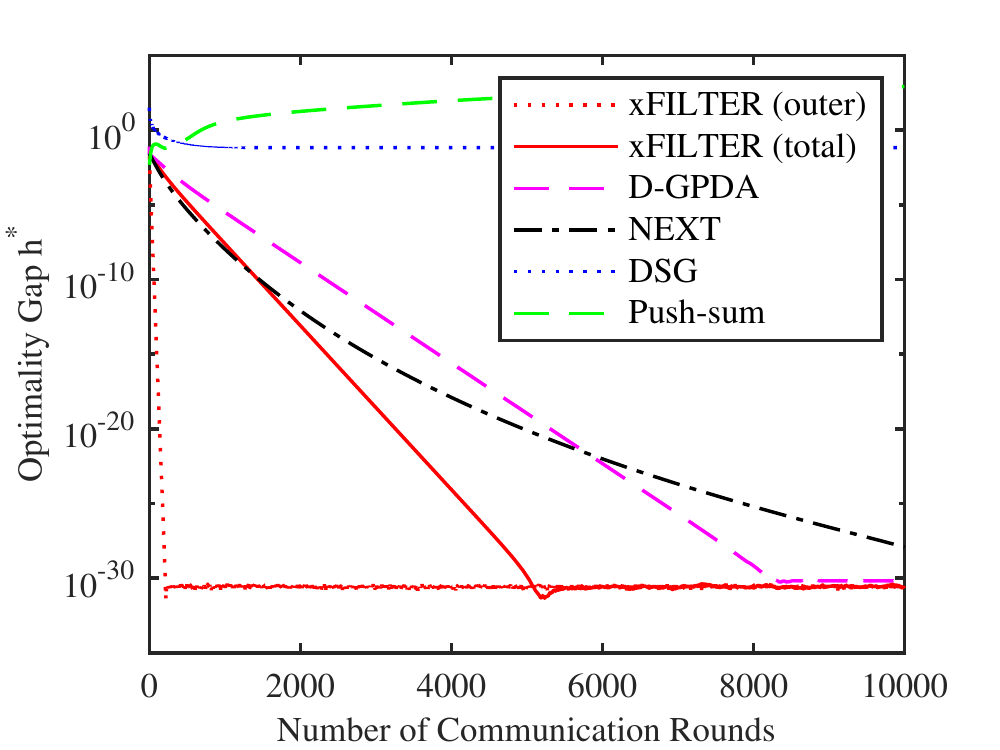}
 		\caption{$M=10, B=20, K=5$} \label{fig:ObjectiveComparison-7}
 	\end{minipage}% 
 	\hfill
 	\begin{minipage}[c]{0.33\linewidth}
 		\includegraphics[width=\linewidth]{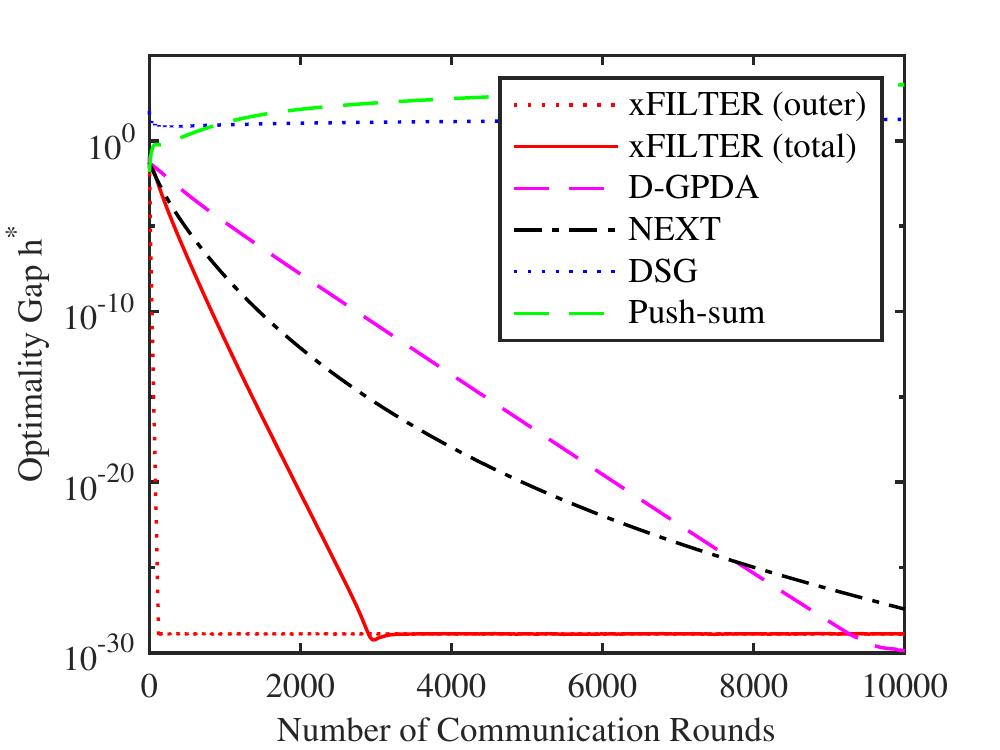}
 		\caption{$M=10, B=20, K=10$} \label{fig:ObjectiveComparison-8}
 	\end{minipage}
 	\hfill
 	\begin{minipage}[c]{0.33\linewidth}
 		\includegraphics[width=\linewidth]{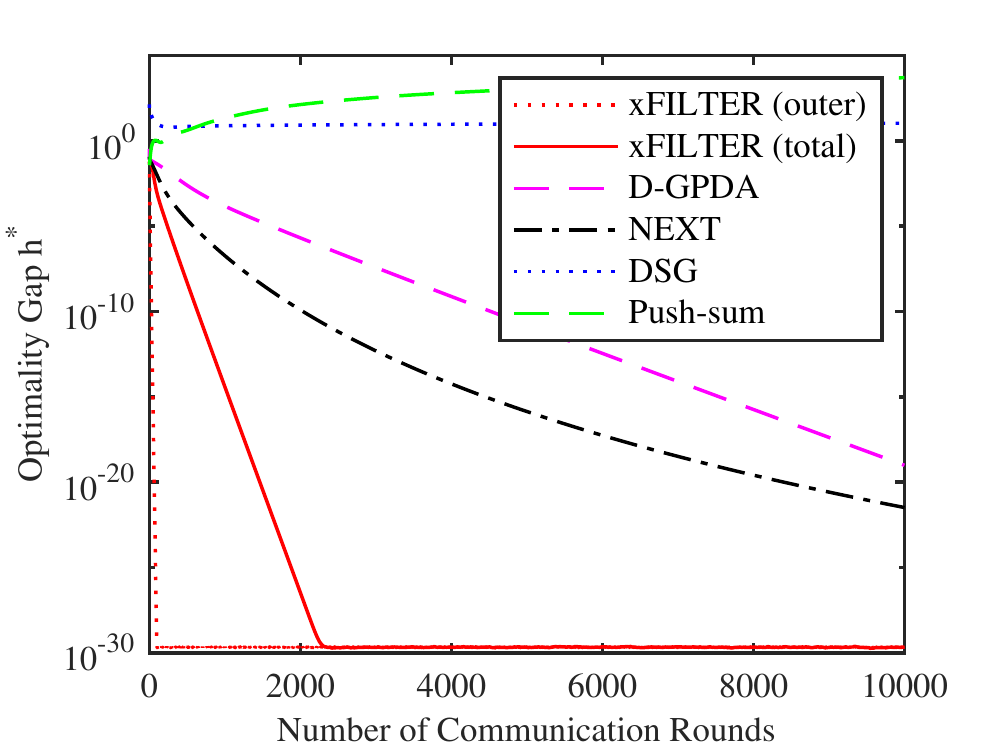}
 		\caption{$M=10, B=20, K=20$} \label{fig:ObjectiveComparison-9}
 	\end{minipage}
 \end{figure}
 
 %\newpage 
 \begin{figure}
% 	\vspace{-0.5cm}
% 	\centering
 	\begin{minipage}[c]{0.5\linewidth}
 		\includegraphics[width=\linewidth]{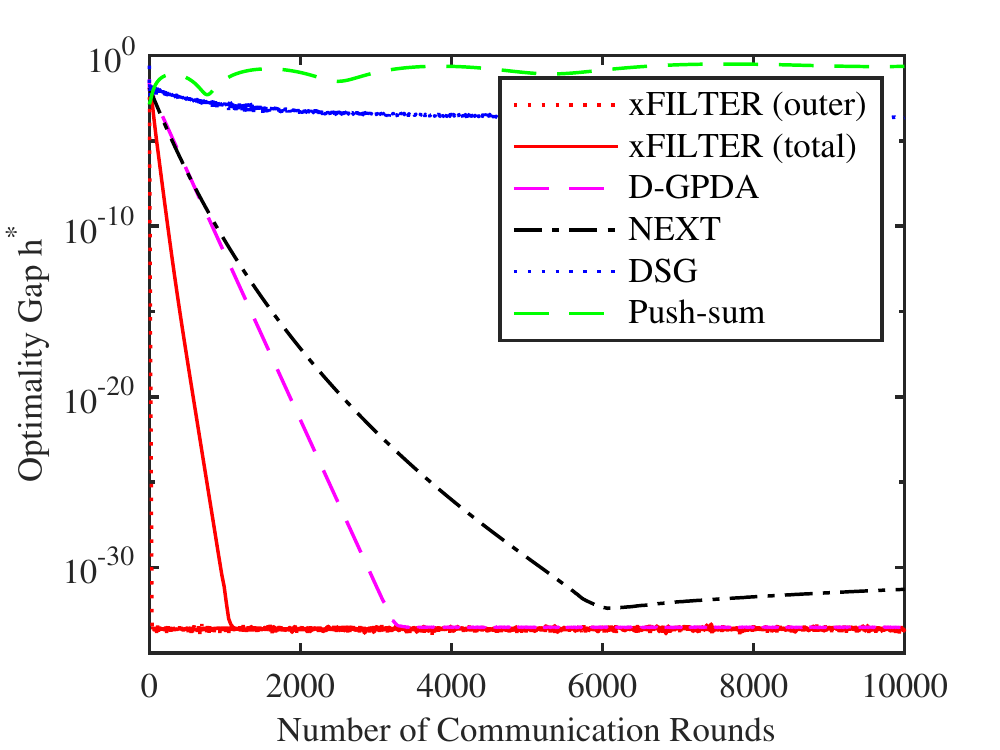}
 		\caption{$M=50, B=2000, K=10$} \label{fig:ObjectiveComparison-10}
 	\end{minipage}% 
 \hfill
 \begin{minipage}[c]{0.5\linewidth}
 	\includegraphics[width=\linewidth]{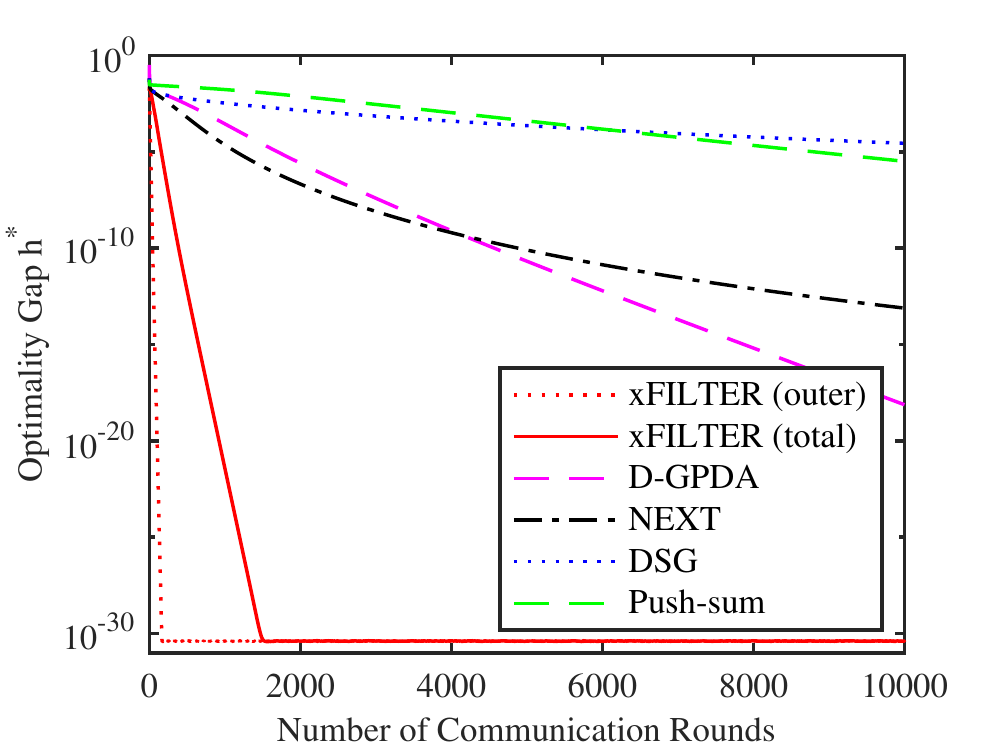}
 	\caption{$M=10, B=56, K=30$} \label{fig:ObjectiveComparison-real}
 \end{minipage}%
 	%\vspace{-2cm}
 \end{figure}

 To compare the convergence performance of the proposed algorithms, we randomly generated $MB$ data points with dimension $K$ and distribute them into $M$ nodes, i.e. each node contains $B$ data points with $K$ features. Then we compare the proposed xFILTER and D-GPDA with the distributed subgradient (DSG) method \cite{Nedic09subgradient}, the Push-sum algorithm \cite{tatarenko2017non}, and the NEXT algorithm \cite{Lorenzo16}. The parameters for NEXT are chosen as $\tau=1, \alpha[0]=0.1$ and $\mu =0.01$ as suggested by \cite{Lorenzo16}, while the parameters for xFILTER are chosen based on \eqref{eq:choice:uniform:2+}.

Simulation results on synthetic data for different $M, B, K$ averaged over $30$ realizations are investigated and shown in Fig. \ref{fig:ObjectiveComparison-1} to Fig. \ref{fig:ObjectiveComparison-10}, where the x-axis denotes the total rounds of communications required, and the y-axis denotes the quality measure  \eqref{eq:criteria:local} proposed in Section \ref{sub:prelim}. 
Note that the curves xFILTER (outer) included in these figures show the number of communication rounds required for xFILTER to perform the ``outer" iterations (which is equivalent to $r$ in Algorithm 2, since in each outer iteration only one round of communication is required in Step {\bf S3}).
	The performance evaluated on real data is also characterized in Fig. \ref{fig:ObjectiveComparison-real}, in which we choose $M=10$, $B=56$, and $K=30$.
	These results show that the proposed algorithms perform well in all parameter settings compared with existing methods.  
	
	We further note that these figures also show (rough) comparison about computation efficiency of different algorithms. Specifically, for D-GPDA, DSG and Push Sum (resp. NEXT), the total rounds of communication is the same as (resp. twice as) the total number of gradient evaluations per node. In contrast,  the total rounds of communication in the {\it outer loop} of xFILTER  is the same as the local gradient evaluations.  Therefore, the comparison between xFILTER (outer) and other algorithms in Fig. \ref{fig:ObjectiveComparison-3} to Fig. \ref{fig:ObjectiveComparison-10} shows the relative computational efficiency of these algorithms. 
		Clearly, xFILTER has a significant advantage over the rest of the algorithms.

Further, we compare the scalability performance of the proposed algorithms with increased network dimension $M$, and the results are shown in Fig. \ref{fig:ScalabilityComparison}, Table \ref{table:iter1} and Table \ref{table:iter2}. In particular, in Fig. \ref{fig:ScalabilityComparison} we compare the total communication rounds required for NEXT and the xFILTER for reaching $h^*_T\le 10^{-10}$ and $h^*_T\le 10^{-15}$, over path graphs with increasing number of nodes. Overall, we see that the xFILTER performs reasonably fast.

 \begin{figure}[t]
 	\centering
 	\begin{minipage}[c]{0.4\linewidth}
 		\includegraphics[width=\linewidth]{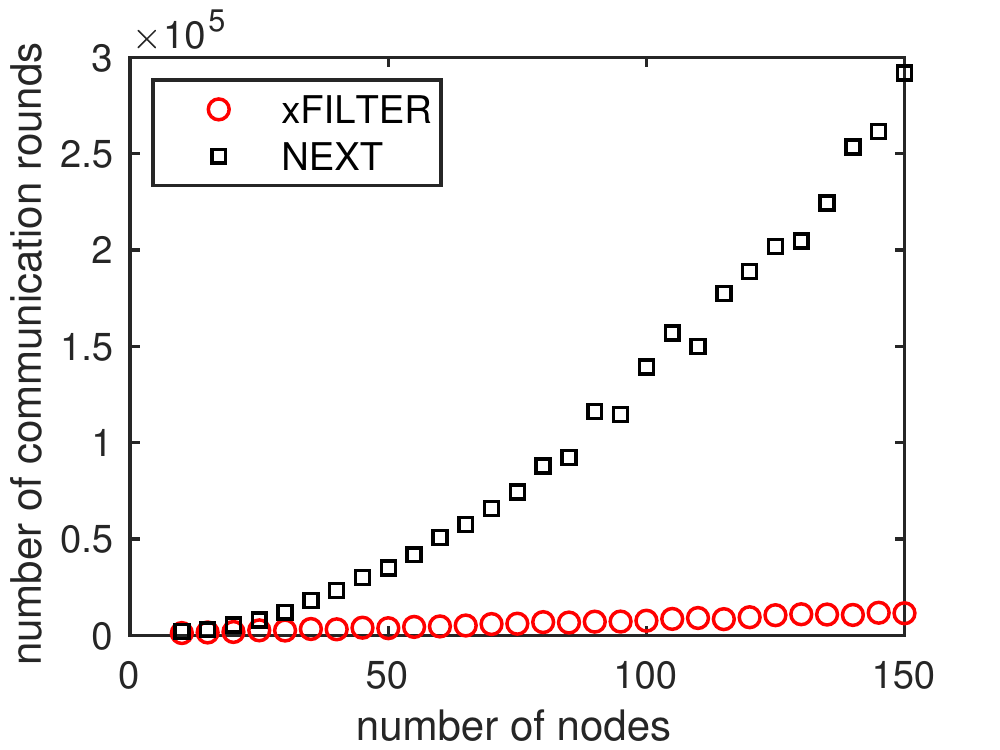}
 		\centering{{\footnotesize (a) $B=10, K=10, \epsilon=10^{-10}$}}  
 	\end{minipage}% 
 \begin{minipage}[c]{0.4\linewidth}
 	\includegraphics[width=\linewidth]{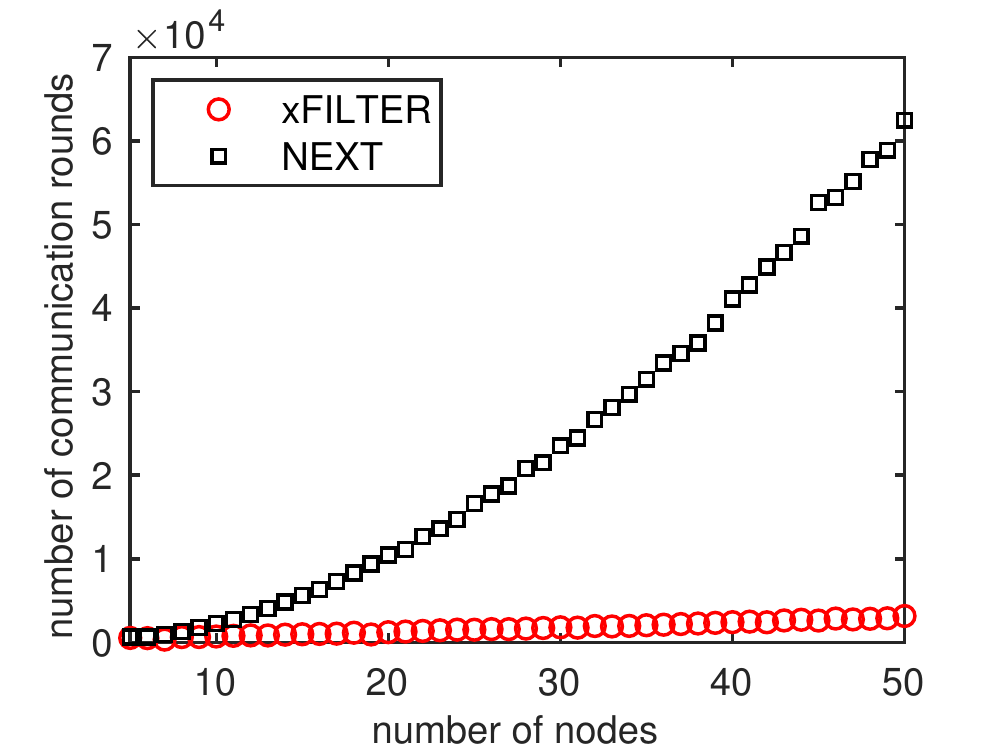}
 	\centering{{\footnotesize (b) $B=200, K=10, \epsilon= 10^{-15}$}}  
 \end{minipage}
 \caption{Comparison of NEXT and xFILTER over path graphs with increasing number of nodes ($M\in [10, \; 150]$ in (a) and $M\in [5, \; 50]$ in (b)). Each point in the figure represents the total number of communication needed to reach $h^*_T\le \epsilon$.  }\label{fig:ScalabilityComparison}
 \end{figure}

We do want to point out that although for the unconstrained problems that we have tested, our proposed algorithms compare relatively favorably with NEXT, NEXT can in fact handle a larger class of problems because it is designed for nonsmooth and constrained nonconvex problems. Further, for all the algorithms we have used, we did not tune the parameters: For xFILTER and D-GPDA, we use the theoretical upper bound suggested in Theorem \ref{thm:final:bound}, and for NEXT we use the parameters suggested in the paper   \cite{Lorenzo16}. For all our tested problems and algorithms, it is possible to fine-tune the stepsizes to make them faster, but since this paper is mostly on the theoretical properties of rate optimal algorithms, we choose not to go down that path.

\begin{table}[t]
	\tabcolsep 0pt 
	\caption{Optimality gap after $200$ rounds of communications $(B=200, K=10)$}
	\begin{center}
		\vspace{-0.3cm}
%		\small
		\def\temptablewidth{1\linewidth}
		{\rule{\temptablewidth}{1pt}}
		\begin{tabular*}{\temptablewidth}{@{\extracolsep{\fill}}c|cccc}
			{number of nodes $M$ } &D-GPDA    &xFILTER   \\ \hline \hline
			10   & $3.96\times 10^{-4}$ & $2.50 \times 10^{-11}$\\
			20  & $5.45\times 10^{-4}$ & $1.92\times 10^{-9}$ \\
			30 & $1.20\times 10^{-4}$  & $4.71\times 10^{-11}$  \\
			40 & $2.95\times 10^{-4}$  & $4.07\times 10^{-10}$  \\
			50  & $3.88\times 10^{-4}$ & $8.47\times 10^{-11}$\\
		\end{tabular*}
		{\rule{\temptablewidth}{1pt}}
%		\vspace{-1cm}
	\end{center}
	\label{table:iter1}
\end{table}

 \begin{table}[t]
 	\tabcolsep 0pt 
 	\caption{Optimality gap after $1000$ rounds of communications $(B=200, K=10)$}
 	\begin{center}
 		\vspace{-0.3cm}
% 		\small
 		\def\temptablewidth{1\linewidth}
 		{\rule{\temptablewidth}{1pt}}
 		\begin{tabular*}{\temptablewidth}{@{\extracolsep{\fill}}c|cccc}
 			{number of nodes $M$ } &D-GPDA    &xFILTER   \\ \hline \hline
 			10   & $8.24\times 10^{-13}$ & $1.93 \times 10^{-33}$\\
 			20  & $9.41\times 10^{-12}$ & $1.43\times 10^{-32}$ \\
 			30 & $2.09\times 10^{-13}$  & $2.26\times 10^{-32}$  \\
 			40 & $1.52\times 10^{-11}$  & $4.19\times 10^{-33}$  \\
 			50  & $2.30\times 10^{-10}$ & $6.48\times 10^{-33}$\\
 		\end{tabular*}
 		{\rule{\temptablewidth}{1pt}}
% 		\vspace{-1cm}
 	\end{center}
 	\label{table:iter2}
 \end{table}

  \subsection{Distributed Neural Network Training}
	In our second experiment, we present some numerical results under a more realistic setting. We consider training a neural network model for fitting the  MNIST data set. The dataset is first randomly partitioned into $10$ subsets, and then gets distributed over $10$ machines.  A fully connected neural network with one hidden layer is used in the experiment. The number of neurons for the hidden layer and the output layer are set as $128$ and $10$, respectively.  The initial weights for the neural network are drawn from a truncated normal distribution centered at zero with variance scaled with the number of input units.
		%, which contains only one out of ten categories. 
		The algorithms are written in Python, and the communication protocol is implemented using the Message Passing Interface (MPI). The empirical performance of the xFILTER is evaluated and compared with the DSG algorithm \cite{nedic2009distributed}. Fig. \ref{fig:DNN} shows that, compared with DSG, the proposed algorithm achieves better communication and computation efficiency, and has improved classification accuracy.
	
	Note that despite the fact that some global parameters (such as the Lipschitz constants) are unknown, the rules provided in \eqref{eq:choice:uniform:2+} or \eqref{eq:choice:uniform:2+:2} still can help us  roughly estimate a set of good parameters. For example, we choose the following parameters 
	\begin{align}
	\Sigma^2= \frac{\sigma}{  \sum_{i}d_i \underline{\lambda}_{\min}(\tcL)} , \quad \Upsilon^2 = \frac{\beta P}{\sum_{i}d_i},
	\end{align}
	and tune the parameter $\beta$ and $\sigma$  by searching from the set  $\{0.1, 0.2, 0.5, 1, 2, 5, \cdots, 100, 200, 500\}$.
	Based on the best practical performance over $10$ runs, we choose  $\beta=100$ and $\sigma=20$ for xFILTER and $\alpha = 0.1$ for DSG.
	
	 \begin{figure}
		\centering
		\begin{minipage}[c]{0.4\linewidth}
			\includegraphics[width=\linewidth]{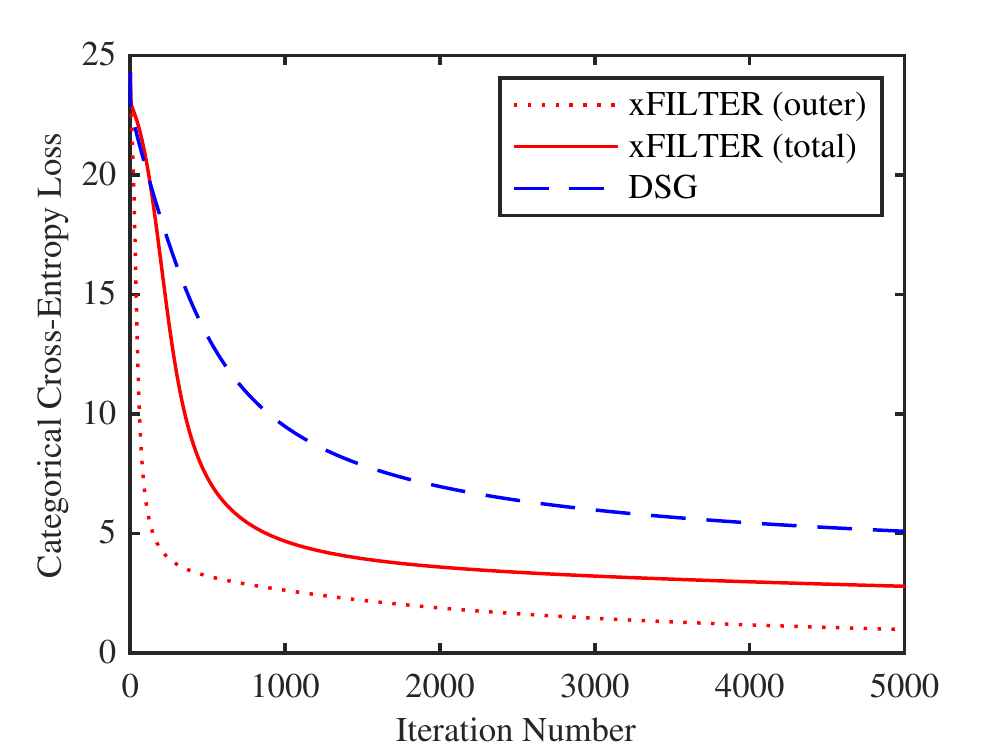}
			\centering{{\footnotesize (a) Training Loss}}  
		\end{minipage}% 
		\begin{minipage}[c]{0.4\linewidth}
			\includegraphics[width=\linewidth]{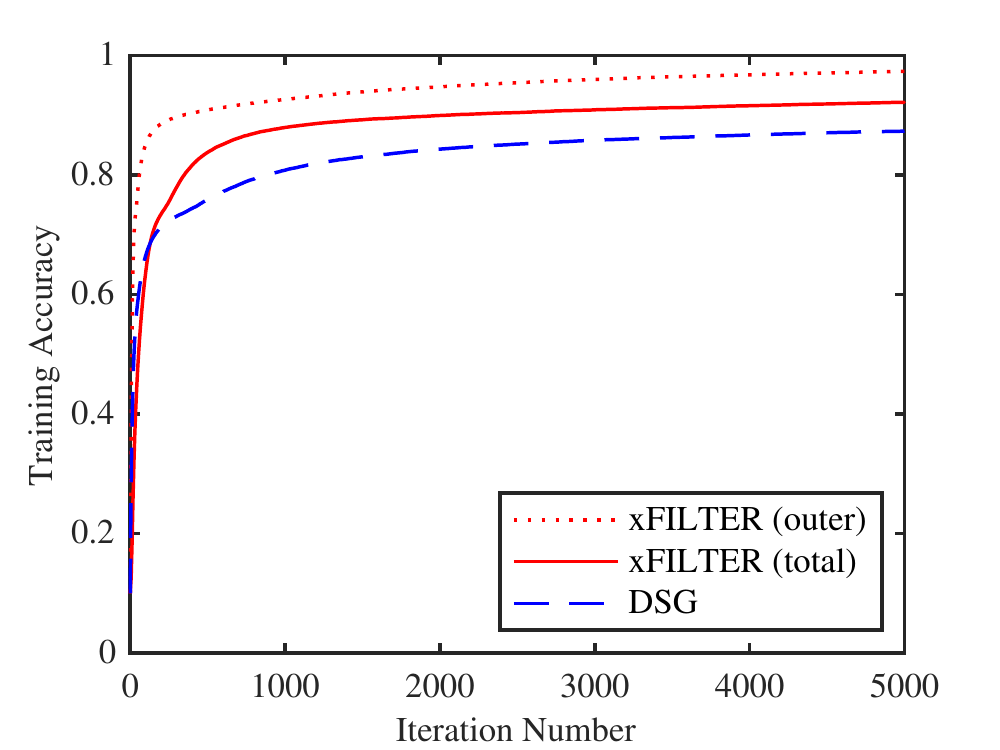}
			\centering{{\footnotesize (b) Training Accuracy}}  
		\end{minipage}
			\caption{Comparison of DSG and xFILTER over path graphs on distributed training neural networks; Plot (a) shows the dynamic of the categorical cross-entropy loss, and plot (b) shows the training classification accuracy. The parameters are chosen based on their best practical performance through grid search. The curves xFILTER (outer)  and xFILTER (total) again represent the number of outer iteration, and the total number of iterations required for xFILTER.}\label{fig:DNN}
%		\vspace{-0.5cm}
	\end{figure}

\section{Conclusion and Future Works}
\begin{figure}[t]
	\centering
	\begin{minipage}[c]{0.7\linewidth}
		\includegraphics[width=\linewidth]{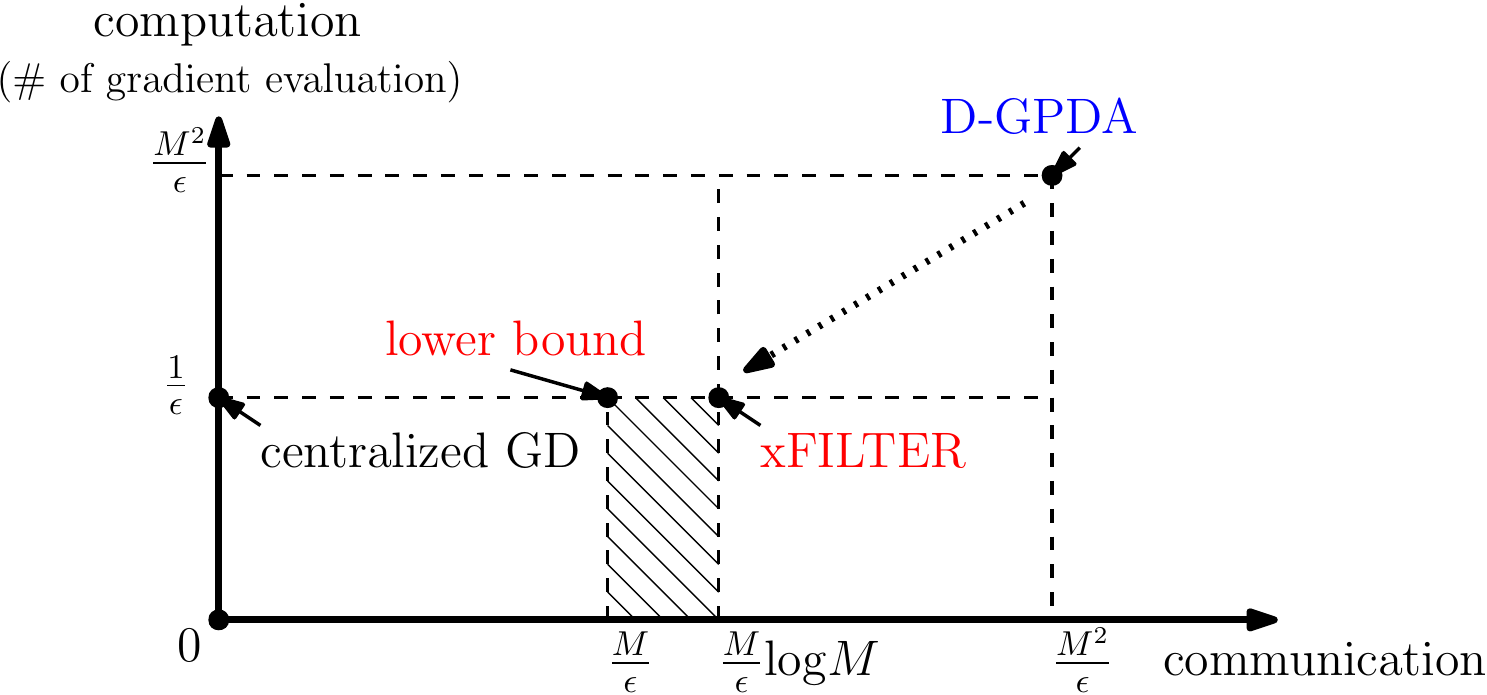}
		\caption{Graphical comparison of various bounds analyzed in this work, illustrated over a path graph with $M$ nodes. }\label{fig:comlexity}
	\end{minipage}
\end{figure}

This paper represents the first work that investigates the performance of optimal first-order non-convex  algorithms for distributed information processing and optimization problems. We first set our scope by defining the problem, network, and algorithm classes $(\cP, \cN, \cA)$ that are under consideration. We then provide a lower complexity bound that characterizes the worst case performance for any first-order distributed algorithm in class $\cA$, and finally propose and analyze two algorithms that are capable of (nearly) achieving the lower bound in various settings.  The various bounds discussed in the work is  illustrated in Fig. \ref{fig:comlexity} through a $M$-node path graph as an example.

To the best of our knowledge, the proposed algorithms are the first and the only available distributed non-convex algorithms in class $\cA$ that can optimally reduce {\it both} the size of the gradient and the consensus error for $(\cP,\cN)$, and achieving the (near) optimal rate performance for problem/network classes $(\cP, \cN)$. However, they still require some global information to initialize the parameters, so it will be of interest to design global information free algorithms that only require local structures to set parameters (just like in the convex case, see discussions in \cite{chang14distributed}). It will also be desirable to merge the inner Chebyshev iteration with the outer dual update to design a {\it single-loop} algorithm, and to extend the proposed algorithms to problems with nonsmooth regularizers and constraints.

\section{Appendix}

\subsection{Proof of Lemma \ref{lm:dual:different}}
\noindent{\bf Proof.}
First we show that for all $r\ge -1$ the following holds for D-GPDA {
	\begin{align}\label{eq:opt}
	\nabla f(x^r) + F^T \lambda^{r} + F^T \Sigma^2 F x^{r+1} + H(x^{r+1}-x^r) = 0.
	\end{align}}
Note that for the initialization \eqref{eq:initialization} we have{
\begin{align*}
&\nabla f(x^{-1}) + (2\Delta + \Upsilon^2) x^0 = 	\nabla f(x^{-1}) + (F^T \Sigma^2 F + H) x^0 = 0.
\end{align*}}
Setting $x^{-1} =0, \lambda^{-1} =0$ and using \eqref{eq:initialization}, we obtain
\begin{align}\label{eq:opt:init}
\nabla f(0) + F^T \lambda^{-1} + F^T \Sigma^2 F x^{0} + H(x^{0}-x^{-1}) = 0.
\end{align}
Further, the optimality condition of the $x$ update \eqref{eq:x:update} suggests that \eqref{eq:opt} holds for all $r\ge 0$, therefore \eqref{eq:opt} is proved.

Second, by using \eqref{eq:opt} and the $y$ update \eqref{eq:mu:update}, we obtain
\begin{align}\label{eq:opt:new:lambda}
F^T \lambda^{r+1} = - \nabla f(x^r) -H(x^{r+1}-x^r), \; \forall~r\ge -1.
\end{align}
Then subtracting the previous iteration leads to{
\begin{align*}
&F^T (\lambda^{r+1}-\lambda^r) = - (\nabla f(x^r) - \nabla f(x^{r-1})) - Hw^{r+1}, \; \forall~r\ge 0.
\end{align*}}
 
Note that the matrix $H\succ 0$, $\Sigma^2\succ0$, then we have 
\begin{align}\label{eq:key}
&H^{-1/2}(\Sigma F)^T  \Sigma^{-1}(\lambda^{r+1}-\lambda^r) = - H^{-1/2} (\nabla f(x^r) - \nabla f(x^{r-1})) -H^{1/2}w^{r+1}.
\end{align}}
Then using the fact that 
$$\Sigma^{-1}(\lambda^{r+1}-\lambda^r) = \Sigma F x^{r+1}\in\mbox{col}(\Sigma F),$$ 
we can square both sides and obtain the following
\begin{align}\label{eq:key:2}
& \underline{\lambda}_{\min}(\Sigma F H^{-1}F^T\Sigma )\|\Sigma^{-1}(\lambda^{r+1}-\lambda^r)\|^2 \nonumber\\
&\le  2 \|H^{-1/2} (\nabla f(x^r) - \nabla f(x^{r-1}))\|^2 + 2  (w^{r+1})^TH w^{r+1}\nonumber\\
&\le  2 \|\Upsilon^{-1}(\nabla f(x^r) - \nabla f(x^{r-1}))\|^2+ 2  (w^{r+1})^T H w^{r+1}\nonumber\\
&\stackrel{\eqref{eq:Lip:extended}}\le \frac{2}{M^2}  \|\Upsilon^{-1} L (x^r-x^{r-1})\|^2+  2  \|w^{r+1}\|^2_{H}, \; \forall~r\ge 0.
\end{align}
This concludes the proof of the first part. 

To show the second part, note that according to \eqref{eq:linear:cheby:error:equiv}, $x^{r+1}$ generated by xFILTER is given by (for all $r\ge -1$){
\begin{align}\label{eq:opt:gpda+}
\hspace{-0.4cm}\nabla f(x^r) + F^T (\lambda^r + \Sigma^2 F x^{r+1}) + \Upsilon^2(x^{r+1}-x^r) = \Upsilon^2 R\epsilon^{r+1}.
\end{align}}
Then use the same analysis steps as in the first part, we arrive at the desired result. 
\QED

%\vspace{-0.5cm}
\subsection{Proof of Lemma \ref{lm:AL:change}}
\noindent{\bf Proof.}	Using the Lipschitz gradient assumption \eqref{eq:Lip:extended}, we have{
	\begin{align}
	&{\sf{AL}}(x^{r+1},\lambda^r)- {\sf{AL}}(x^r,\lambda^r) \le \langle \nabla f(x^r) + f^T \lambda^r+ F^T \Sigma^2F  x^r, x^{r+1}-x^r \rangle\nonumber\\
	&\quad +\frac{1}{2M} \|x^{r+1}-x^r\|_L^2 +\frac{1}{2}\|\Sigma F (x^{r+1}-x^r)\|^2\nonumber\\
	&=  \langle \nabla f(x^r) + F^T \lambda^r+ A^T \Sigma^2F  x^{r+1}, x^{r+1}-x^r\rangle \nonumber\\
	&\quad + \langle H (x^{r+1}-x^r), x^{r+1}-x^r \rangle  +\frac{1}{2M} \|x^{r+1}-x^r\|_L^2 \nonumber\\
	&\quad +\frac{1}{2}\|\Sigma F (x^{r+1}-x^r)\|^2 - \|x^{r+1}-x^r\|_{H + F^T \Sigma^2 F} \nonumber\\
	&\stackrel{\eqref{eq:opt},\eqref{eq:sum:A:B}}\le - (x^{r+1}-x^r)^T \left(\frac{\Delta}{2}-\frac{L}{2 M}+\Upsilon^2\right) (x^{r+1}-x^r). %\quad \forall~r\ge 0.
%	\vspace{-1cm}
	\end{align}}
	Using the update rule of the dual variable, and combine the above inequality, we obtain
	\begin{align}
	\sf{AL}^{r+1} - \sf{AL}^r &\le - \frac{1}{2}\|x^{r+1}-x^r\|^2_{\Delta+2\Upsilon^2 -  L/M} + \langle \lambda^{r+1}-\lambda^r, A x^{r+1}\rangle \nonumber\\
	& =  -\frac{1}{2}\|x^{r+1}-x^r\|^2_{\Delta+2\Upsilon^2 -  L/M} + \|\Sigma^{-1}(\lambda^{r+1}-\lambda^r)\|^2\nonumber
%	& \le -\frac{1}{2}\|x^{r+1}-x^r\|^2_{\Delta+2\Upsilon^2 - L/M} + \kappa \left(\frac{2}{M^2} \|\Upsilon^{-1} L (x^r-x^{r-1})\|^2+  2  \|w^{r+1}\|^2_{H}\right).
	\end{align}
Combined with Lemma \ref{lm:dual:different} we complete the first part. 

The second part follows similar steps. The modifications are that $H$ is replaced by $\Upsilon^2$, and that there is an additional error term in the optimality condition; cf. \eqref{eq:linear:cheby:error:equiv}. \QED	

%\vspace{-0.3cm}
\subsection{Proof of Lemma \ref{lm:Potential} and Lemma \ref{lm:Potential+}}
\noindent{\bf Proof.} Using the optimality condition \eqref{eq:opt}, we have
\begin{align*}
\langle F^T \lambda^{r+1}  + \nabla f(x^r) +H(x^{r+1}-x^r), x^{r+1}- x\rangle  =0,\; \forall~r\ge -1
\end{align*}
This implies that for all $r\ge 0$
\begin{align*}
\langle F^T (\lambda^{r+1}-\lambda^r)  + \nabla f(x^r) -\nabla f(x^{r-1})+ H w^{r+1}, x^{r+1}- x^{r}\rangle  =0.
\end{align*}
It follows that
\begin{align}\label{eq:second:change}
&\frac{1}{2}\|\Sigma Fx^{r+1}\|^2 +\frac{1}{2}\|x^{r+1}-x^r\|_{H}^2 \le \frac{1}{2}\|\Sigma  F x^r\|^2 + \frac{1}{2}\|x^r-x^{r-1}\|_{H} - \frac{1}{2}\|w^{r+1}\|^2_{H} \\\
&\quad + \frac{1}{2 {M}}\|x^{r+1}-x^r\|^2_{L}+ \frac{1}{2{M}}\|x^{r}-x^{r-1}\|^2_{L}, \quad \forall~r\ge 0.\nonumber
\end{align}
Then combining Lemma \ref{lm:AL:change} and \eqref{eq:second:change},  for all $r\ge 0$ we have{
\begin{align}
\hspace{-0.2cm}&P^{r+1}- P^r \le - \left(\frac{c}{2}- 2\kappa \right) \|w^{r+1}\|_{H}^2\nonumber\\
\hspace{-0.2cm}& -\frac{1}{2}(x^{r+1}-x^r)^T \bigg(\Delta+2 \Upsilon^2 -  \frac{L}{M} - \frac{4\kappa}{M^2}L \Upsilon^{-2}L  - \frac{ 2 c L}{M} \bigg) (x^{r+1}-x^r) \nonumber.
\end{align}}
Therefore, in order to make the potential function decrease, we need to follow \eqref{eq:constants}. 

To show a similar result for the xFILTER, consider the following optimality condition derived from \eqref{eq:opt:gpda+} {
\begin{align*}
\langle F^T \lambda^{r+1}  + \nabla f(x^r) +\Upsilon^2(x^{r+1}-x^r) - \Upsilon^2 R\epsilon^{r+1}, x^{r+1}- x\rangle  =0, \; \forall~x.
\end{align*}}
Following similar steps as in  \eqref{eq:second:change}, and use \eqref{eq:epsilon:bound}, we have
\begin{align}\label{eq:second:change+}
\frac{1}{2}\|\Sigma Fx^{r+1}\|^2 +\frac{1}{2}\|x^{r+1}-x^r\|_{\Upsilon^2}^2 &\le\frac{1}{2} \|\Sigma  F x^r\|^2 + \frac{1}{2}\|x^r-x^{r-1}\|_{\Upsilon^2} - \frac{1}{2}\|w^{r+1}\|^2_{\Upsilon^2} \\
&\quad + {1}/(2M)\|x^{r+1}-x^r\|^2_{L}+ {1}/(2M)\|x^{r}-x^{r-1}\|^2_{L} \nonumber\\
&\quad + 1/4\|x^{r+1}-x^r\|_{\Upsilon^2  R}^2 + 1/4\|x^r-x^{r-1}\|_{\Upsilon^2  R}^2, \quad \forall~r\ge -1.\nonumber
\end{align}
Then combining Lemma \ref{lm:AL:change}, \eqref{eq:second:change}, and the estimate of the size of $\epsilon$ in \eqref{eq:epsilon:bound}, we have{
	\begin{align*}
&\widetilde{P}^{r+1}- \widetilde{P}^r\le -\frac{1}{2}(x^{r+1}-x^r)^T  V (x^{r+1}-x^r) - \left(\frac{\widetilde{c}}{2}- 3\widetilde{\kappa} \right) \|w^{r+1}\|_{\Upsilon^2}^2.
	\end{align*}}
with
\begin{align*}
V:= \bigg(\Upsilon^2 R -  (1+2\tilde{c})\frac{L}{M} - \frac{6\tilde{\kappa}}{M^2}L \Upsilon^{-2}L  - \frac{\Upsilon^2 R(24\tilde{\kappa}+6+16\tilde{c})}{16} \bigg).
\end{align*}
Therefore in order to make the potential function decrease, we need to follow \eqref{eq:constants+}. \QED	

\subsection{Proof of Lemma \ref{lm:bounded}}
\noindent{\bf Proof.} 
For D-GPDA, we can express the AL as (for all $r\ge 0$) {
\begin{align*}
\hspace{-0.3cm}&{\sf{AL}}^{r+1}- f(x^{r+1}) %= \langle \lambda^{r+1}, A x^{r+1}\rangle +\frac{1}{2}\|\Sigma A x^{r+1}\|^2\nonumber\\
=\langle \lambda^{r+1}, \Sigma^{-2}(\lambda^{r+1}-\lambda^r)\rangle +\frac{1}{2}\|\Sigma F x^{r+1}\|^2\nonumber\\
\hspace{-0.3cm}&=\frac{1}{2}\hspace{-0.1cm}\left(\|\Sigma^{-1}\lambda^{r+1}\|^2\hspace{-0.1cm}-\hspace{-0.1cm}\|\Sigma^{-1}\lambda^r\|^2\hspace{-0.1cm}+\|\Sigma^{-1}(\lambda^{r+1}-\lambda^r)\|^2+\|\Sigma F x^{r+1}\|^2\right).
\end{align*}}
Since $\inf_{x}f(x)=\underline{f}$ is lower bounded, let us define 
\begin{align*}
\widehat{\sf{AL}}^{r+1}\hspace{-0.4cm}:= {\sf{AL}}^{r+1} - \underline{f}, \;  \widehat{f}(x)  := f(x)- \underline{f}\ge 0, \;\widehat{P}^{r+1} \hspace{-0.2cm} := P^{r+1}- \underline{f}.
\end{align*}
Therefore, summing over $r=-1\cdots, T$, we obtain
\begin{align*}
&\sum_{r=-1}^{T}\widehat{\sf{AL}}^{r+1} = \frac{1}{2}\left(\|\Sigma^{-1}\lambda^{T+1}\|^2 -\|\Sigma^{-1}\lambda^{-1}\|^2\right)\nonumber\\
&+\sum_{r=-1}^{T}\left(\widehat{f}(x^{r+1})+\frac{1}{2}\|\Sigma F x^{r+1}\|^2 + \frac{1}{2}\|\Sigma^{-1}(\lambda^{r+1}-\lambda^r)\|^2\right) \nonumber.
\end{align*}
Using the initialization $\lambda^{-1}=0$, then the above sum is lower bounded by zero. %since   
This fact implies that the sum of $\widehat{P}^{r+1}$ is also lower bounded by zero (note, the remaining terms in the potential function are all nonnegative)
\begin{align*}
\sum_{r=0}^{T}\widehat{P}^{r+1} \ge 0, \quad \forall~T>0, 
\end{align*}
Note that if the parameters of the system are chosen according to \eqref{eq:constants}, then $P^{r+1}$ is nonincreasing, which implies that its shifted version $\widehat{P}^{r+1}$ is also nonincreasing. Combined with the nonnegativity of the sum of the shifted potential function, we can conclude that 
\begin{align}
\widehat{P}^{r+1} \ge 0,\quad \mbox{and}\quad P^{r+1}\ge \inf f(x), \quad \forall~r\ge 0.
\end{align}
Next we compute $P^0$. By using \eqref{eq:sum:A:B}, we have 
\begin{align}\label{eq:P0}
\hspace{-0.3cm}P^0
%&={\sf{AL}}(x^{0},\lambda^{0})+ \frac{2\kappa}{M^2} \|  \Upsilon^{-1}  L x^{0}\|^2 + \frac{c}{2}\left(\|\Sigma Ax^{0}\|^2 +\|x^{0}\|_{H}^2+ \frac{1}{ {M}}\|x^{0}\|^2_{L}\right)\nonumber\\
&= {\sf{AL}}^0+ \frac{2\kappa}{M^2} \|  \Upsilon^{-1}  L x^{0}\|^2 + \frac{c}{2}\left(\|x^{0}\|_{2\Delta+ \Upsilon^2 +L/M}^2\right)\\
{\sf{AL}}^0 & \le  f(x^0) + 2 \|\Sigma F x^0\|^2\nonumber\\
x^0 & \stackrel{\eqref{eq:initialization}}= 	(2\Delta + \Upsilon^2)^{-1} \frac{1}{M}[\nabla f_1(0);  \cdots; \nabla f_M(0)]\nonumber\\
& =(2\Delta + \Upsilon^2)^{-1} \frac{1}{M}d_0
\end{align}
where in the last equality we have used the definition of $d_0$ in \eqref{eq:d0}. 
Use the above relations, we have
\begin{align}
P^0 \le f(x^0) + (x^{0})^T Z x^0
\end{align}
%To express  $P_{c}(x^{0},x^{-1},\lambda^{0})$ in terms of $d_0$, we proceed as follows. First we note that the coefficients in $\|x_0\|^2$ is given by
with the matrix $Z$ defined as{
\begin{align*}
Z&= \frac{2\kappa}{M^2} L \Upsilon^{-2} L +\frac{c(2\Delta+\Upsilon^2+L/M)}{2} + 2 F\Sigma^2 F \preceq 4c (2\Delta+\Upsilon^2)
\end{align*}}
where the last inequality follows from our choice of parameters in \eqref{eq:constants:choices}, and the fact $c\ge 1$.
Note that{
\begin{align}
&4c \|x^0\|^2_{2\Delta+\Upsilon^2} \le \frac{4 c}{M^2} d^T_0  (2\Delta + \Upsilon^2)^{-1}(2\Delta + \Upsilon^2) (2\Delta + \Upsilon^2)^{-1} d_0\nonumber\\
&= \frac{4 c}{ M^2} d^T_0  (2\Delta + \Upsilon^2)^{-1} d_0\le \frac{2}{M} d_0^T L^{-1} d_0 %\label{eq:c1}
\end{align}}
where the last inequality comes from the choice of the parameters \eqref{eq:constants:choices}, which implies that $ 2\Delta+ \Upsilon^2\succeq  2\frac{L c}{M}$.
These constants combined with \eqref{eq:P0} shows the desired result. 

For xFILTER, the proof for the lower boundedness is the same. To bound the size of $\widetilde{P}^0$, first note that we again have
 {
	\begin{align*}
	\hspace{-0.3cm}&{\sf{AL}}^{r+1}- f(x^{r+1}) =\frac{1}{2}\hspace{-0.1cm}\left(\|\Sigma^{-1}\lambda^{r+1}\|^2\hspace{-0.1cm}-\hspace{-0.1cm}\|\Sigma^{-1}\lambda^r\|^2\hspace{-0.1cm}+\|\Sigma^{-1}(\lambda^{r+1}-\lambda^r)\|^2+\|\Sigma F x^{r+1}\|^2\right).
	\end{align*}}
By letting $r=-1$, and use the fact that $x^{-1}=0$ and $\lambda^{-1}=0$, we obtain
	\begin{align}\label{eq:al0}
\hspace{-0.3cm}&{\sf{AL}}^{0}- f(x^{0}) =\frac{1}{2}\hspace{-0.1cm}\left(2\|\Sigma^{-1}\lambda^{0}\|^2+\|\Sigma F x^{0}\|^2\right) = \frac{3}{2}\|\Sigma^{-1}\lambda^{0}\|^2 .
\end{align}}

Then we have
	\begin{align}
	&\widetilde{P}^0= {\sf{AL}}^0+ \frac{3\widetilde{\kappa}}{M^2} \|  \Upsilon^{-1}  L x^{0}\|^2 + \frac{3}{8} \tilde{\kappa}\|x^0\|^2_{\Upsilon^2 R }\nonumber\\
	&\quad + \frac{\widetilde{c}}{2}\left(\|\Sigma F x^{0}\|^2+ \|x^{0}\|^2_{\Upsilon^2+\Upsilon^2 R/4+L/M}\right),\label{eq:P0+}\\
	&{\sf{AL}}^0 \le  f(x^0) + 2 \|\Sigma F x^0\|^2, \; x^{-1}=0, \; \lambda^{-1}=0,\\
	&x^0 \stackrel{\eqref{eq:linear:cheby:error:equiv}}= 	R^{-1}\Upsilon^{-2}\nabla f(0) - \epsilon^0, \quad \widetilde{\epsilon}^{-1} \stackrel{\eqref{eq:epsilon:tilde}}= R^{-1}\Upsilon^{-2}\nabla f(0). \label{eq:epsilon:tilde:0}
	%&= (2\Delta + \Upsilon^2)^{-1} \frac{1}{M}[\nabla f_1(0); \nabla f_2(0); \cdots; \nabla f_M(0)] \nonumber\\
%	: =(2\Delta + \Upsilon^2)^{-1} \frac{1}{M}d_0.
	\end{align}
Use the above relation, we have
\begin{align}
\widetilde{P}^0\le f(x^0) + (x^{0})^T \widetilde{Z} x^0
\end{align}
%To express  $P_{c}(x^{0},x^{-1},\lambda^{0})$ in terms of $d_0$, we proceed as follows. First we note that the coefficients in $\|x_0\|^2$ is given by
with the matrix $Z$ defined as
	\begin{align*}
	\widetilde{Z}&= \frac{3\tilde{\kappa}}{M^2} L \Upsilon^{-2} L + \bigg(\frac{3}{8} \tilde{\kappa} + \tilde c\bigg)\Upsilon^2 R + \frac{\tilde c L}{2M} + 2 F\Sigma^2 F \preceq 3\Upsilon^2 R
	\end{align*}
where the last inequality follows from our choice of parameters in \eqref{eq:constants:choices+}. Therefore we have
\begin{align}
(x^{0})^T \widetilde{Z} x^0 &\le 3 (x^{0})^T \Upsilon^2  R x^0 \nonumber\\
& \le 3 (\nabla f(0) - \Upsilon^2  R \epsilon^0)^T R^{-1} \Upsilon^{-2}(\nabla f(0) - \Upsilon^2  R \epsilon^0)\nonumber\\
&\stackrel{\rm (i)}\le 6 (\nabla f(0))^T  R^{-1} \Upsilon^{-2}  \nabla f(0) + 6 (\epsilon^0)^T \Upsilon^2  R \epsilon^0\nonumber\\
&\le 3M (\nabla f(0))^T  L^{-1} \nabla f(0) + \frac{3}{8M}\|x^0\|_{\Upsilon^2  R}^2 %\epsilon^0%6\eta \xi^{-2}/\lambda_{\min}(R) \|R \tilde{\epsilon}^{-1}\|^2\nonumber\\
%&\stackrel{\eqref{eq:epsilon:bound}\le \frac{3}{M} d_0^T  L^{-1} d_0 +\frac{6 \eta}{\xi^2\lambda_{\min}(R)} \|\nabla f(0)\|^2 \nonumber.
\end{align}
where in ${\rm (i)}$ we have used the Cauchy-Swartz inequality; the last inequality uses \eqref{eq:epsilon:bound}, the choice of the parameters \eqref{eq:constants:choices+} (which implies $ \Upsilon^2 R\ge 4L/M$). The above series of inequalities imply that
 \begin{align*}
2\|x^0\|_{\Upsilon^2 R}^2\le \left(3-\frac{3}{8M}\right)\|x^0\|^2_{\Upsilon^2  R}\le  3M (\nabla f(0))^T  L^{-1} \nabla f(0). 
 \end{align*}
Therefore overall we have 
\begin{align}
&(x^{0})^T \widetilde{Z} x^0\le 3 (x^{0})^T \Upsilon^2  R x^0 \le 5 M (\nabla f(0))^T  L^{-1} \nabla f(0).
\end{align}
Finally, by observing $\frac{1}{M^2} d^T_0 d_0 =\|\nabla f(0)\|^2$, the desired result is obtained.
\QED

\subsection{Proof of Theorem \ref{thm:final:bound}}
\noindent{\bf Proof.} To show the first part, we consider the optimality condition \eqref{eq:opt:new:lambda}, and multiply both sides of it  by the all one vector, and use the fact that $1^T A^T =0$ to obtain
\begin{align*}
1^T \nabla f(x^r) + 1^TH (x^{r+1}-x^r) =0, \; \forall~r\ge -1.
\end{align*}
Squaring both sides and rearranging terms, we have{
\begin{align*}
&\bigg\|\frac{1}{M}\sum_{i=1}^{M}\nabla f_i(x^{r}_i)\bigg\|^2 \le (x^{r+1}-x^r)^T   H 1 1^T  H (x^{r+1}-x^r)\nonumber\\
%& \le (x^{r+1}-x^r)^T   H    (x^{r+1}-x^r) \lambda_{\max}H ^{1/2} 1 1^TH^{1/2})\nonumber\\
& \le (x^{r+1}-x^r)^T   H   (x^{r+1}-x^r) \times 1^T H 1 \nonumber\\
& \le \|x^{r+1}-x^r\|^2_{H}     \times \bigg(4 \sum_{(i,j)i \sim j} {\sigma^2_{ij}} + \sum_{i=1}^{M}\beta^2_i\bigg), \; \forall~r\ge -1.
\end{align*}}
Combining with \eqref{eq:final:descent}, we obtain, for all $r\ge 0${
\begin{align}\label{eq:key:bound}
&\bigg\|\frac{1}{M}\sum_{i=1}^{M}\nabla f_i(x^{r}_i)\bigg\|^2 \le\|x^{r+1}-x^r\|^2_{H}  \left(4 \sum_{e\in\cE} \sigma^2_{e} + \sum_{i=1}^{M}\beta^2_i\right)\nonumber\\
&\le 8\left(P^r- P^{r+1}\right)   \bigg(4 \sum_{e\in\cE} \sigma^2_{e} + \sum_{i=1}^{M}\beta^2_i\bigg) .
\end{align}}
where in the last inequality we used $2(\Delta+\Upsilon^2)\succeq H$. 
We then bound the consensus error. Lemma \ref{lm:dual:different} implies
\begin{align}
&\|\Sigma F x^{r+1}\|^2 \le \kappa \left(\frac{2}{M^2} \|\Upsilon^{-1} L (x^r-x^{r-1})\|^2+  2  \|w^{r+1}\|^2_{H}\right) \nonumber\\
&\stackrel{\eqref{eq:constants:choices:kappa}}\le 2 \kappa \left(4 \|w^{r+1}\|^2_H + \frac{2}{M^2} \|\Upsilon^{-1} L (x^{r+1}-x^{r})\|^2\right).
\end{align}
Therefore 
\begin{align}
&\|\Sigma F x^{r}\|^2 \le 4 \kappa \left(4 \|w^{r+1}\|^2_H + \frac{2}{M^2} \|\Upsilon^{-1} L (x^{r+1}-x^{r})\|^2\right)+ 2 \|\Sigma F (x^{r+1}-x^r)\|^2.
\end{align}
Combining with \eqref{eq:final:descent}, and using the fact that [cf. \eqref{eq:constants:choices}]
\begin{align}
\Delta+ \Upsilon^2\succeq \frac{8\kappa L \Upsilon^{-2} L}{M^2}, \quad 2\Delta\succeq F \Sigma^2 F
\end{align}
we have
\begin{align}\label{eq:bound:consensus}
&\|\Sigma F x^{r}\|^2  \le 16\kappa \|w^{r+1}\|^2_H + 5\|x^{r+1}-x^r\|^2_{\Delta+\Upsilon^2} \stackrel{\eqref{eq:final:descent}}\le 20 (P^r- P^{r+1}).
\end{align}
Also note that by the definition of $e(T)$ we have
\begin{align}
T \times  e(T)\le \sum_{r=1}^{T}\left(\|\Sigma F x^{r}\|^2+\bigg\|\frac{1}{M}\sum_{i=1}^{M}\nabla f_i(x^{r}_i)\bigg\|^2\right)
\end{align}
Then combining the above with \eqref{eq:final:descent} and \eqref{eq:bound:consensus}, and the fact that the potential function is lower bounded by $\underline{f}$, we obtain the desired result. 

To show the result for xFILTER, multiply both sides of the optimality condition \eqref{eq:opt:gpda+} by the all one vector, and use the fact that $F 1 =0$ to obtain
\begin{align}
1^T \nabla f(x^r) + 1^T \Upsilon^2(x^{r+1}-x^r) = 1^T \Upsilon^2  R\epsilon^{r+1}.
\end{align}
Squaring both sides and rearranging terms we have 
{
	\begin{align*}
	\bigg\|\frac{1}{M}\sum_{i=1}^{M}\nabla f_i(x^{r}_i)\bigg\|^2 &\le 2 (x^{r+1}-x^r)^T   \Upsilon^2 1 1^T  \Upsilon^2 (x^{r+1}-x^r) + 2 (\epsilon^{r+1})^T \Upsilon^2 R 1 1^T \Upsilon^2 R \epsilon^{r+1}\nonumber\\
%	& \le 2(x^{r+1}-x^r)^T   \Upsilon^2    (x^{r+1}-x^r) \lambda_{\max}H ^{1/2} 1 1^TH^{1/2})\nonumber\\
	& \stackrel{\eqref{eq:epsilon:bound}}\le 2 (x^{r+1}-x^r)^T   \Upsilon^2   (x^{r+1}-x^r) \times 1^T \Upsilon^2 1 +  M/(4M) \|x^{r+1}-x^r\|^2_{\Upsilon^2  R} \nonumber\\
	& \le \|x^{r+1}-x^r\|^2_{\Upsilon^2 R}  \times  2 \bigg(1+ \sum_{i=1}^{M}\beta^2_i\bigg), \; \forall~r\ge -1.
	\end{align*}}
where in the last inequality we have used the fact that $\Upsilon^2 R = \Upsilon^2 + F^T \Sigma^2 F \succeq \Upsilon^2$. 

To bound the consensus error, we first use \eqref{eq:epsilon:bound} and obtain 
\begin{align*}
\|\Upsilon^2  R(\epsilon^{r+1}-\epsilon^r)\|^2\le \frac{1}{4M}\|x^{r+1}-x^r\|^2_{\Upsilon^2 R}+ \frac{1}{4M}\|x^{r}-x^{r-1}\|^2_{\Upsilon^2 R}.
\end{align*}
Similarly as the first part, we use Lemma \ref{lm:dual:different} and obtain {
\begin{align}\label{eq:temp2}
&\|\Sigma F x^{r+1}\|^2 \nonumber\\
&\le 3\widetilde{\kappa} \left( \|x^{r+1}-x^r\|^2_{\frac{\Upsilon^2 R}{4M}}+   \|w^{r+1}\|^2_{\Upsilon^2} + \|x^r-x^{r-1}\|^2_{\frac{\Upsilon^2 R}{4M}+ \frac{L \Upsilon^{-2} L}{M^2}}\right)\nonumber\\
& \le 2 \|x^{r+1}-x^r\|^2_{\Upsilon^2 R} + 3\tilde{\kappa} \|w^{r+1}\|^2_{\Upsilon^2} + 2\|x^{r}-x^{r-1}\|_{\Upsilon^2 R}^2, \;\;\forall~r\ge 0
\end{align}}
where the last inequality comes from \eqref{eq:constants:choices+}, that
\begin{align}
2 \Upsilon^2  R\succeq 3\widetilde{\kappa} \left(\frac{L \Upsilon^{-2} L}{M^2} +{\Upsilon^2  R}\right). %\succeq  3\widetilde{\kappa} \left(\frac{L Y^{-2} L}{M^2} +\Upsilon^2\right).
\end{align}
By combining \eqref{eq:temp2} and the following inequality 
$$\|\Sigma F x^{r}\|^2 \le 2	\|\Sigma F (x^{r+1}-x^{r})\|^2 + 2\|\Sigma F x^{r+1}\|^2,$$ 
we have{
	\begin{align*}
	\|\Sigma F x^{r}\|^2 &\le 4 \|x^{r+1}-x^r\|^2_{\Upsilon^2  R + F^T \Sigma^2 F} + 6\tilde{\kappa} \|w^{r+1}\|^2_{\Upsilon^2} + 4\|x^{r}-x^{r-1}\|_{\Upsilon^2  R}^2 \\
	&\stackrel{\eqref{eq:final:descent+}}\le 64(\widetilde{P}^r-\widetilde{P}^{r+1}) + 64(\widetilde{P}^{r-1}-\widetilde{P}^{r}), \quad \forall~r\ge 1\\
	\|\Sigma F x^{0}\|^2 & \le 64(\widetilde{P}^0-\widetilde{P}^{1})  + 4\|x^0\|^2_{\Upsilon^2  R}.
	\end{align*}}
So overall we have that{
\begin{align}
&\sum_{r=0}^{T_r}\left(\bigg\|\frac{1}{M}\sum_{i=1}^{M}\nabla f_i(x^{r}_i)\bigg\|^2  + 	\|\Sigma F x^{r}\|^2\right)\nonumber\\
&\le 64\left( 1 + \sum_{i=1}^{M}\beta^2_i + 1\right) \sum_{r=1}^{T_r}((\widetilde{P}^r-\widetilde{P}^{r+1})  + (\widetilde{P}^{r-1}-\widetilde{P}^{r}))  \nonumber\\
&\quad +  64(\widetilde{P}^0-\widetilde{P}^{1})  + 4\|x^0\|^2_{\Upsilon^2 R}\nonumber\\
&\le 128\left( 1 + \sum_{i=1}^{M}\beta^2_i + 2\right) (\widetilde{P}^0 - \underline{f}) + 4\|x^0\|^2_{\Upsilon^2 R}.
\end{align}}
where the last inequality utilizes the descent property of $\widetilde{P}^0$ in Lemma \ref{lm:Potential+}.
Note that from \eqref{eq:al0}, \eqref{eq:P0+} and use $\widetilde{c}=8\widetilde{\kappa}$ in \eqref{eq:constants:choices:kappa+}, we obtain  
	\begin{align}
	\widetilde{P}^0 \ge f(x^0) +   {\widetilde{\kappa}}\|x^0\|^2_{\Upsilon^2  R}.
	\end{align}
	Therefore From \eqref{eq:P0:bound+} and Lemma \ref{lm:bounded} we have that{\small
		\begin{align*}
		4\|x^0\|^2_{\Upsilon^2 R}\ &\le\frac{4\left( \widetilde{P}^0 -f(x^0)\right)}{\widetilde{\kappa}} \\
		&\stackrel{\eqref{eq:P0:bound+}}\le \frac{4 \left( f(x^0)+  \frac{5}{M} d_0^{\top} L^{-1} d_0 - \underline{f} \right)}{\widetilde{\kappa}}:= \frac{4\widetilde{C}_1}{\widetilde{\kappa}}.
		\end{align*}}
	Combining the above two relations leads to {\small
		\begin{align*}
		&\frac{1}{T_r}\sum_{r=0}^{T_r}\left(\big\|\frac{\sum_{i=1}^{M}\nabla f_i(x^{r}_i)}{M}\big\|^2  + 	\|\Sigma F x^{r}\|^2\right)\\
		&\le 128\left( \sum_{i=1}^{M}\beta^2_i + 3 + \frac{1}{32\widetilde{\kappa}}\right) \widetilde{C}_1/T_r.
		\end{align*}}
	This completes the proof. 
\QED

{\small
\bibliographystyle{IEEEtran}
\bibliography{ref_2019,biblio,ref_pprox_pda,refs,ref_zero1,ref_zero2,Surbib.bib,nsf_refs_2014_OS,distributed_opt,biblio,ref-Haoran,ref}}

\begin{thebibliography}{10}
\providecommand{\url}[1]{#1}
\csname url@samestyle\endcsname
\providecommand{\newblock}{\relax}
\providecommand{\bibinfo}[2]{#2}
\providecommand{\BIBentrySTDinterwordspacing}{\spaceskip=0pt\relax}
\providecommand{\BIBentryALTinterwordstretchfactor}{4}
\providecommand{\BIBentryALTinterwordspacing}{\spaceskip=\fontdimen2\font plus
\BIBentryALTinterwordstretchfactor\fontdimen3\font minus
  \fontdimen4\font\relax}
\providecommand{\BIBforeignlanguage}[2]{{%
\expandafter\ifx\csname l@#1\endcsname\relax
\typeout{** WARNING: IEEEtran.bst: No hyphenation pattern has been}%
\typeout{** loaded for the language `#1'. Using the pattern for}%
\typeout{** the default language instead.}%
\else
\language=\csname l@#1\endcsname
\fi
#2}}
\providecommand{\BIBdecl}{\relax}
\BIBdecl

\bibitem{Lian17decentralized}
X.~Lian, C.~Zhang, H.~Zhang, C.-J. Hsieh, W.~Zhang, and J.~Liu, ``Can
  decentralized algorithms outperform centralized algorithms? a case study for
  decentralized parallel stochastic gradient descent,'' in \emph{Advances in
  Neural Information Processing Systems}, 2017.

\bibitem{Forero11}
P.~A. Forero, A.~Cano, and G.~B. Giannakis, ``Distributed clustering using
  wireless sensor networks,'' \emph{IEEE Journal of Selected Topics in Signal
  Processing}, vol.~5, no.~4, pp. 707--724, Aug 2011.

\bibitem{Wai15_icassp}
T.-H.~C. H.-T.~Wai and A.~Scaglione, ``A consensus-based decentralized
  algorithm for non-convex optimization with application to dictionary
  learning,'' in \emph{the Proceedings of IEEE International Conference on
  Acoustics, Speech and Signal Processing}, 2015.

\bibitem{Nedic09subgradient}
A.~Nedi{\'c} and A.~Ozdaglar, ``Distributed subgradient methods for multi-agent
  optimization,'' \emph{IEEE Transactions on Automatic Control}, vol.~54,
  no.~1, pp. 48--61, 2009.

\bibitem{nedic2015distributed}
A.~Nedic and A.~Olshevsky, ``Distributed optimization over time-varying
  directed graphs,'' \emph{IEEE Transactions on Automatic Control}, vol.~60,
  no.~3, pp. 601--615, 2015.

\bibitem{shi2014extra}
W.~Shi, Q.~Ling, G.~Wu, and W.~Yin, ``Extra: An exact first-order algorithm for
  decentralized consensus optimization,'' \emph{SIAM Journal on Optimization},
  vol.~25, no.~2, pp. 944--966, 2014.

\bibitem{jakovetic2015linear}
D.~Jakoveti{\'c}, J.~M. Moura, and J.~Xavier, ``Linear convergence rate of a
  class of distributed augmented lagrangian algorithms,'' \emph{IEEE
  Transactions on Automatic Control}, vol.~60, no.~4, pp. 922--936, 2015.

\bibitem{BoydADMMsurvey2011}
S.~Boyd, N.~Parikh, E.~Chu, B.~Peleato, and J.~Eckstein, ``Distributed
  optimization and statistical learning via the alternating direction method of
  multipliers,'' \emph{Foundations and Trends in Machine Learning}, vol.~3,
  no.~1, pp. 1--122, 2011.

\bibitem{Schizas09}
I.~Schizas, G.~Mateos, and G.~Giannakis, ``Distributed {LMS} for
  consensus-based in-network adaptive processing,,'' \emph{IEEE Transactions on
  Signal Processing}, vol.~57, no.~6, pp. 2365 -- 2382, 2009.

\bibitem{bianchi2013convergence}
P.~Bianchi and J.~Jakubowicz, ``Convergence of a multi-agent projected
  stochastic gradient algorithm for non-convex optimization,'' \emph{IEEE
  Transactions on Automatic Control}, vol.~58, no.~2, pp. 391--405, 2013.

\bibitem{Zhu-Martinez2}
M.~Zhu and S.~Mart\'{\i}nez, ``An approximate dual subgradient algorithm for
  distributed non-convex constrained optimization,'' \emph{IEEE Transactions on
  Automatic Control}, vol.~58, no.~6, pp. 1534--1539, June 2013.

\bibitem{hong14nonconvex_admm}
M.~Hong, Z.-Q. Luo, and M.~Razaviyayn, ``Convergence analysis of alternating
  direction method of multipliers for a family of nonconvex problems,''
  \emph{SIAM Journal On Optimization}, vol.~26, no.~1, pp. 337--364, 2016.

\bibitem{Lorenzo16}
P.~Di~Lorenzo and G.~Scutari, ``Next: In-network nonconvex optimization,''
  \emph{IEEE Transactions on Signal and Information Processing over Networks},
  vol.~2, no.~2, pp. 120--136, 2016.

\bibitem{Hajinezhad17inexact}
D.~Hajinezhad and M.~Hong, ``Perturbed proximal primal dual algorithm for
  nonconvex nonsmooth optimization,'' \emph{Mathematical Programming}, vol.
  176, no. 1-2, pp. 207--245, July 2019.

\bibitem{hong17icml}
M.~Hong, D.~Hajinezhad, and M.-M. Zhao, ``{Prox-PDA}: The proximal primal-dual
  algorithm for fast distributed nonconvex optimization and learning over
  networks,'' in \emph{the Proceedings of the 34th International Conference on
  Machine Learning (ICML)}, 2017.

\bibitem{Hajinezhad17zeroth}
D.~Hajinezhad, M.~Hong, and A.~Garcia, ``Zone: Zeroth order nonconvex
  multi-agent optimization over networks,'' \emph{IEEE Transactions on
  Automatic Control}, 2019.

\bibitem{Daneshmand_et_al_Asilomar16}
A.~Daneshmand, G.~Scutari, and F.~Facchinei, ``Distributed dictionary
  learning,'' in \emph{Proceedings of the Asilomar Conference on Signals,
  Systems, and Computers}, Nov. 6--9, 2016.

\bibitem{Daneshmand_et_al_ICASSP15}
A.~Daneshmand, Y.~Sun, G.~Scutari, and F.~Facchinei, ``Distributed dictionary
  learning over networks,'' in \emph{Proceedings of the IEEE International
  Conference on Acoustics, Speech and Signal Processing}, March 5-9 2017.

\bibitem{Zeng19distributedGD}
J.~Zeng and W.~Yin, ``On nonconvex decentralized gradient descent,'' \emph{IEEE
  Transactions on Signal Processing}, vol.~66, no.~11, pp. 2834--2848, June
  2018.

\bibitem{Jiang2017}
Z.~Jiang, A.~Balu, C.~Hegde, and S.~Sarkar, ``Collaborative deep learning in
  fixed topology networks,'' in \emph{Advances in Neural Information Processing
  Systems}, 2017.

\bibitem{vlaski2019distributed}
S.~Vlaski and A.~H. Sayed, ``Distributed learning in non-convex
  environments--part i: Agreement at a linear rate,'' \emph{arXiv preprint
  arXiv:1907.01848}, 2019.

\bibitem{vlaski2019distributed2}
------, ``Distributed learning in non-convex environments--part ii: Polynomial
  escape from saddle-points,'' \emph{arXiv preprint arXiv:1907.01849}, 2019.

\bibitem{swenson2019annealing}
B.~Swenson, S.~Kar, H.~V. Poor, and J.~Moura, ``Annealing for distributed
  global optimization,'' \emph{arXiv preprint arXiv:1903.07258}, 2019.

\bibitem{hong172ndorder}
M.~Hong, J.~D. Lee, and M.~Razaviyayn, ``Gradient primal-dual algorithm
  converges to second-order stationary solutions for nonconvex distributed
  optimization,'' in \emph{the Proceedings of the 35th International Conference
  on Machine Learning (ICML)}, 2018.

\bibitem{daneshmand2018second}
A.~Daneshmand, G.~Scutari, and V.~Kungurtsev, ``Second-order guarantees of
  gradient algorithms over networks,'' in \emph{Proceedings of the 56th Annual
  Allerton Conference on Communication, Control, and Computing (Allerton)},
  2018.

\bibitem{swenson2019distributed}
B.~Swenson, R.~Murray, H.~V. Poor, and S.~Kar, ``Distributed gradient descent:
  Nonconvergence to saddle points and the stable-manifold theorem,''
  \emph{arXiv preprint arXiv:1908.02747}, 2019.

\bibitem{duenner2018distributed}
C.~Duenner, A.~Lucchi, M.~Gargiani, A.~Bian, T.~Hofmann, and M.~Jaggi, ``A
  distributed second-order algorithm you can trust,'' in \emph{Proceedings of
  the International Conference on Machine Learning (ICML)}, 2018.

\bibitem{fang2018distributed}
C.-H. Fang, S.~B. Kylasa, F.~Roosta-Khorasani, M.~W. Mahoney, and A.~Grama,
  ``Distributed second-order convex optimization,'' \emph{arXiv preprint
  arXiv:1807.07132}, 2018.

\bibitem{nesterov05}
Y.~Nesterov, ``Smooth minimization of nonsmooth functions,'' \emph{Mathematical
  Programming}, vol. 103, pp. 127--152, 2005.

\bibitem{nesterov83}
------, ``A method of solving a convex programming problem with convergence
  rate $o(1/k^2)$,'' \emph{Soviet Mathematics Doklady}, vol.~27, pp. 372--376,
  1983.

\bibitem{nemirovsky83}
A.~Nemirovsky and D.~Yudin, ``Problem complexity and method efficiency in
  optimization,'' in \emph{Interscience Series in Discrete Mathematics}.\hskip
  1em plus 0.5em minus 0.4em\relax Wiley, 1983.

\bibitem{nesterov04}
Y.~Nesterov, \emph{Introductory lectures on convex optimization: A basic
  course}.\hskip 1em plus 0.5em minus 0.4em\relax Springer, 2004.

\bibitem{Beck:2009:FIS:1658360.1658364}
A.~Beck and M.~Teboulle, ``A fast iterative shrinkage-thresholding algorithm
  for linear inverse problems,'' \emph{SIAM Journal on Imgaging Science},
  vol.~2, no.~1, pp. 183 -- 202, 2009.

\bibitem{Ouyang15}
Y.~Ouyang, Y.~Chen, G.~Lan, and J.~E.~Pasiliao, ``An accelerated linearized
  alternating direction method of multipliers,'' \emph{SIAM Journal on Imaging
  Sciences}, vol.~8, no.~1, pp. 644--681, 2015.

\bibitem{tseng08acc}
P.~Tseng, ``On accelerated proximal gradient methods for convex-concave
  optimization,'' 2008, preprint.

\bibitem{Jakovetic14}
D.~Jakovetic, J.~Xavier, and J.~M.~F. Moura, ``Fast distributed gradient
  methods,'' \emph{IEEE Transactions on Automatic Control}, vol.~59, no.~5, pp.
  1131--1146, May 2014.

\bibitem{scaman2017optimal}
K.~Scaman, F.~Bach, S.~Bubeck, Y.~Lee, and L.~Massouli{\'e}, ``Optimal
  algorithms for smooth and strongly convex distributed optimization in
  networks,'' \emph{arXiv preprint arXiv:1702.08704}, 2017.

\bibitem{uribe2017optimal}
C.~Uribe, S.~Lee, A.~Gasnikov, and A.~Nedi{\'c}, ``Optimal algorithms for
  distributed optimization,'' \emph{arXiv preprint arXiv:1712.00232}, 2017.

\bibitem{scaman2018optimal}
K.~Scaman, F.~Bach, S.~Bubeck, L.~Massouli{\'e}, and Y.~T. Lee, ``Optimal
  algorithms for non-smooth distributed optimization in networks,'' in
  \emph{Advances in Neural Information Processing Systems}, 2018, pp.
  2740--2749.

\bibitem{cartis2010complexity}
C.~Cartis, N.~Gould, and P.~Toint, ``On the complexity of steepest descent,
  newton's and regularized newton's methods for nonconvex unconstrained
  optimization problems,'' \emph{{SIAM} journal on optimization}, vol.~20,
  no.~6, pp. 2833--2852, 2010.

\bibitem{carmon2017lower}
Y.~Carmon, J.~C. Duchi, O.~Hinder, and A.~Sidford, ``Lower bounds for finding
  stationary points i,'' \emph{Mathematical Programming}, Jun 2019.

\bibitem{tian2018asy}
Y.~Tian, Y.~Sun, B.~Du, and G.~Scutari, ``Asy-sonata: Achieving geometric
  convergence for distributed asynchronous optimization,'' \emph{arXiv preprint
  arXiv:1803.10359}, 2018.

\bibitem{Daneshmand18}
A.~Daneshmand, Y.~Sun, and G.~Scutari, ``Convergence rate of distributed convex
  and nonconvex optimization methods with gradient tracking,'' 2018, purdue
  University, Tech. Rep.

\bibitem{fu16}
X.~Fu, K.~Huang, N.~Sidiropolous, A.~M.-S. So, and M.~Hong, ``Scalable and
  optimal generalized canonical correlation analysis via alternating
  optimization,'' 2016, submitted to NIPS 2016.

\bibitem{Chung97}
F.~R.~K. Chung, \emph{Spectral Graph Theory}.\hskip 1em plus 0.5em minus
  0.4em\relax The American Mathematical Society, 1997.

\bibitem{Butler2016}
S.~Butler, \emph{Algebraic aspects of the normalized Laplacian}.\hskip 1em plus
  0.5em minus 0.4em\relax Cham: Springer International Publishing, 2016, pp.
  295--315.

\bibitem{Duchi12}
J.~C. Duchi, A.~Agarwal, and M.~J. Wainwright, ``Dual averaging for distributed
  optimization: Convergence analysis and network scaling,'' \emph{IEEE
  Transactions on Automatic Control}, vol.~57, no.~3, pp. 592--606, March 2012.

\bibitem{UZAWA58}
H.~Uzawa, ``Iterative methods in concave programming,'' in \emph{Studies in
  Linear and Nonlinear Programming}.\hskip 1em plus 0.5em minus 0.4em\relax
  Stanford University Press, 1958, p. 154–165.

\bibitem{Nedic2009saddle}
A.~Nedi{\'c} and A.~Ozdaglar, ``Subgradient methods for saddle-point
  problems,'' \emph{Journal of optimization theory and applications}, vol. 142,
  no.~1, pp. 205--228, 2009.

\bibitem{rockafellar1976augmented}
R.~T. Rockafellar, ``Augmented lagrangians and applications of the proximal
  point algorithm in convex programming,'' \emph{Mathematics of operations
  research}, vol.~1, no.~2, pp. 97--116, 1976.

\bibitem{wright_proximal}
S.~J.Wright, ``Implementing proximal point methods for linear programming,''
  \emph{Journal of Optimization Theory and Applications}, vol.~65, no.~3, pp.
  531--554, Jun 1990.

\bibitem{Tian14}
D.~Tian, H.~Mansour, A.~Knyazev, and A.~Vetro, ``Chebyshev and conjugate
  gradient filters for graph image denoising,'' in \emph{Proceedings of the
  IEEE International Conference on Multimedia and Expo Workshops (ICMEW)},
  2014.

\bibitem{Gadde13}
A.~Gadde, S.~K. Narang, and A.~Ortega, ``Bilateral filter: Graph spectral
  interpretation and extensions,'' in \emph{Proceedings of the IEEE
  International Conference on Image Processing}.

\bibitem{Tsynkov07}
V.~S. Ryaben'kii and S.~V. Tsynkov, \emph{A Theoretical Introduction to
  Numerical Analysis}.\hskip 1em plus 0.5em minus 0.4em\relax CRC Press, 2007.

\bibitem{Samarskij89}
A.~A. Samarskij and E.~S. Nikolaev, \emph{Numerical Methods for Grid Equations
  Volume II Iterative Methods}.\hskip 1em plus 0.5em minus 0.4em\relax
  Springer, 1989.

\bibitem{boyd2006randomized}
S.~Boyd, A.~Ghosh, B.~Prabhakar, and D.~Shah, ``Randomized gossip algorithms,''
  \emph{IEEE Transactions on Information Theory}, vol.~52, no.~6, pp.
  2508--2530, 2006.

\bibitem{antoniadis2011penalized}
A.~Antoniadis, I.~Gijbels, and M.~Nikolova, ``Penalized likelihood regression
  for generalized linear models with non-quadratic penalties,'' \emph{Annals of
  the Institute of Statistical Mathematics}, vol.~63, no.~3, pp. 585--615,
  2011.

\bibitem{tatarenko2017non}
T.~Tatarenko and B.~Touri, ``Non-convex distributed optimization,'' \emph{IEEE
  Transactions on Automatic Control}, vol.~62, no.~8, pp. 3744--3757, 2017.

\bibitem{nedic2009distributed}
A.~Nedi{\'c}, A.~Olshevsky, A.~Ozdaglar, and J.~N. Tsitsiklis, ``On distributed
  averaging algorithms and quantization effects,'' \emph{IEEE Transactions on
  Automatic Control}, vol.~54, no.~11, pp. 2506--2517, 2009.

\bibitem{chang14distributed}
T.-H. Chang, M.~Hong, and X.~Wang, ``Multi-agent distributed optimization via
  inexact consensus {ADMM},'' \emph{IEEE Transactions on Signal Processing},
  vol.~63, no.~2, pp. 482--497, Jan 2015.

\end{thebibliography}
\newpage
 
\end{document}